\documentclass[leqno]{amsart}
  \usepackage{amsmath}
  \usepackage{amssymb,amsxtra,amsthm,amscd}
  \usepackage{blaom}
 \theoremstyle{plain}% That is, in italics
 \newtheorem*{equivalenceProblem}{Equivalence Problem}
 \newtheorem*{obstructionProblem}{Obstruction Problem}
 \renewcommand{\equiv}{:=}
 \newcommand{\ii}{\mathrm{I}\!\mathrm{I}}
  
  \newcommand{\gl}{\mathfrak{gl}}
  
  \newcommand{\levi}{\nabla^{\mathrm{L}\text{-}\mathrm{C}}}
  \DeclareMathOperator{\cocurvature}{cocurv}                                          
  \newcommand{\sym}{\mathrm{sym}}
  \newcommand{\alt}{\mathrm{alt}}

  \DeclareMathOperator{\coRicci}{coRicci}
  \DeclareMathOperator{\sskew}{skew}
  \newcommand{\transpose}{\mathrm{t}}
  \renewcommand{\paragraph}[1]{{\indent\em #1.}}  
\begin{document}
\title[Cartan's method of equivalence]%
{Lie algebroids and\\ Cartan's method of equivalence}
\author{Anthony D.~Blaom}%\\ \mbox{} \\ {\sffamily DRAFT}}%
\date{\today}%
\address{22 Ridge Road, Waiheke Island, New Zealand. {\tt anthony.blaom@gmail.com}}%
\thispagestyle{empty}
\begin{abstract}
  \'Elie Cartan's general equivalence problem is recast in the
  language of Lie algebroids. The resulting formalism, being
  coordinate and model-free, allows for a full geometric
  interpretation of Cartan's method of equivalence via reduction and
  prolongation. We show how to construct certain normal forms (Cartan
  algebroids) for objects of finite-type, and are able to interpret
  these directly as `infinitesimal symmetries deformed by curvature.'

  Details are developed for transitive structures but rudiments of the
  theory include {\em intransitive} structures (intransitive symmetry
  deformations).  Detailed illustrations include subriemannian contact
  structures and conformal geometry.
%
%   Cartan's method of equivalence constructs the local invariants of
%   geometric structures, realizing them as obstructions to symmetry;
%   these invariants generalize the familiar curvature invariant of a
%   Riemannian structure.  Here Cartan's method is recast in the
%   language of Lie algebroids.  The resulting formalism is fully
%   invariant, for it is not only coordinate-free but {\em
%     model\/}-free.
  
%   Details are developed for transitive finite-type geometric
%   structures but rudiments of the theory include {\em intransitive}
%   structures (intransitive symmetry deformations).  Detailed
%   illustrations include subriemannian contact structures and conformal
%   geometry.
\end{abstract}
\maketitle
%% Curie quote:
\vspace{-\baselineskip}
{\small
\begin{center}
   {\em c'est la dissym\'etrie qui cr\'ee le ph\'enom\`ene}%\newline
   %\mbox{~}
   \qquad--- Pierre Curie, 1894.
\end{center}}
\vspace{\baselineskip}%
\tableofcontents
\newpage
%%%%%%%%%%%%%%%%%%%%%%%%%%%%%%%%% Body starts here %%%%%%%%%%%%%%%%%%%%%%%
\section{A new setting for Cartan's method of equivalence}\lab{toodle}
This paper has its origins in attempts to understand, in as
invariant a language as possible, \'Elie Cartan's assertion that
finite-type geometric structures are `symmetries deformed by
curvature.' Having identified in \cite{Blaom_05} a model of geometric
structures well-suited to this viewpoint --- namely Lie algebroids
equipped with a suitably compatible linear connection, called {\df
  Cartan algebroids} --- we turn to the problem of realizing this
model in practice. To do so requires us to revisit Cartan's method of
equivalence, and even to reformulate the basic `problem of
equivalence' that it addresses.  The advantages of this reformulation
are both theoretical and computational, as we shall demonstrate.

The introduction to this paper is in two parts. In this first section we recast Cartan's equivalence problem in the language of Lie algebroids and formulate specific objectives for the paper. We describe the basic geometric objects with which the paper will be concerned: certain infinitesimal models of geometric structures, and Cartan algebroids, which amount to normal forms thereof.

In Sect.\,\ref{normalform} we outline those elements of Cartan's method which allow us, in the Lie algebroid setting, to associate, with any given object of finite-type, an intrinsically defined Cartan algebroid. Limitations of the method and an outline of the remainder of the paper will also be given there.

% Throughout the paper the reader may like to keep the following guiding goal in mind: to understand, in a nice geometric way, the invariants obstructing the existence of symmetry, in a finite-type geometric structure, such as a Riemannian or conformal structure.

\subsection{The equivalence problem}
Cartan's method of equivalence is a procedure for determining when two objects are equivalent under a change of coordinates. The method applies to an astonishing variety of objects: differential equations, polynomials and variational problems; tensor structures on smooth manifolds, such as Riemannian, conformal, symplectic, complex and almost K\"ahler structures; differential operators and associated smooth manifold structures, such as affine and projective structures; and so on. For an introduction to the method and a survey of applications see \cite{Olver_95,Gardner_89,Ivey_Landsberg_03,Bryant_etal_91,Montgomery_02}.  

What makes the method so general is that the equivalence problem is
formulated in terms of certain secondary data of universal form, this
data encoding essential information about the objects, rather than in
terms of the objects themselves. In Cartan's original approach the
secondary data is a collection of one-forms (a coframe), defined
pointwise up to extra `group parameters.' In subsequent
reformulations, the secondary data is a $G$-structure (see, e.g.,
\cite{Sternberg_64}) or an exterior differential system (see, e.g.,
\cite{Bryant_etal_91}). While in practice the construction of the
secondary data can seem ad hoc, it can often be given the following
interpretation: it is a geometrization of the PDE's one would write
down to find the symmetries (self-equivalences) of the given object,
whether or not these symmetries actually exist.

In our Lie algebroid reformulation of Cartan's method, the secondary
data --- here called an {\df infinitesimal geometric structure} ---
can generally be understood as a geometrization of the PDE's one
writes down to find the {\em infinitesimal} symmetries of the given
object. This geometrization is coordinate-free and does not depend on
any choice of group, group quotient, or other fixed model. We motivate
the formal definition with the following example.

\subsection{The 1-symmetries of a Riemannian metric}\lab{firsty}
Consider the problem of finding the infinitesimal symmetries of a Riemannian metric $\sigma$,  i.e., its Killing fields, defined on a smooth $n$-dimensional manifold $M$. Call the 1-jet of a vector field $ V$ on $M$, evaluated at some point $m\in M$, a {\df 1-symmetry} of $\sigma$, whenever $\sigma$
has, at $ m\in M$, vanishing Lie derivative along $ V$. Then
the collection of all 1-symmetries is a subbundle
${\mathfrak g}\subset J^1 (TM)$ and $ V$ is a Killing field (i.e., genuine 
symmetry) if and only if its first-order prolongation
$J^1 V$ is a section of ${\mathfrak g}$. The sections of $J^1(TM)$ of
the form $J^1 V$ for some $V$ are called {\df holonomic}. Whence:
\begin{conditions}
\item\lab{intro1} {\em The Killing fields of $\sigma$ are in
    one-to-one correspondence with the holonomic sections of
    $\mathfrak g$.}
\end{conditions}
Moreover, as we show later, $\sigma$ can be recovered from ${\mathfrak g}$ up to a constant factor, so that little is lost by restricting attention to $\mathfrak g$.

The important observation to make here is that ${\mathfrak g}$ is a Lie algebroid over $M$, as is the tangent bundle $TM$, and its first
jet $J^1(TM)$. In fact, adopting the natural language associated with
these objects, we have:
\begin{conditions}
\item{\em The bundle ${\mathfrak g}\subset J^1 (TM)$ of 1-symmetries 
of a Riemannian metric $\sigma$ is the isotropy subalgebroid of 
$\sigma$ under the representation of $J^1 (TM)$ on $\symmetric^2 (TM)$ 
determined by the adjoint representation of $J^1 (TM)$ on $ TM$.}
\end{conditions}
The terms `isotropy' and `adjoint representation' are natural generalizations to Lie algebroids of familiar Lie algebra notions. The adjoint representation of a Lie algebroid is described in in Sect.\,\ref{section2}. Isotropy subalgebroids are defined in Sect.\,\ref{sectionQ}.

A {\df Lie algebroid} over a smooth manifold $M$ is vector bundle
${\mathfrak g}$ on $M$, together with a Lie bracket on its space of
sections, and a vector bundle morphism $\# \colon {\mathfrak g}
\rightarrow TM$ called the {\df anchor}, satisfying certain conditions
making Lie algebroids generalizations of both tangent bundles
($\#\colon TM \rightarrow TM$ the identity) and Lie algebras ($M$ a
single point).  The bundle of $k$-jets of sections of any Lie
algebroid is itself a Lie algebroid. 

Lie algebroids over $M$ also generalize the infinitesimal actions of
Lie algebras on $M$ (see \ref{cart} below), in which case the image of
the anchor $\# \colon {\mathfrak g} \rightarrow TM$ is the
distribution tangent to the foliation of $M$ by orbits. The image of
the anchor of an arbitrary Lie algebroid is always tangent to some
foliation, the leaves of which are accordingly called {\df orbits}; a
Lie algebroid with surjective anchor is called {\df transitive}.

\subsection{Infinitesimal geometric structures and their 
symmetries}\lab{structure}%\mbox{} \newline%
Infinitesimal geometric structures, as defined below, generalize the vector bundle of $1$-symmetries of a Riemannian metric: 
\begin{definition}
Let ${\mathfrak t}$ be any Lie algebroid over $M$ (the tangent bundle
$TM$ in the simplest case).  Then an {\df infinitesimal geometric structure} on
${\mathfrak t}$ is any subalgebroid ${\mathfrak g}\subset
J^1{\mathfrak t}$.    
\end{definition}
In practice the infinitesimal geometric structure associated with a given object --- whether it be a tensor, an operator, or whatever --- is constructed as the isotropy subalgebroid of an appropriate Lie algebroid representation. Further examples, both transitive and intransitive, are given in Sect.\,\ref{sectionQ}.  Infinitesimal geometric structures are this paper's principal objects of investigation.

The {\df structure kernel} of ${\mathfrak g}$ is the kernel
${\mathfrak h} \subset T^*\!M \otimes TM $ of the restriction
$a\colon{\mathfrak g}\rightarrow{\mathfrak t}$ of the projection $
J^1{\mathfrak t}\rightarrow{\mathfrak t}$, this morphism being the
anchor of $\mathfrak g$ when $\mathfrak t=TM$.  This is an analogue of
the structure group $G$ of a $G$-structure. The structure kernel of
$\mathfrak g$ is a subalgebroid whenever it has constant rank.

The {\df image} of ${\mathfrak g}$ is the image of $ a$.  If
$a\colon{\mathfrak g}\rightarrow{\mathfrak t}$ is surjective, we call
${\mathfrak g}$ a {\df surjective} infinitesimal geometric structure.
If $\mathfrak t=TM$, then `surjective' is synonymous with
`transitive.'  For example, the bundle ${\mathfrak g}\subset J^1(TM)$
of 1-symmetries of a Riemannian metric is surjective.  We will see
that every Poisson structure has an associated infinitesimal geometric
structure ${\mathfrak g} \subset J^1(T^*\!M)$ which {\em is}
surjective but generally {\em not} transitive. A simple example of an
infinitesimal geometric structure failing to be surjective or
transitive is the {\em joint} isotropy subalgebroid of a Riemannian
metric and a vector field $V$ with non-degenerate energy
$\frac{1}{2}\|V\|^2$. 

Structures sometimes viewed as transitive are actually intransitive,
when the notion is invariantly formulated. For example, almost complex
structures are generically intransitive structures.

Associated with any $G$-structure on $M$ is a corresponding infinitesimal geometric structure on $TM$, but this structure is {\em always} surjective (equivalently, transitive).

Here now, in Lie algebroid language, is an analogue of Cartan's problem of equivalence: 
\begin{equivalenceProblem}
  Given smooth manifolds $M_1$ and $M_2$ and infinitesimal geometric
  structures ${\mathfrak g}_1 \subset J^1 (TM_1)$ and ${\mathfrak g}_2
  \subset J^1 (TM_2)$, decide whether there exists a diffeomorphism
  $\phi \colon M_1 \rightarrow M_2$, with associated tangent map $T
  \phi \colon TM_1 \rightarrow TM_2$, such that the corresponding Lie
  algebroid isomorphism $J^1(T \phi):J^1 (TM_1)\rightarrow J^1 (TM_2)$
  maps ${\mathfrak g}_1$ isomorphically onto ${\mathfrak g}_2$.
\end{equivalenceProblem}
\begin{remark}
  If we want to formulate a more general equivalence problem,
  replacing $TM_1$ and $TM_2$ with general Lie algebroids ${\mathfrak
    t}_1$ and ${\mathfrak t}_2$, then we are forced to allow for {\em
    arbitrary} Lie algebroid morphisms ${\mathfrak t}_1 \rightarrow
  {\mathfrak t}_2$ instead of the tangent maps $T\phi \colon TM_1
  \rightarrow TM_2$, or we must restrict the class of infinitesimal
  geometric structures considered so that morphisms of `coordinate
  change type' can be defined. However, for the restricted purposes
  of the the present paper, further elaboration on this point will be
  unnecessary.
\end{remark}

\subsection{Cartan's method}
Having formulated the equivalence problem in terms of appropriate
secondary data, Cartan's method attempts to put the data into an
appropriate normal form. In the original one-form approach, for
example, this normal form is a coframe, on a possibly larger space,
but with group parameters eliminated. See, e.g., \cite{Gardner_89} or
\cite{Olver_95}.

The normalizing algorithm involves two fundamental operations, known as {\df reduction} and {\df prolongation}. If the secondary data is regarded as a system of PDE's, then reduction amounts to identifying coordinate changes that decouple certain of the equations from the others, these latter being rendered redundant. Prolongation means increasing the number of variables by introducing derivatives of the independent variables as new independent variables; new equations are added to account for the equality of mixed partial derivatives.

When the normalizing procedure succeeds, the normal form obtained delivers certain basic invariants, and certain `derived' invariants can be constructed from these. One may test for the equivalence of two objects by comparing the values taken on by the invariants. (Specifically, one compares the associated `classifying manifold'; see, e.g., \cite[Chapter 8]{Olver_95}). Objects for which the normalizing procedure succeeds have, at most, finite-dimensional symmetry, and so are called {\df finite-type} objects.

The basic invariants of a finite-type object may be understood as
obstructions to the existence of a maximal set of symmetries
(self-equivalences); they are embodied in the `curvature' of the
associated normal form, when this normal form is understood as a
`symmetry deformed by curvature.' In common practice, this
interpretation is often obscured, however.

For objects of {\df infinite-type}, the normalizing procedure fails to terminate and an altogether different criterion for equivalence must be applied. Infinite-type objects can be identified by applying Cartan's `involutivity' test and are described extensively in \cite{Bryant_etal_91,Ivey_Landsberg_03}. They will not be studied here.

\subsection{The symmetries of infinitesimal geometric structures}
Analogous to the Killing fields of a Riemannian structure are the {\df symmetries}\footnote{The notions of structure kernel, image, surjectivity, and symmetry all make sense for {\em arbitrary} subsets ${\mathfrak g}\subset J^1{\mathfrak t}$.  No further comment will accompany their application in this extended context.} of an infinitesimal geometric structure ${\mathfrak g} \subset J^1{\mathfrak t}$.  These are those sections $V$ of ${\mathfrak t}$ whose prolongations $ J^1 V\subset J^1 {\mathfrak t}$ are sections of ${\mathfrak g}$.  Evidently, the symmetries of $\mathfrak g$ are in one-to-one correspondence with the {\df holonomic} sections of ${\mathfrak g}\subset J^1{\mathfrak t}$ --- those sections that are prolongations of something.  Symmetries are necessarily sections of the image of ${\mathfrak g}$ and are closed under the Lie algebroid bracket. 

We will not be presenting a complete solution to the equivalence problem here. Rather, our main focus is the following: %
\begin{obstructionProblem}
  Given an infinitesimal geometric structure ${\mathfrak g} \subset J^1{\mathfrak t}$, find the obstructions to the existence of symmetries of $\mathfrak g$, in the sense defined above.
\end{obstructionProblem}

We now turn our attention to the normal forms that we shall be constructing out of infinitesimal geometric structures. We begin by describing such a normal form explicitly in the case of Riemannian geometry, and by showing concretely how curvature invariants enter as the obstruction to symmetry.
\subsection{A normal form for Riemannian geometry}\lab{cakeWalk}%
Let ${\mathfrak g} \subset J^1(TM)$ be the bundle of 1-symmetries of a Riemannian metric $\sigma$, as described in \ref{firsty} above. Note that the Levi-Cevita connection $\nabla$ associated with $\sigma$ may be regarded as a splitting of a canonical exact sequence, \begin{equation*}
         0\rightarrow T^*\!M\otimes TM\rightarrow J^1 (TM)\rightarrow 
         TM\rightarrow 0.
\end{equation*}
We have a corresponding exact sequence,
\begin{equation*}
  0\longrightarrow  {\mathfrak h} \longrightarrow {\mathfrak g} \longrightarrow  TM \longrightarrow 0,
\end{equation*}
where ${\mathfrak h}\subset T^*\!M\otimes TM$
denotes the $\mathfrak{o}(n)$-bundle of all $\sigma$-skew-symmetric
tangent space endomorphisms, which is $\nabla$-invariant. 
By the $\nabla$-invariance of $\sigma$, the splitting subbundle of $J^1(TM)$ lies inside $\mathfrak g$
 and we obtain an identification ${\mathfrak g}\cong
TM\oplus{\mathfrak h}$. On $\mathfrak g$, define a linear connection $\nabla^{(1)}$ by
\begin{equation*}
        \nabla^{(1)}_U(V\oplus\phi)=(\nabla_UV+\phi (U))\oplus 
        (\nabla_U\phi+\curvature\nabla\, (U, V)).
\end{equation*}
Then we will eventually prove that:
\begin{conditions}
\item\lab{hutt} {\em     A section of $ X\subset{\mathfrak g}$ is holonomic if and only if it is $\nabla^{(1)}$-parallel.}
\end{conditions}
\noindent In particular, by \eqrefs{firsty}{intro1}, the curvature of
$\nabla^{(1)}$ is the local obstruction to maximal symmetry, i.e., to
the existence, locally, of a set of Killing fields of maximal
dimension.

If $\curvature \nabla^{(1)}=0 $, and we assume for simplicity that $M$ is simply connected, then ${\mathfrak g}\cong {\mathfrak g}_0 \times M$, where ${\mathfrak g}_0$ is the $\frac{1}{2}n(n+1)$-dimensional vector space of $\nabla^{(1)} $-parallel sections, which is naturally isomorphic to the Lie algebra of Killing fields.

With the help of the Bianchi identities for linear connections, one
computes
\begin{equation*}
        \curvature\nabla^{(1)} (U_1, U_2)(V\oplus\phi)=0\oplus 
        (-(\nabla_V\curvature\nabla+\phi\cdot\curvature\nabla)(U_1,U_2)),
\end{equation*}
implying that $\nabla^{(1)}$ is flat if and only if $\curvature\nabla$ is both $\nabla$-invariant and ${\mathfrak h}$-invariant.  Now ${\mathfrak h} $-invariance implies, by purely algebraic arguments, that $\curvature\nabla$ has only a scalar component; $\nabla$-invariance then implies constant scalar curvature.  Whence from \eqref{hutt} one recovers the standard criterion for maximal local homogeneity of a Riemannian manifold.

\subsection{Cartan algebroids: symmetries deformed by curvature}\lab{cart}%
A linear connection $\nabla$ on a Lie algebroid $\mathfrak g$ is a {\df Cartan connection} if it is suitably compatible with the Lie algebroid structure \cite{Blaom_05}. The pair $({\mathfrak g}, \nabla)$ is then a {\df Cartan algebroid}.  The formal definition and basic properties are reviewed in Sect.\,\ref{geomet}.

In the Riemannian example above, the pair $({\mathfrak g},
\nabla^{(1)})$ is a Cartan algebroid, and we saw that ${\mathfrak
  g}\cong {\mathfrak g}_0 \times M$, for some Lie algebra ${\mathfrak
  g}_0$, when $\nabla^{(1)} $ is globally flat. In fact, whenever
${\mathfrak g}_0$ is any Lie algebra, acting smoothly on a manifold
$M$, then the trivial bundle ${\mathfrak g}_0 \times M$ inherits the
structure of a Lie algebroid, called an {\df action} (or {\df
  transformation) algebroid}, and the trivial flat connection is a
Cartan connection. Conversely, any Cartan algebroid with a flat Cartan
connection is locally an action algebroid (Theorem
\ref{dissymmetry}, Sect.\,\ref{geomet}). It is in this sense that
Cartan algebroids are infinitesimal symmetries deformed by
curvature. Moreover, the orbits of the Cartan algebroid may be
regarded as deformations of orbits of a symmetry.

In \cite{Blaom_05} we described how Cartan algebroids may be viewed as
infinitesimal, model-free, and possibly intransitive versions of
classical Cartan geometries, and mentioned other alternative models
contained in the literature. Recently, Crampin \cite{Crampin_09} has
delineated the relationship between transitive Cartan algebroids and
adjoint tractor bundles \cite{Cap_Gover_02}, which like Cartan
algebroids are infinitesimal objects, but unlike them are based on a
transitive model fixed a priori.

One consequence of choosing a model-free approach is worth repeating
here. Generally, `curvature' has referred to the local deviation from
an underlying model --- typically $\mathbb R^n$ or a homogeneous space
$G/H$ --- which is fixed.  In a model-free approach, such as the one
described here, all potential models are created equal, and curvature
merely measures the local deviation from {\em some} maximally
symmetric space.  From this point of view, Euclidean space, hyperbolic
space, and spheres, are all `flat' Riemannian manifolds.  If Cartan's
method determines that a space is flat in the weaker sense, then the
particular flat space one has at hand is simply part of the method's
output.  In model-based approaches model `mutation' (or something
similar) may be needed to detect all possible cases \cite{Sharpe_97}.

In Appendix \ref{global}, we explain how flat Cartan algebroids may be regarded as infinitesimal analogues of Lie pseudogroups of
transformations, and discuss the global analogues of arbitrary Cartan
algebroids, {\df Cartan groupoids}, which may be regarded as `curved' Lie pseudogroups.

\section{An outline of Cartan's method \`a la Lie algebroids}\lab{normalform}%
In this second introductory section we outline the construction of Cartan algebroids from infinitesimal geometric structures.
\subsection{Cartan connections via generators}\lab{neoGenerator}
To solve the obstruction problem formulated in Sect.\,\ref{toodle}, we attempt to reduce it to a special case solved in the theorem below.  To formulate the result, recall that the linear connections $\nabla$ on a vector bundle ${\mathfrak t}$ are in one-to-one correspondence with the splittings $s\colon{\mathfrak t}\rightarrow J^1{\mathfrak
  t}$ of a natural exact sequence
\begin{equation*}
         0\rightarrow T^*\!M\otimes{\mathfrak t}\monomorphism J^1{\mathfrak 
         t}\rightarrow{\mathfrak t}\rightarrow 0.
\end{equation*}
Now suppose ${\mathfrak t}$ is a Lie algebroid and ${\mathfrak g}
\subset J^1{\mathfrak t} $ an infinitesimal geometric structure. Then
we call $\nabla$ a {\df generator} of ${\mathfrak g}\subset
J^1{\mathfrak t}$ if $s({\mathfrak t_1})\subset{\mathfrak g}$, where
${\mathfrak t}_1\subset{\mathfrak t}$ denotes the image of ${\mathfrak
  g}$.  Generators are certain `preferred connections,' which
generally exist, but need not be unique.  For example, the Levi-Cevita
connection is a generator for the bundle ${\mathfrak g}\subset J^1
(TM)$ of 1-symmetries of a Riemannian metric $\sigma$ (but so is any
linear connection $\nabla$ on $M$ such that $\bar\nabla \sigma=0$,
where $\bar\nabla_U V:=\nabla_V U + [U,V]$).  Generators are
indispensable in explicit computations.

The following crucial observation is not difficult to establish. (For
a proof see Sect.\,\ref{sectionR}.)
\begin{theorem}
                 Let ${\mathfrak g}\subset J^1{\mathfrak t}$ be an 
                 infinitesimal geometric structure and assume the projection $ 
                 a\colon{\mathfrak g}\rightarrow{\mathfrak t}$ has constant 
                 rank.  Then ${\mathfrak g}$ has a unique generator $\nabla$ 
                 if and only if it is surjective and has  structure 
                 kernel ${\mathfrak h}=0$.  In that case $\nabla$ is a Cartan 
                 connection on ${\mathfrak t}$ whose parallel sections are 
                 precisely the symmetries of ${\mathfrak g}$.
\end{theorem}
\noindent In particular, $\curvature\nabla$ is then the local 
obstruction to maximal symmetry.

When geometric structures do not satisfy the hypotheses of Theorem \ref{neoGenerator}, one tries to correct this with an appropriate sequence of prolongation and reduction operations described next.

% A specific algorithm is offered in \ref{hopeful}
% below.

% It should be emphasized that because we will work in a model and 
% coordinate-free way, the above sequence, and in particular $\mathfrak 
% g^{(n)}$ and its Cartan connection, are {\em canonical} objects 
% associated with the initial data $\mathfrak g\subset J^1\mathfrak t$.  
% Analogous sequences constructed via model-based approaches amount to 
% mere {\em representations} of this canonical sequence.  To construct 
% such a representation, one may have to make certain (a priori) 
% assumptions about the model (e.g., existence of invariant splittings 
% of representations of the structure group) or make certain 
% non-canonical choices (caused, for example, by multiplicities in 
% relevant representations of the structure group).  Additionally, in an 
% innately transitive approach reduction cannot be fully carried out, 
% leading to redundancies in the final description.

\subsection{Prolongation}\lab{calendula}
The {\df prolongation} of an infinitesimal geometric structure
${\mathfrak g}\subset J^1{\mathfrak t}$ is a natural `lift' of
${\mathfrak g}$ to a subset ${\mathfrak g}^{(1)}\subset J^1{\mathfrak
  g}$: Noting that $J^1 {\mathfrak g}$ is a subbundle of
$J^1(J^1{\mathfrak t})$, and the existence of a natural inclusion $J^2
{\mathfrak t} \subset J^1(J^1 {\mathfrak t}) $, we define
\begin{equation*}
  {\mathfrak g} ^{(1)}\equiv J^1 {\mathfrak g} \cap J ^2 {\mathfrak t}.
\end{equation*}
It turns out that ${\mathfrak g}^{(1)}$ is an infinitesimal
geometric structure on $\mathfrak g$ whenever ${\mathfrak g} ^{(1)} $
has constant rank. Most importantly, there
is a one-to-one correspondence between symmetries of ${\mathfrak g}$
and symmetries of ${\mathfrak g}^{(1)}$, furnished by prolongation of
sections:
\begin{proposition}
  A section $W \subset {\mathfrak t} $ is a symmetry of $\mathfrak g$ if and
  only if $J^1 W \subset {\mathfrak g} $ is a symmetry of ${\mathfrak
    g} ^{(1)} $.
\end{proposition}
\noindent%
We prove this proposition in Sect.\,\ref{sectionX}.

As an example, consider once again the bundle $\mathfrak{g}\subset
J^1(TM)$ of 1-symmetries of a Riemannian metric $\sigma$. Then it is
possible to show that ${\mathfrak g}^{(1)}\subset J^1 {\mathfrak g}$
is surjective with a trivial structure kernel. Replacing
$\mathfrak{g}$ and $\mathfrak{t}$ in Theorem \ref{neoGenerator} with
${\mathfrak{g}^{(1)}}$ and $\mathfrak{g}$, one establishes the
existence of the Cartan connection $\nabla ^{(1)}$ given explicitly in
\ref{cart} above, as the unique generator of
$\mathfrak{g}^{(1)}$. (For details, see \ref{illustration}.)

When ${\mathfrak g} $ is transitive one can characterize the
prolongation as the joint isotropy of the projection $a\colon
{\mathfrak g} \rightarrow {\mathfrak t} $, viewed as a section $a
\subset {\mathfrak g}^* \otimes {\mathfrak t} $ (the {\df tautological
  one-form}) and its exterior derivative $da \subset
\alternating^2({\mathfrak g})\otimes {\mathfrak t}$ (the {\df
  torsion}). This is explained further in Sect.\,\ref{sectionX}. The
isotropy characterization greatly simplifies computations in the
transitive case. The isotropy characterization is deduced from general
considerations outlined in Sect.\,\ref{sectionX} and these should be
useful in tackling intransitive cases as well. An in-depth analysis of
the transitive case, suitable for making concrete computations, is
given in Sect.\,\ref{method}.

\subsection{Reduction}\lab{reduction}%
Let ${\mathfrak g}\subset J^1{\mathfrak t}$ be an infinitesimal
geometric structure.  By a {\df reduction} of ${\mathfrak g}$ we shall
mean any subalgebroid ${\mathfrak g}'\subset{\mathfrak g}$ with the
same symmetries as ${\mathfrak g}$; it suffices to check that
symmetries of ${\mathfrak g}$ are symmetries of ${\mathfrak g}'$.  In
contrast to prolongation, there is no unique way to construct
reductions.  Notice, however, that if ${\mathfrak g}'\subset{\mathfrak
  g}$ is a reduction and ${\mathfrak g}''\subset{\mathfrak g}$ merely
a subalgebroid satisfying ${\mathfrak g}'\subset{\mathfrak
  g}''\subset{\mathfrak g}$, then ${\mathfrak g}''$ is automatically a
reduction of ${\mathfrak g}$ also. We say ${\mathfrak g}''$ is a {\df
  cruder} reduction than ${\mathfrak g}'$.

We now describe the most important reduction techniques: elementary 
reduction and $\Theta$-reduction.

\subsection{Elementary reduction} \lab{cradle}
Returning to Cartan's method described above, we emphasize that 
transitivity is not a hypothesis of Theorem \ref{neoGenerator} (and 
that Cartan algebroids can be {\em intransitive}).  Rather, one 
requires surjectivity.  If an infinitesimal geometric structure 
${\mathfrak g}\subset J^1{\mathfrak t}$ is {\em not} surjective, we 
may attempt to make it so by passing to the {\df elementary reduction} 
${\mathfrak g}_1$ of ${\mathfrak g}$.  By definition,  
\begin{equation*}
                {\mathfrak g}_1:={\mathfrak g}\cap J^1{\mathfrak 
                t}_1,
\end{equation*}
where ${\mathfrak t}_1\subset{\mathfrak t}$ denotes the image of 
${\mathfrak g}$.  Assuming $\mathfrak t_1\subset \mathfrak t$ and 
$\mathfrak g_1\subset {\mathfrak g}$ have constant rank, they are 
subalgebroids.  In particular, ${\mathfrak g}_1\subset J^1{\mathfrak 
t}_1$ becomes an infinitesimal geometric structure.  Moreover, one 
easily proves:
\begin{proposition}
  If the elementary reduction ${\mathfrak g}_1$ of ${\mathfrak g}$ has
  constant rank then it is a reduction of ${\mathfrak g}$ in the sense
  above. If $\mathfrak g$ is surjective then ${\mathfrak
    g}_1={\mathfrak g}$. Conversely, if ${\mathfrak
    g}_1={\mathfrak g}$ then $\mathfrak g$ is contained in $J^1
  {\mathfrak t}_1$ and is surjective, as an infinitesimal geometric
  structure on ${\mathfrak t}_1$.
\end{proposition}

Because surjectivity is built into the definition of 
$G$-structures, elementary reduction never appears in that setting. 
Elementary reduction is described  further in Sect.\,\ref{elementary}, 
together with a cruder alternative called {\df  image reduction}.

\subsection{$\Theta$-reduction}\lab{koota}
If an infinitesimal geometric structure ${\mathfrak g}\subset
J^1{\mathfrak t}$ is already surjective but has a non-trivial
structure kernel, then, with a view to applying Theorem
\ref{neoGenerator}, one can try to shrink the structure kernel by
prolonging (see above).  However, the prolongation ${\mathfrak
  g}^{(1)}$ generally fails to be surjective itself.  While one might
attempt to correct this by turning to the elementary reduction of
${\mathfrak g}^{(1)}$, there is an alternative that is computationally
more attractive.  One anticipates the lack of surjectivity of
${\mathfrak g}^{(1)}$ by first replacing ${\mathfrak g}$ by its {\df
  $\Theta$-reduction}.  By definition, this is the image of
${\mathfrak g}^{(1)}\subset J^1{\mathfrak g}$, i.e., the set
${\mathfrak g}_1^{(1)}:=p({\mathfrak g}^{(1)})\subset \mathfrak{g}$,
where $ p\colon J^1{\mathfrak g}\rightarrow{\mathfrak g}$ is the
natural projection.  The point is that $\Theta$-reductions can be
computed without directly computing the larger space ${\mathfrak
  g}^{(1)}$.  In Sect.\,\ref{sectionY} we prove:
\begin{proposition}
		 If ${\mathfrak g}_1^{(1)}$ and ${\mathfrak g}^{(1)}$ have 
		 constant rank, then the $\Theta$-reduction ${\mathfrak 
		 g}_1^{(1)}$ of ${\mathfrak g}$ is a reduction of ${\mathfrak 
		 g}$ in the sense of \ref{reduction}.  One has ${\mathfrak 
		 g}_1^{(1)}={\mathfrak g}$ if and only if ${\mathfrak 
		 g}^{(1)}$ is surjective.
\end{proposition}

Classically, `reduction' has usually meant {\df reduction by torsion},
a notion we define in \ref{intrinsic} for arbitrary surjective
infinitesimal geometric structures.  If ${\mathfrak t}=TM$, then
reduction by torsion coincides with $\Theta$-reduction.  More
generally, however, torsion reduction is a cruder reduction technique.
In some cases $\Theta$-reduction is not only more efficient but also
more computationally convenient; conformal geometry
(Sect.\,\ref{point}) is a case in point.  While
$\Theta$-reductions are always defined, we have restricted their
detailed analysis to the case of surjective structures over {\em
  transitive} Lie algebroids. This analysis appears in
Sect.\,\ref{method}.

% In general the elementary reduction ${\mathfrak g}_1$ fails to be 
% surjective.  However, one may then turn to the elementary reduction 
% $\mathfrak g_2$ of $\mathfrak g_1$ and, iterating further, obtain a 
% descending sequence $\mathfrak g,\mathfrak g_1,\mathfrak g_2, \ldots$ 
% of reductions, and a corresponding descending sequence of images 
% ${\mathfrak t},{\mathfrak t}_1,{\mathfrak t}_2,\ldots$.  If 
% ${\mathfrak g}_j$ and ${\mathfrak t}_j$ have constant rank for every 
% $j$, then ${\mathfrak g}$ is called {\df non-singular}.  In that case 
% the sequences must stabilize: There is some finite $n$ such that 
% ${\mathfrak g}_j={\mathfrak g}_n$ and ${\mathfrak t}_j={\mathfrak 
% t}_n$ for all $ j\ge n$.  In \ref{} below we call ${\mathfrak g}_n$ 
% the {\df surjective reduction} of ${\mathfrak g}$.

% \begin{proposition}
%   If ${\mathfrak g}$ is a non-singular infinitesimal geometric 
%    structure on $\mathfrak t$, then its reduction $\mathfrak g_n$ is a 
%   reduced infinitesimal geometric structure on $\mathfrak t_n$ of 
%   ${\mathfrak t}$ having the same symmetries as ${\mathfrak g}$.
% \end{proposition}

\subsection{A specific normalizing algorithm and its limitations}\lab{hopeful}%
We now offer a specific algorithm for constructing a Cartan algebroid out of an infinitesimal geometric structure of finite-type. First, we define two auxiliary
procedures. 

By Proposition \ref{cradle}, the following procedure, which we shall
call {\sffamily surjectify ${\mathfrak g} $}, forces ${\mathfrak g}$
to be surjective:%
{\sffamily
\begin{tabbing}
{\qquad} \= {\qquad} \kill
do while ${\mathfrak g} $ is not surjective \> \\
\> replace $ {\mathfrak g} $  with ${\mathfrak g}_1$ 
   {\normalfont (elementary reduction)}\\
end do. \>
\end{tabbing}
}%
\noindent Next, we let {\sffamily strongly surjectify $\mathfrak g$}
denote the following procedure making $\mathfrak g$ and ${\mathfrak g}
^{(1)} $ simultaneously surjective (by Propositions \ref{cradle} and
\ref{koota}): {\sffamily
\begin{tabbing}
{\qquad} \= {\qquad} \kill
do while $\mathfrak g ^{(1)} $ is not surjective \> \\
\> surjectify $\mathfrak g$  \\
\> replace $ {\mathfrak g} $  with ${\mathfrak g}_1 ^{(1)}$  
   {\normalfont ($\Theta $-reduction)}\\
\> surjectify $\mathfrak g$\\
end do. \>
\end{tabbing}
}%
\noindent To describe an implementation of this procedure it evidently suffices to describe $\Theta $-reduction in the special surjective case. 

One might attempt to normalize an infinitesimal geometric structure using elementary reduction and prolongation alone. In practice, however, it is generally easier to apply the following algorithm:
{\sffamily
\begin{tabbing}
{\qquad} \= {\qquad} \kill
surjectify $\mathfrak g$ \> \\
% if ${\mathfrak h}= 0$ apply {\normalfont Theorem \ref{neoGenerator}} 
%   and stop \> \\
repeat until stop encountered \> \\
\> if ${\mathfrak h}= 0$ apply {\normalfont Theorem \ref{neoGenerator}} 
  and stop \\
\> strongly surjectify $\mathfrak g$ \\
\> if ${\mathfrak h}= 0$ apply {\normalfont Theorem \ref{neoGenerator}} 
  and stop \\
\> replace $\mathfrak g$ with $\mathfrak g ^{(1)} $  
   {\normalfont (prolongation)}\\ 
end repeat. \> 
\end{tabbing}
}%
\noindent Notice that prolongation is delayed as long as possible.  We
now describe the ways in which the above algorithm can fail.

Firstly, an execution of {\sffamily surjectify $\mathfrak g$ } or
{\sffamily strongly surjectify $\mathfrak g$ } could fail because
$\mathfrak g$, at some iteration of these procedures' {\sffamily
  do-while} loops, loses the constancy of its rank. While prolongation
of $\mathfrak g$ might resolve this kind of singularity (by recovering
rank constancy), this requires a prolongation theory for `variable
rank Lie algebroids' (or Lie pseudoalgebras) which is not provided
here. Similarly, $\mathfrak g$ might lose rank constancy at some
iteration of the {\sffamily repeat-until} loop of the main algorithm.

Even if all singularities are successfully resolved, it may happen
that no {\sffamily stop} is ever encountered in the {\sffamily
  repeat-until} loop, which then becomes perpetual. However, such cases can be detected by applying (a Lie algebroid version of) Cartan's involutivity test; they occur for objects of infinite-type which are not described here.

Another possibility for failure concerns $\Theta$-reduction. In
practice, it seems rather complicated to implement without making the
added assumption that the base algebroid ${\mathfrak t}$ of $\mathfrak
g \subset J^1 {\mathfrak t} $ is {\em transitive}. One way to handle
intransitivity and singularities might be to {\em restrict}, in some
way, the infinitesimal geometric structure ${\mathfrak g}$ to each
orbit of ${\mathfrak t}$. This restriction will sit over a transitive
base (this being a Lie algebroid over the orbit) but is not simply the
pullback in the category of Lie algebroids, for one wants an
infinitesimal geometric structure over the base, not merely a Lie
algebroid.  Also, one needs to understand how conclusions regarding
the restricted structure combine with transverse information to solve
the original problem.  Fortunately, a splitting theory for Lie
algebroids exists \cite{Fernandes_02} and this possibly reduces the
transverse problem to the case of an isolated singularity
(zero-dimensional orbit).  None of this is explored here either.

If the Cartan algorithm above succeeds it ends in an application of
Theorem \ref{neoGenerator}, this delivering a Cartan algebroid whose
parallel sections are in natural one-to-one correspondence with the
symmetries of $\mathfrak g$. We then say that $\mathfrak g$ {\df
  admits an associated Cartan algebroid}.

\subsection{Paper outline}
In Sect.\,\ref{section2} we review basic Lie algebroid notions and
establish attendant notation. In particular, we describe the
generalizations of linear (Koszul) connections afforded by Lie
algebroids, these connections amounting to deformations of Lie
algebroid representations, which are also defined. We recall the
definition of the adjoint representation of $J^1{\mathfrak g} $ on a Lie
algebroid ${\mathfrak g}$, and write down the bracket on
$J^1{\mathfrak g}$ explicitly.  We introduce the notion of {\df
  associated connections}, which are ubiquitous throughout, and
developed further in Sect.\,\ref{sectionR}.

Sect.\,\ref{geomet} summarizes features of Cartan algebroids
established in \cite{Blaom_05}, stating, in particular, the result
that Cartan algebroids are deformations of infinitesimal symmetries
(Theorem \ref{dissymmetry}).

Sect.\,\ref{sectionQ} gives many examples of infinitesimal geometric
structures, arising as the isotropy subalgebroids associated with
various structures in differential geometry. We also explain how to
associate an infinitesimal geometric structure with a Cartan algebroid
or a classical $G$-structure. From our description of affine
structures, it will be clear how one may associate an infinitesimal
geometric structure with an arbitrary (but suitably regular)
differential operator on $M$. We demonstrate in practice how one
computes the image of an infinitesimal geometric structure, without
resorting to local coordinate calculations.

Associated with an infinitesimal geometric structure ${\mathfrak g} \subset J^1 t$, with structure kernel ${\mathfrak h} $ and image ${\mathfrak t}_1$, is an exact sequence
\begin{equation*}
0 \rightarrow   {\mathfrak h} \rightarrow {\mathfrak g} \rightarrow {\mathfrak t}_1 \rightarrow 0.
\end{equation*}
A generator $\nabla$ of $\mathfrak g$, as defined in
\ref{neoGenerator} above, above amounts to a splitting of this
sequence, determining an identification ${\mathfrak g} \cong
{\mathfrak t}_1 \oplus {\mathfrak h}$, which is very useful in
computations. Sect.\,\ref{sectionR} characterizes the linear
connections $\nabla$ on ${\mathfrak t}$ arising in this way and writes
down the Lie algebroid structure induced on ${\mathfrak t}_1 \oplus
{\mathfrak h}$. Additionally, we express the geometric Bianchi
identities in terms of generators and discuss the invariant differential
operators associated with a structure. This section can be lightly
scanned on a first reading.

Sect.\,\ref{elementary} shows how to explicitly compute an elementary
reduction of $\mathfrak g$ once a generator is known. A simple
application to functions on a Riemannian three-manifold illustrates
the technique.

Sections \ref{sectionX} and \ref{sectionY} are important theoretical
parts of the paper, describing prolongation, torsion reduction, and
$\Theta $-reduction, and the relationships between these. An immediate
application is Theorem \ref{goosey}, which furnishes general
conditions under which an infinitesimal geometric structure
${\mathfrak g} \subset J^1{\mathfrak t} $ is its own associated Cartan
algebroid, Riemannian geometry being a case in point. Keen to
illustrate the result in this and other concrete examples, we include
additional detail in the special case ${\mathfrak t}=TM $ (Theorem
\ref{goosey2}) that is not established until Sect.\,\ref{method}. The
reader preferring a more linear presentation may skip to
Sect.\,\ref{method} immediately after \ref{goosey}, returning to
\ref{goosey2} and the remainder of the paper thereafter.

We have not attempted substantially novel applications of Cartan's
method in the present work. In particular, our applications to
intransitive phenomena are fairly superficial. We hope to correct this
deficiency elsewhere. For the purposes of making comparisons with
other approaches, we give a detailed application to subriemannian
contact three-manifolds in Sect.\,\ref{contact}. An analysis using
$G$-structures is to be found in \cite{Hughen_95,Montgomery_02}. In
addition to constructing obstructions to symmetry, we go on to
construct the invariant differential operators.

In Sect.\,\ref{method} we return to prolongation, explaining in detail
how to `prolong' a generator, and hence how compute prolongations in
practice. This analysis includes the general case ${\mathfrak t} \ne
TM$, but we must assume ${\mathfrak t}$ is {\em transitive} to avoid
complications. A detailed section on conformal geometry,
Sect.\,\ref{point}, illustrates the more general prolongation
results.

\section{Preliminary notions}\lab{section2}
% This section is a rapid review of the theory of connections, from the 
% point of view of Lie algebroid representations, our main purpose being 
% to establish notation and terminology.  Of particular relevance will 
% be the {\df adjoint} representation of $ J^1{\mathfrak g}$ on a Lie 
% algebroid ${\mathfrak g}$, and a concrete description of the bracket 
% on $ J^1{\mathfrak g}$.  We introduce the notion of {\df associated}
% connections, which are developed further in Sect.\,\ref{sectionR}.

For an introduction to Lie groupoids and algebroids, see
\cite{CannasdaSilva_Weinstein_99} or \cite{Mackenzie_05}.
%\vspace{\baselineskip}%
%\noindent 
All constructions in this paper are made in the $ C^\infty$ 
category.
\subsection*{Notation} We use $\alternating^k(V)\cong 
\mathsf\Lambda^k(V^*)$ and $\symmetric^k(V)\cong\mathsf S^k(V^*)$ to 
denote the spaces of ${\mathbb R}$-valued alternating and symmetric $ 
k$-linear maps on a vector space $ V$.  Similar notation applies to 
the tensor algebra of a vector bundle $ E$ over $M$.  If $\sigma$ is a 
section of $ E$, then this is indicated by writing $\sigma\in \Gamma 
(E)$ or $\sigma\subset E$.  Thus $\sigma\subset\alternating^2 
(TM)\otimes E$ means $\sigma$
is an $ E$-valued differential two-form on $M$.%
\subsection{Lie algebroids}\lab{Lie}
A {\df Lie algebroid} over $M$ consists of a vector bundle ${\mathfrak g}$ over $M$, a Lie bracket $[\,\cdot\,,\,\cdot\,]$ on the space of sections $ 
\Gamma ({\mathfrak g})$, together with a vector bundle morphism 
$\#\colon{\mathfrak g}\rightarrow TM$, called the {\df anchor}.  One 
requires that the bracket satisfy a Leibnitz identity,
\begin{equation*}
        [ X,fY]=f[X, Y]+df(\#X)Y,
\end{equation*}
where $ f$ is an arbitrary smooth function.  The tangent bundle $TM$,
equipped with the Jacobi-Lie bracket on vector fields, and identity
map as anchor, is a Lie algebroid; in the general case one demands
that the induced map of sections $\#\colon \Gamma ({\mathfrak
  g})\rightarrow \Gamma (TM)$ be a Lie algebra homomorphism. 
  
A Lie algebroid is {\em transitive} if its anchor is surjective. The
image of the anchor is an integrable distribution. The leaves of the
integrating foliation, which is possibly singular, are the {\df orbits}
of the Lie algebroid. A subbundle $\mathfrak h\subset\mathfrak
g$ is a {\df subalgebroid} if $\Gamma(\mathfrak
h)\subset\Gamma(\mathfrak g)$ is a subalgebra.  Lie {\em algebras} are
simply the Lie algebroids over a single point.\newline

Linear connections on $M$ may be interpreted as 
representations of the tangent bundle $TM$, deformed by
curvature. This point of view leads to the more general notion of a
$\mathfrak g${\df -connection}, defined in \ref{representations}
below.  Proposition \ref{representations} gives an `analytic'
characterization of these connections, generalizing the usual one,
which the reader may take as an alternative definition.

\subsection{The definition of $\gl(E)$ for a vector bundle
  $E$}\lab{glE}%
A representation of a Lie algebra ${\mathfrak g}$ is a vector
space $E$, together with a Lie algebra homomorphism ${\mathfrak
  g}\rightarrow\gl (E):=\homomorphism (E, E)$.  Turning now to the
generalization of $\gl(E)$ relevant to Lie {\em algebroid}
representations, let $E$ be a vector {\em bundle} over the same base $M$,
and consider the exact sequence
\begin{equation*}
          0\rightarrow T^*\!M\otimes E\monomorphism J^1 E\rightarrow E
\rightarrow 0.
\end{equation*}
Here $T^*\!M\otimes E\monomorphism J^1 E$ is the inclusion which, as a
map on sections, sends\footnote{An opposite sign convention is adopted
  in our paper \cite{Blaom_05} and elsewhere. The present convention
  avoids unpleasant sign changes in \eqrefs{adjective}{gull} and
  \eqrefs{garden}{bird}.}   $df\otimes\sigma$ to
$fJ^1\sigma-J^1(f\sigma)$. Applying $\homomorphism (\,\cdot\,, E)$ to the sequence,
and identifying $\homomorphism (T^*M\otimes E, E)$ with
$TM\otimes\homomorphism (E, E)$, we obtain a second exact sequence
\begin{equation*}
   0\rightarrow\homomorphism (E, E)\monomorphism\homomorphism (J^1 E,E)
   \xrightarrow{\#}TM\otimes\homomorphism(E, E)\rightarrow 0.
\end{equation*}
Noticing that there is natural inclusion ${TM}\monomorphism {TM}
\otimes \homomorphism(E,E)$ via the map $v\mapsto v\otimes\identity$,
we define $\gl(E)\subset \homomorphism(J^1 E,E)$ to be the preimage of
$TM$ under the surjective arrow $\#$, and obtain a third exact
sequence
\begin{equation*}
   0\rightarrow\homomorphism (E, E)\monomorphism \gl(E)
   \xrightarrow{\#}TM\rightarrow 0.
\end{equation*}
\begin{proposition}Regard each section $D$ of $\homomorphism(J^1 E,E)$ as a differential operator $D\colon\Gamma(E)\rightarrow \Gamma(E)$. Then:
 \begin{conditions}
 \item A section $D\subset \homomorphism(J^1 E,E)$ lies in $\gl(E)$ if and only if there exists a vector field $V$ such that
 \begin{equation*}
         D(f\sigma)=fD\sigma+ df(V)\sigma
 \end{equation*}
 for all sections $\sigma$ of $E$ and functions $f$; in that case
 $V=\#D$.
 \item The operator commutator bracket,
\begin{equation*}
  [D_1,D_2]_{\mathfrak{gl}(E)}\,\sigma:=D_1D_2 \sigma-D_2D_1 \sigma,
\end{equation*}
makes $\gl(E)$ into a Lie algebroid with anchor $\#$.
 \end{conditions}
\end{proposition}
\begin{remark}
  $\gl(E)$ is in fact a realization of the Lie algebroid of the Lie
  groupoid $\operatorname{GL}(E)$ of isomorphisms between fibres of
  the vector bundle $E$, i.e., `relative frames' on $E$. So elements
  of $\gl(E)$ have the interpretation of `infinitesimal moving
  frames.' For details, see \cite{Mackenzie_05} (where
  $\mathfrak{gl}(E)$ is denoted ${\mathcal D}(E)$).
\end{remark}

\subsection{Lie algebroid representations and $\mathfrak g$-connections}%
  \lab{representations} A {\df representation} of a Lie algebroid 
  ${\mathfrak g}$ on $ E$ is a morphism ${\mathfrak 
  g}\rightarrow{\mathfrak{gl}}(E)$ of Lie algebroids.  When $M$ is a 
  single point one recovers the usual representations of a Lie {\em 
  algebra}. Deforming the representation notion we arrive at 
  the following:
\begin{definition}
         Let ${\mathfrak g}$ be any Lie algebroid over $M$. A ${\mathfrak 
         g}$-{\df connection} on a vector bundle $ E$ over $M$ is a vector 
         bundle morphism
         \begin{equation*}
                \nabla\colon{\mathfrak g}\rightarrow{\mathfrak{gl}}(E)
         \end{equation*}
         that is {\em not} required to be a Lie algebroid morphism,
         but {\em is} nevertheless required to respect the anchors,
         $\#\nabla=\nabla\#$.
\end{definition}
\noindent%
Suppose $X$ is a section of $\mathfrak g$. When the section
$\nabla (X)$ of ${\mathfrak{gl}}(E)\subset\homomorphism (J^1 E, E)$ is
to be viewed as a differential operator,  we instead write
$\nabla_X$, i.e., $\nabla_XV:=\nabla (X)(J^1 V)$. In view of the
preceding characterization of the sections of ${\mathfrak{gl}}(E)$ as
derivations, we have the Leibnitz identity
\begin{equation*}
                \nabla_X(f\sigma)=f\nabla_X\sigma+ df(\#X)\sigma;\qquad 
                X\subset{\mathfrak g},\sigma\subset E.
\end{equation*}
Conversely:
\begin{proposition}
  Every vector bundle morphism $\nabla\colon\mathfrak 
  g\rightarrow\homomorphism(J^1 E,E)$ that is Leibnitz in the above 
  sense is a $\mathfrak g$-connection.
\end{proposition}

If $\nabla$ is a $\mathfrak g$-connection, then the formula
\begin{equation*}
        \curvature\nabla\, (X, Y):=[\nabla (X),\nabla 
(Y)]_{\mathfrak{gl}(E)}-\nabla (\,[X, Y]_{\mathfrak g}),
\end{equation*}
defining the Lie algebroid curvature of the map
$\nabla\colon\mathfrak g\rightarrow\gl(E)$ takes a familiar form:
\begin{equation*}
        \curvature\nabla\,(X, Y) 
        Z=\nabla_X\nabla_YZ-\nabla_Y\nabla_XZ-\nabla_{[X, Y]_{\mathfrak g}}Z.
\end{equation*}
\noindent The $\mathfrak g$-connection $\nabla $ is a $\mathfrak
g$-representation when $\curvature \nabla = 0$.
\begin{example}
  If ${\mathfrak g}$ is a Lie algebroid and $ E\subset{\mathfrak g}$
  is a subalgebroid contained in the kernel of its anchor then a
  canonical representation $\rho$ of ${\mathfrak g}$ on $ E$ is well
  defined by $\rho_XY:=[X, Y]_{\mathfrak g}$.  Important cases in
  point are the kernel of the anchor itself, and the structure kernel
  of an infinitesimal geometric structure, when these have constant
  rank.
\end{example}

\subsection{Linear connections}\lab{linear}
Using the language of the preceding discussion, a linear connection 
$\nabla$ on $E$ is just a $ TM$-connection on $E$, becoming a 
$TM$-representation when $\nabla$ is flat.  It is an elementary fact 
that the linear connections on $E$ are in one-to-one correspondence 
with the splittings $s\colon E\rightarrow J^1 E$ of the exact 
sequence
\begin{equation*}
          0\rightarrow T^*\!M\otimes E\monomorphism J^1 E\rightarrow E
\rightarrow 0.
\end{equation*}
The splitting associated with a linear
connection $\nabla$ on $ E$ is, as a map on sections, given by
\begin{equation}
         s\sigma=J^1\sigma+\nabla\sigma;\qquad\sigma\subset 
         E.\mathlab{flower}
\end{equation}
Here $\nabla\sigma\subset T^*\!M\otimes E$ is defined by $(\nabla 
\sigma)(V):=\nabla_V\sigma$.

\subsection{Jet bundles}\lab{jets}
If ${\mathfrak g}$ is any Lie algebroid over $M$ with anchor $\#$, 
then the vector bundle $ J^k{\mathfrak g}$ of $k$-jets of sections of 
${\mathfrak g}$ is another Lie algebroid.  Its anchor is the composite
\begin{equation*}
        J^k{\mathfrak
    g}\rightarrow{\mathfrak g}\xrightarrow{\#}TM.
\end{equation*}
Its bracket is uniquely determined by requiring that $ k$th-order
prolongation $$J^k\colon \Gamma ({\mathfrak g})\rightarrow \Gamma
(J^k{\mathfrak g})$$ be a morphism of Lie algebras (but see also
\ref{garden} below).  One may view $J^k$ as a functor from the
category of Lie algebroids to itself \cite{Crainic_Fernandes_05}.

As this paper will demonstrate, many naturally occurring Lie 
algebroids in differential geometry can be constructed from $TM$ using 
just two operations: application of the prolongation functor 
$J^1(\,\cdot\,)$; and passage to a subalgebroid.
\subsection{The adjoint representation}\lab{adjective}
The generalization of Lie algebra adjoint representations to a Lie
algebroid $\mathfrak g$ is not a self-representation; it is rather a
representation of $ J^1{\mathfrak g}$ on ${\mathfrak g}$. This representation is well-defined by
\begin{equation*}
        \adjoint_{J^1 X}^{\mathfrak g}Y=[X, Y]_{\mathfrak g}.
\end{equation*}
Using the identity
\begin{equation}
        [ J^1 X, J^1 Y]_{J^1{\mathfrak g}}=J^1[ X, Y]_{\mathfrak 
        g},\mathlab{Scrooge}
\end{equation}
one shows that $\adjoint^{\mathfrak g}$ is indeed a representation (and not merely a $J^1 {\mathfrak g}$-connection).

We note that 
\begin{equation}
        \adjoint_\phi^{\mathfrak g}X=\phi (\#X);\qquad X\subset{\mathfrak g},
        \mathlab{gull}
\end{equation}
for all sections $\phi\subset T^*\!M\otimes{\mathfrak g}\subset 
J^1{\mathfrak g}$. If $ a\colon{\mathfrak g}\rightarrow{\mathfrak h}$ 
is a morphism of Lie algebroids, then one has the identity
\begin{equation}
         a\Big(\,\adjoint_\xi^{\mathfrak g}X\,\Big)=\adjoint_{(J^1 
         a)\xi}^{\mathfrak h}(aX);\qquad\xi\subset J^1{\mathfrak g},\, 
         X\subset{\mathfrak g}.\mathlab{America}
\end{equation}

\subsection{The bracket on $J^1(\,\cdot\,)$ of a Lie algebroid}\lab{garden}
Recall that the bracket on $ J^1{\mathfrak g}$ is implicitly defined
by the requirement \eqrefs{adjective}{Scrooge}. With the help of the
adjoint representation, we now describe this bracket concretely.

Although the exact sequence
\begin{equation}
         0\rightarrow T^*\!M\otimes{\mathfrak g}\monomorphism J^1{\mathfrak 
         g}\rightarrow{\mathfrak g}\rightarrow 0\mathlab{vegetable}
\end{equation}
possesses no canonical splitting, the corresponding exact sequence of 
section spaces,
\begin{equation*}
         0\rightarrow  \Gamma (T^*\!M\otimes{\mathfrak g})\monomorphism
         \Gamma (J^1{\mathfrak 
         g})\rightarrow\Gamma ({\mathfrak g})\rightarrow 0,
\end{equation*}
is split by $ J^1\colon \Gamma ({\mathfrak g})\rightarrow \Gamma 
(J^1{\mathfrak g})$,   delivering a canonical identification
\begin{equation*}
         \Gamma (J^1{\mathfrak g})\cong \Gamma ({\mathfrak g})\oplus \Gamma 
         (T^*\!M\otimes{\mathfrak g}).
\end{equation*}
Under this identification, the Lie  algebra $ \Gamma (J^1{\mathfrak 
g})$ is a semidirect product described in the proposition below.

In addition to having the adjoint representation of $ J^1 {\mathfrak g}$ 
on ${\mathfrak g}$, we have a  representation of $ J^1 
{\mathfrak g}$ on $ TM$,  given by the composite
\begin{gather}
                 J^1{\mathfrak g}\xrightarrow{J^1\#}{ J^1 
                 (TM)}\xrightarrow{\adjoint^{TM}}\gl 
                 (TM),\mathlab{korgit}\\
                 \text{i.e.,}\quad J^1 X\cdot V=[\#X, V];\qquad
                 X\subset {\mathfrak g},\,V\subset TM.\notag
\end{gather}
So we can construct a natural representation of $J^1 {\mathfrak g}$ 
on $T^*\!M\otimes{\mathfrak g}$; it is given by
\begin{equation}
         (J^1 X\cdot\phi) V=[X,\phi (V)]_{\mathfrak g}-\phi ([\#X, 
         V]_{TM});\qquad V\subset TM.\mathlab{exclamation}
\end{equation}
On the other hand, $ T^*\!M\otimes{\mathfrak g}$ is the structure
kernel of $ J^1{\mathfrak g}$ so that $ J^1{\mathfrak g}$ acts on $
T^*\!M\otimes{\mathfrak g}$ via bracket; see the example in
\ref{representations}.  A consequence of the following result is that
these two representations coincide.
\begin{proposition}
  The subalgebroid $T^*\!M\otimes{\mathfrak g}\subset J^1{\mathfrak
    g}$ inherits the bracket
\begin{equation}
                [\phi_1,\phi_2]_{T^*\!M\otimes{\mathfrak g}}\,(V)=\phi_1(\# 
                \phi_2 V)-\phi_2(\# \phi_1 V);\qquad V\subset TM,\mathlab{bird}
\end{equation}
where $\#$ denotes the anchor of ${\mathfrak g}$.  Sections of 
$J^1{\mathfrak g}$ are of the form $ J^1 X+\phi$ for uniquely 
determined sections $ X\subset{\mathfrak g}$ and $\phi\subset 
T^*\!M\otimes{\mathfrak g}$.  The bracket on $ J^1{\mathfrak g}$ is 
given by
\begin{equation*}
        [ J^1 X+\phi, J^1 Y+\psi]_{J^1{\mathfrak g}}=
        J^1[ X, Y]_{\mathfrak g}+J^1 X\cdot\psi-J^1 
        Y\cdot\phi+[\phi,\psi]_{T^*\!M\otimes{\mathfrak g}},
\end{equation*}
where the action of $ J^1{\mathfrak g}$ on $T^*\!M\otimes 
{\mathfrak g}$ is defined by \eqref{exclamation}.
\end{proposition}
\noindent To prove the proposition one uses 
\eqrefs{adjective}{Scrooge} and the fact that sections of 
$T^*\!M\otimes{\mathfrak g}$ are finitely generated by those of the 
form $df\otimes X=f J^1 X - J^1 (fX)$.

\subsection{Dual connections, torsion, and associated 
connections}\lab{dual}
Let $\mathfrak g$ be a Lie algebroid and $\nabla$ a ${\mathfrak
  g}$-connection on itself. We define the
{\df dual} of $\nabla$ to be the ${\mathfrak g}$-connection $\nabla^*$
on ${\mathfrak g}$ defined by
\begin{equation*}
        \nabla^*_XY:=\nabla_YX+[X, Y].
\end{equation*}
One has `duality' in the sense that $\nabla^{**}=\nabla$.

The {\df torsion} of $\nabla$ is the section, $\torsion\nabla$, of
$\alternating^2 ({\mathfrak g})\otimes{\mathfrak g}$ measuring the
difference between $\nabla$ and its dual:
\begin{equation*}
        \torsion\nabla\, (X, Y):=\nabla_XY-\nabla^*_XY
        =\nabla_X Y-\nabla_YX-[X, Y].
\end{equation*}
The torsion or curvature of $\nabla$ can be expressed in terms  of 
the torsion and curvature of $\nabla^*$ (and, by duality, vice versa):
\begin{proposition}
        Let $\nabla$ be a ${\mathfrak g}$-connection on ${\mathfrak g}$, and 
        $\nabla^*$ its dual. Then:
        \begin{gather}
           \torsion\nabla=-\torsion\nabla^*\\
           \left\{
           \begin{split}
             \curvature\nabla\,(X, Y)Z&=(\nabla_Z^*\torsion\nabla^*)(X, Y)
                       -\curvature\nabla^*\,(Z, X)Y\\
                       &+\curvature\nabla^*\,(Z, Y)X;
                       \qquad X,Y,Z\subset\mathfrak g.
           \end{split}\right.\mathlab{dual2}
        \end{gather}
%         In particular:
%         \begin{conditions}
%         \item If $\curvature\nabla^*=0$, then $\curvature\nabla=0$ if
%           and only if $\torsion\nabla^*$ is $\nabla^*$-parallel.
%         \item\lab{dual4} If $\nabla=\nabla^*$ (no torsion), then
%             \begin{equation*}
%               \curvature\nabla\,(X,Y)Z+\curvature\nabla\,(Y,Z)X
%               +\curvature\nabla\,(Z,X)Y=0,\mathlab{BianchiI}
%             \end{equation*}
%             for arbitrary $X,Y,Z\in\mathfrak g$.
%         \end{conditions}
\end{proposition}

We now introduce two important connections generalizing the dual of a 
linear connection on $TM$.  They are examples of {\df associated 
connections}, defined more generally in \ref{neoAssociated}.  
% Associated connections will be ubiquitous in our study of 
% infinitesimal geometric structures.

Let $\nabla$ be an arbitrary {\em linear} (i.e., $TM$-) connection on 
a Lie algebroid ${\mathfrak g}$.  The {\df associated ${\mathfrak 
g}$-connection} on ${\mathfrak g}$ is defined by
\begin{equation*}
        \bar\nabla_XY=\nabla_{\#Y}X+[X,Y]_{\mathfrak g};\qquad X,\, 
        Y\subset{\mathfrak g}.
\end{equation*}
The {\df associated ${\mathfrak g}$-connection} on  $ TM$ is defined by
\begin{equation*}
                \bar\nabla_XV=\#\nabla_VX+[\#X,V]_{TM};
                        \qquad X\subset{\mathfrak g},\,V\subset TM.
\end{equation*}
The anchor $\#\colon{\mathfrak g}\rightarrow TM$ is equivariant with 
respect to these connections: 
$\bar\nabla_X\#Y=\#\bar\nabla_XY$.  If ${\mathfrak g}=TM$ both 
connections coincide with the dual of $\nabla$. 
\section{Cartan algebroids}\lab{geomet}
This section summarizes some results of \cite{Blaom_05} where Cartan 
algebroids were introduced.  We also define the {\df cocurvature} of 
an arbitrary linear connection on a Lie algebroid, a notion to 
reappear in the context of generators of infinitesimal geometric 
structures.  Some further remarks linking {\em flat} Cartan algebroids 
with Lie pseudogroups appear in Appendix \ref{global}.

As mentioned in Sect.\,\ref{toodle}, Cartan algebroids may be regarded
as deformations of action algebroids (see Theorem \ref{dissymmetry}
below). We recall the definition of the latter first.

\subsection{Action algebroids}\lab{actions}
Let ${\mathfrak g}_0 $ be a finite-dimensional Lie algebra with
bracket $[\,\cdot\,,\,\cdot\,]_{{\mathfrak g}_0}$ acting smoothly on a
manifold $M$. That is, prescribe some Lie algebra homomorphism
$\rho\colon {\mathfrak g}_0 \rightarrow \Gamma (TM)$. We may regard
such actions as abstract infinitesimal symmetries.  The trivial bundle ${\mathfrak g}\equiv{\mathfrak g}_0\times M$ over $M$ has a natural Lie algebroid
structure.  This is the associated {\df action} (or {\df
  transformation}) {\df algebroid}.
                  
The anchor of this Lie algebroid is the `action map'
$\#\colon{\mathfrak g}_0\times M\rightarrow TM$ sending $ (\xi, m)$ to
$\rho\xi(m)$.  The Lie bracket on $\Gamma ({\mathfrak g}_0\times M)$
is a natural extension of the bracket on ${\mathfrak g}_0$, regarded
as the subspace ${\mathfrak g}_0\subset \Gamma ({\mathfrak g}_0\times
M)$ of constant sections.  To describe it, let $\tau
(\,\cdot\,,\,\cdot\,)$ denote the naive extension of this bracket,
i.e., the one that is bilinear with respect to all smooth functions
(and consequently {\em not} Leibnitz):
\begin{equation*}
        \tau (X, Y)(m):=[X(m), Y(m)]_{{\mathfrak g}_0};\qquad 
        X,Y\in\Gamma(\mathfrak g_0\times M).
\end{equation*}
And let $\nabla$ denote the canonical flat connection on ${\mathfrak
  g}_0\times M$.  Then the Lie bracket on $\Gamma ({\mathfrak
  g}_0\times M)$ is defined by
\begin{equation}
        [ X, Y]:=\nabla_{\#X}Y-\nabla_{\#Y}X+\tau (X,Y).\mathlab{trousers}
\end{equation}
Notice that if $\bar\nabla $ denotes the associated $\mathfrak g$-connection on $\mathfrak g$ then $\torsion \bar\nabla=\tau $.
\subsection{Cartan connections}\lab{algebra}
Let $\nabla$ be a linear connection on a Lie algebroid ${\mathfrak
  g}$.  Then $\nabla$ is a {\df Cartan connection} if the
corresponding splitting 
\begin{gather}
  s_\nabla\colon {\mathfrak g} \rightarrow J^1 {\mathfrak g}\notag\\
  s_\nabla \sigma :=J^1 \sigma + \nabla \sigma \mathlab{qookie}
\end{gather}
of the exact sequence of Lie algebroids,
\begin{equation*}
                 0\rightarrow T^*\!M\otimes{\mathfrak g}\monomorphism J^1{\mathfrak 
                 g}\rightarrow{\mathfrak g}\rightarrow 0,
\end{equation*}
is a Lie algebroid morphism.
\noindent A {\df Cartan algebroid} is a Lie algebroid equipped with a 
Cartan connection.  A {\df morphism} of Cartan algebroids is simply a 
connection-preserving morphism of the underlying Lie algebroids.  

It follows immediately from the definition that the ${\mathfrak
  g}$-connections $\bar\nabla$ on ${\mathfrak g}$ and $TM$, associated
with a Cartan connection $\nabla$, are flat i.e, are $\mathfrak
g$-representations.  In particular, every Cartan algebroid ${\mathfrak
  g}$ has a canonical self-representation.  For a converse statement,
see the corollary below.
%
% We do not know any general existence results for Cartan connections on
% a prescribed Lie algebroid $\mathfrak g$. If $\mathfrak g$ is a
% transitive subalgebroid of $J^1(TM)$, and $M$ is {\em contractible},
% then existence follows from a general result of
% Mackenzie\footnote{Warning: In the cited work `connection' means {\df
%     splitting}.} \cite[Theorem 8.2.1]{Mackenzie_05}. Mackenzie's
% classification of Lie algebroid extensions probably delivers results
% more generally but we have not attempted such an analysis here.  While
% Cartan's method of equivalence constructs Cartan connections, the Lie
% algebroid on which the connection is constructed is not under direct
% control, being prolonged or reduced as appropriate.

\subsection{Cocurvature}\lab{cocurvature}
Associated with an arbitrary linear connection $\nabla$ on ${\mathfrak
  g}$ is what we call its {\df cocurvature}.  This is a section of
$\alternating^2 ({\mathfrak g})\otimes T^*\!M\otimes{\mathfrak g}$,
denoted $\cocurvature\nabla$, measuring the lack of `compatibility' of
$\nabla$ with the Lie algebroid structure of ${\mathfrak g}$:
\begin{equation}\begin{split}
        \cocurvature\nabla\,( X, Y)V&:=\nabla_V[X, Y]-[\nabla_VX, 
        Y]-[X,\nabla_VY]\\&+\nabla_{\bar\nabla_XV}Y-\nabla_{\bar\nabla_YV}X.\mathlab{coc}
\end{split}
\end{equation}
Here $\bar\nabla$ denotes the associated ${\mathfrak g}$-connection on 
$ TM$; see \ref{dual} above.  The name `cocurvature' comes from 
\eqref{costly} below.

Define $ s_\nabla$ as in \eqref{qookie} above.  Then the description
of the bracket on $ J^1{\mathfrak g}$ in \ref{garden} establishes the
result \eqref{dgi} below. The other conclusions of the following
proposition are then straightforward.
\begin{proposition}\mbox{}
        \begin{conditions}
                            \item\lab{dgi}  $-\cocurvature\nabla (X, Y)
                                =\curvature s_\nabla(X, Y):=[s_\nabla X, s_\nabla 
                                Y]_{J^1{\mathfrak g}}-s_\nabla[ X, Y]_{\mathfrak g}.$
                        
                                \item\lab{kp} $\nabla$ is a Cartan connection if and 
                                only if $\cocurvature\nabla=0$.
                
              \item\lab{determine} For any sections $ X, Y, Z\subset{\mathfrak g}$
                and $V\subset TM$ one has
                \begin{align*}
                  \cocurvature\nabla(X,Y)\#Z&=-\curvature\bar\nabla (X, Y)Z,\\
                  \#\cocurvature\nabla(X,Y)V&=-\curvature\bar\nabla(\#X,\#Y)V,
                \end{align*}
                where $\bar\nabla$ denotes the  associated ${\mathfrak
                  g}$-connection on ${\mathfrak g}$ in the first
                formula, and on $TM$ in the second.
                
                \item\lab{costly} In particular, if ${\mathfrak g}=TM$, then
                \begin{equation*}
                        \cocurvature\nabla=-\curvature\bar\nabla,
                \end{equation*}
                                where $\bar\nabla$ denotes the dual linear connection 
                                on $TM$.
        \end{conditions}
\end{proposition}

As simple consequences of  \eqref{determine} we have:
\begin{corollary}\mbox{}
  \begin{conditions}
        \item\lab{sort1} Suppose ${\mathfrak g}$ is transitive.  Then 
        $\nabla$ is a Cartan connection on ${\mathfrak g}$ if and only if 
        the associated ${\mathfrak g}$-connection $\bar\nabla$ on 
        ${\mathfrak g}$ is flat.
    
        \item\lab{sort2} Suppose ${\mathfrak g}$ has an injective anchor.  
        Then $\nabla$ is a Cartan connection if and only if the associated 
        ${\mathfrak g}$-connection $\bar\nabla$ on $ TM$ is flat.
  \end{conditions}
\end{corollary}
Although we shall make no use of the fact here, it is worth remarking 
that a Cartan connection $\nabla$ on a transitive Lie algebroid 
${\mathfrak g}$ is {\em uniquely} determined by the corresponding 
self-representation $\bar\nabla$; see \cite[Proposition 
6.1]{Blaom_05}.

\subsection{Basic examples of Cartan algebroids}\lab{basics}\mbox{} 
We now list some elementary examples of Cartan algebroids.  Example 
\eqref{hen} explains the choice of name `Cartan algebroid.'
\begin{conditions}
\item Every action algebroid $\mathfrak g_0\times M$, equipped with 
its canonical flat connection $\nabla$, is a Cartan algebroid.  
Locally this is the only flat example.  See \ref{dissymmetry} below.

\item As we sketch in Appendix \ref{global}, every Lie pseudogroup of 
transformations in $M$ has a flat Cartan algebroid as its 
infinitesimalization.

\item\lab{placid}  According to Proposition \eqrefs{cocurvature}{costly}  a 
linear connection $\nabla$ on $ TM$ is a Cartan connection if and only 
if its dual $\nabla^*$ is flat, i.e., is an infinitesimal parallelism 
on $ M$. By duality, every Cartan connection on $ TM$  arises as the 
dual of some infinitesimal parallelism. See also  \ref{parallelisms}.

\item\lab{pointFive} If $ M$ is a Lie group, then the flat linear 
connection $\nabla$ on $ TM$ corresponding to left (or right) 
trivialization of $ TM$ is, as a special case of \eqref{placid}, a 
Cartan connection on $ TM$.

\item\lab{goose} The distribution $D$ tangent to some regular
  foliation ${\mathcal F}$ on $ M$ is a subalgebroid of $ TM$.
  According to Corollary \eqrefs{cocurvature}{sort2}, a linear
  connection $\nabla$ on $ D$ is a Cartan connection if and only if
  the $D$-connection $\bar\nabla$ on $TM$, defined by
  $\bar\nabla_UV=\nabla_VU+[U, V]$, $U\subset D,\, V\subset TM$, is
  flat.

\item In \ref{linear2} elementary arguments will show that every 
{\em torsion-free} linear connection $\nabla$ on $TM$ determines a Cartan 
connection on $J^1 (TM)$.  The parallel sections of this Cartan 
connection are the prolonged infinitesimal isometries of $\nabla$.

\item\lab{hen} Let $M$ be a classical Cartan geometry modeled on some
  homogeneous space $ G_0/H_0$ (see, e.g., \cite{Sharpe_97}).  If
  $\pi\colon P\rightarrow M$ denotes the associated principal
  $H_0$-bundle, then according to \cite{Blaom_05}, the vector bundle
  $TP/H_0$ is a Lie algebroid supporting a Cartan connection
  determined by the classical Cartan connection on $P$. (Note,
  however, that flatness of the former is not sufficient for flatness
  of the latter; see \cite{Blaom_05} and our comments on curvature in
  \ref{cart}.)
\end{conditions}
\subsection{The symmetric part of a Cartan algebroid}%
\lab{substructure}%
An arbitrary Cartan algebroid $\mathfrak g $ has a canonical 
subalgebroid isomorphic to an action Lie algebroid.  Indeed, let 
$\nabla$ denote the Cartan connection and let ${\mathfrak 
g}_0\subset\Gamma(\mathfrak g)$ be the subspace of $\nabla$-parallel 
sections, which is finite-dimensional.  Then vanishing cocurvature 
ensures that ${\mathfrak g}_0\subset\Gamma(\mathfrak g)$ is a Lie 
subalgebra, and we obtain an action of $\mathfrak g_0$ on $M$ given by
\begin{gather*}
  \mathfrak g_0\times M\rightarrow TM\\
  (X,m)\mapsto \# X(m).  
\end{gather*}
Equipping the action algebroid $\mathfrak g_0\times M$ with  its 
canonical flat connection, we obtain a morphism of Cartan algebroids,
\begin{gather}
  \mathfrak g_0\times M\rightarrow\mathfrak g\label{kernel}\\
  (X,m)\mapsto X(m).\notag
\end{gather}
Assuming $M$ is connected, this morphism is injective because 
$\nabla$-parallel sections vanishing at a point vanish everywhere.  
We call the image of the monomorphism \eqref{kernel} the {\df symmetric part} of $\mathfrak g $.

\subsection{Curvature as the local obstruction to symmetry}\lab{dissymmetry}%
A Cartan algebroid $\mathfrak g$ is {\df globally flat} if it is
isomorphic to an action algebroid ${\mathfrak g}_0\times M$, equipped
with its canonical flat connection --- or, equivalently, if it
coincides with its symmetric part.  We call $\mathfrak g$ {\df flat}
if every point of $M$ has an open neighborhood $U$ on which the
restriction $\mathfrak g|_U$ is globally flat.\footnote{In
  \cite{Blaom_05} we used {\df symmetric} and {\df locally symmetric}
  in place of {\df globally flat} and {\df flat}, respectively.}

The following theorem shows that a Cartan algebroid may be viewed as 
an infinitesimal symmetry deformed by curvature.
\begin{theorem}[\cite{Blaom_05}]
  Let $\mathfrak g $ be a Cartan algebroid with Cartan connection
  $\nabla$, defined over a connected manifold $M$. Then $\mathfrak g$
  is flat if and only if $\curvature\nabla=0$.  When $M$ is
  simply-connected,  flatness already implies  global flatness.
  
  In the globally flat case the bracket on the Lie algebra ${\mathfrak 
  g}_0$ of $\nabla$-parallel sections is given by
  \begin{equation}
    [\xi,\eta]_{{\mathfrak g}_0}=\torsion\bar\nabla 
    (\xi,\eta),\mathlab{bracket}
  \end{equation}
  where $\bar\nabla$ denotes the  associated representation of ${\mathfrak 
  g}$ on itself: $$\bar\nabla_XY=\nabla_{\#Y}X+[X, Y].$$
\end{theorem}
\begin{proof}
  The necessity of vanishing curvature is immediate.  To establish the 
  assertions in the first paragraph it suffices to show that 
  \eqrefs{substructure}{kernel} is an isomorphism whenever 
  $\curvature\nabla=0$ and $M$ is simply-connected.  Indeed, in that 
  case $\nabla$ determines a trivialization of the bundle $\mathfrak 
  g$ in which constant sections correspond to the $\nabla$-parallel 
  sections of $\mathfrak g$ --- that is, to elements of $\mathfrak 
  g_0$.  In particular, $\mathfrak g_0\times M$ and $\mathfrak g$ will 
  have the same rank, implying the monomorphism 
  \eqrefs{substructure}{kernel} is an isomorphism.
  
  The formula \eqref{bracket} holds in the globally flat case because
  it holds for any action algebroid, as is readily established.
\end{proof}
\begin{example}
  Every Lie group possesses a dual pair of flat linear connections 
  $\nabla,\nabla^*$ corresponding to the left and right 
  trivializations of the tangent bundle (see 
  \eqrefs{basics}{pointFive} above).  Conversely, whenever a 
  simply-connected manifold $M$ supports a linear connection $\nabla$ 
  on $TM$ such that $\nabla$ and its dual $\nabla^*$ are simultaneously flat, 
  then $\nabla$ is a flat Cartan connection on $ TM$ and the theorem 
  above delivers an isomorphism $ TM\cong{\mathfrak g}_0\times M$, 
  where ${\mathfrak g}_0$ is the Lie algebra of $\nabla$-parallel 
  vector fields.  This isomorphism amounts to a ${\mathfrak 
  g}_0$-valued Mauer-Cartan form on $M$, integrating to a Lie group 
  structure under a suitable completeness hypothesis.  The Lie bracket 
  on $\mathfrak g_0$ is given by $[U,V]=-\torsion\nabla (U, V)$.
\end{example}
\noindent For the application of Theorem \ref{dissymmetry} to examples 
\eqrefs{basics}{goose} and \eqrefs{basics}{hen}, see \cite{Blaom_05}.

\section{Examples of infinitesimal geometric structures}\lab{sectionQ}
In this section we describe the infinitesimal geometric structures
associated with Riemannian structures, vector fields on a Riemannian
manifold, parallelisms, almost complex structures, Poisson structures,
subgeometries of an Abelian Lie group, affine structures, projective
structures, classical $G$-structures, and Cartan
algebroids. 

Subriemannian contact structures and conformal structures are
described separately in Sections\,\ref{contact} and
\ref{point}. Conformal parallelisms are described in
\ref{conformalParallelism}.
%
% This section describes several examples of infinitesimal geometric
% structures, including those associated with Riemannian structures,
% almost complex structures, Poisson structures, affine structures,
% projective structures, arbitrary $G$-structures, and Cartan
% algebroids.  From our description of affine structures, it will be
% clear how one may associate an infinitesimal geometric structure with
% an arbitrary (but suitably regular) differential operator on $M$.
% Since every classical Cartan geometry has an associated Cartan
% algebroid (see \eqrefs{basics}{hen} above) every such geometry has an
% associated infinitesimal geometric structure as well.  

% Examples \ref{vectorField}  and \ref{}??of a non-surjective infinitesimal geometric structure is given
% in \ref{vectorField}.

% Conformal structures and subriemannian contact structures are 
% described separately in Sections\,\ref{contact} and \ref{point}.
%
\subsection{Isotropy}\lab{isotropy}
Most infinitesimal geometric structures occurring in nature are best 
understood as isotropy (or {\em joint} isotropy) subalgebroids, of 
certain jet-bundle representations.  In the case of Riemannian 
geometry it is the isotropy of the Riemannian metric 
$\sigma\subset\symmetric^2 (TM)$, i.e., of a {\em section} of some vector 
bundle; in conformal geometry, it is the isotropy of a rank-one {\em 
subbundle} $\langle\sigma\rangle\subset\symmetric^2 (TM)$; in 
projective geometry, it is the isotropy of an {\em affine} subbundle 
of $ J^1 (TM)^*\otimes T^*\!M\otimes TM$.  The following definition of 
isotropy is general enough to cover all these possibilities.

Let $\rho\colon{\mathfrak g}\rightarrow\gl(E)$ denote some
representation of a Lie algebroid ${\mathfrak g}$.  Let $\Sigma\subset
E$ denote any affine subbundle of $ E$ (a single section of $E$ in the
simplest case) and $\Sigma_0\subset E$ the corresponding vector
subbundle parallel to $\Sigma$ (resp., the zero section). Assume
either that $\Sigma_0$ is $\mathfrak g$-invariant or that
$\Sigma=\Sigma_0$. Then the {\df isotropy} of $\Sigma$ is the
collection of all elements $x\in{\mathfrak g}$ for which
\begin{equation*}
        \sigma\subset\Sigma\implies\rho_x\sigma\in\Sigma_0,
\end{equation*}
for arbitrary local sections $\sigma\subset E$; here
$\rho_x\sigma:=\rho(x)(J^1\sigma(m))$ where $ m\in M$ is the base
point of $ x$.

The isotropy of $\Sigma$ is a subset of ${\mathfrak g}$ intersecting 
fibers in subspaces whose dimensions may vary, i.e., is a 
`variable-rank subbundle'; sections of this bundle are closed under 
the bracket of ${\mathfrak g}$.  When this rank is constant the 
isotropy is a bona fide subbundle and consequently a subalgebroid, 
called the {\df isotropy subalgebroid} of $\Sigma$.  A section $X$ of 
${\mathfrak g}$ is then a section of the isotropy if and only if 
$\rho_X\sigma$ is a section of $\Sigma_0$ for all sections 
$\sigma\subset\Sigma$.

% \subsection{On structure kernels}
% Let ${\mathfrak h}$ denote its structure kernel of an infinitesimal
% geometric structure ${\mathfrak g}\subset J^1{\mathfrak t}$.  Since
% $J^1\mathfrak t\rightarrow\mathfrak t$ has kernel
% $T^*\!M\otimes\mathfrak t$, we have ${\mathfrak h} \subset T^*\!M
% \otimes {\mathfrak t}$.  If ${\mathfrak h}$ has constant rank, then
% ${\mathfrak h}\subset T^*\!M\otimes\mathfrak t$ is a subalgebroid.
% Since $T^*\!M\otimes\mathfrak t$ is totally intransitive, this simply
% means ${\mathfrak h}(m)$ is a subalgebra of $T_m^*M\otimes{\mathfrak
%   t} (m)$ at each point $ m\in M$, which is also true in the variable
% rank case.  In classical parlance, each fiber of ${\mathfrak h}$ is a
% {\df tableau} \cite{Bryant_etal_91}.  Recall from Proposition
% \ref{garden} that the bracket on $T^*\!M\otimes{\mathfrak t}$ is given
% by
% \begin{equation*}
%                 [\phi_1,\phi_2]_{T^*\!M\otimes{\mathfrak t}}(U)=\phi_1(\# 
%                 \phi_2 U)-\phi_2(\# \phi_1 U);\qquad U\subset TM.
% \end{equation*}

\subsection{Riemannian structures}\lab{Riemannian}
The adjoint representation of $ J^1 (TM)$ on $ TM$ determines a 
representation of $ J^1 (TM)$ on $\symmetric^2 (TM)$: A section 
$X\subset J^1 (TM) $ acts on a section $\sigma\subset\symmetric^2 
(TM)$ according to
\begin{equation}
        (X\cdot\sigma) (V_1, V_2):={\mathcal L}_{\#X}(\sigma (V_1, 
         V_2))-\sigma (\adjoint_X^{TM}V_1, V_2)-\sigma 
         (V_1,\adjoint_X^{TM}V_2).\mathlab{tsk}
\end{equation}
Since this simply means $ J^1 V\cdot\sigma={\mathcal L}_V\sigma$, the 
isotropy ${\mathfrak g}\subset J^1 (TM)$ of a Riemannian metric 
$\sigma\subset\symmetric^2 (TM)$ is indeed its associated bundle of 
$1$-symmetries, as described in \ref{firsty}.

According to  \eqref{tsk}, the subalgebroid $T^*\!M\otimes TM\subset 
J^1 (TM)$ acts on $\symmetric^2 (TM)$ via
\begin{equation*}
         (\phi\cdot\sigma) (V_1, V_2)=-\sigma (\phi V_1, V_2)-\sigma 
         (V_1,\phi V_2);\qquad\phi\subset T^*\!M\otimes TM.
\end{equation*}
The structure kernel of ${\mathfrak g}$ is the isotropy ${\mathfrak 
h}\subset T^*\!M\otimes TM$ of $\sigma$ under this restricted 
representation.  So ${\mathfrak h}$ is the bundle of 
$\sigma$-skew-symmetric endomorphisms of tangent spaces, a Lie algebra 
bundle modeled on ${\mathfrak o}(n)$, $n:=\dimension M$.  In 
particular, all conformally equivalent metrics give the same structure 
kernel.

One way to see that ${\mathfrak g}$ is surjective (i.e., transitive)
is to apply the algebraic lemma \ref{algebraicLemma} in Appendix
\ref{miscellany} to the morphism $ X\mapsto X\cdot\sigma\colon J^1
(TM)\rightarrow\symmetric^2 (TM)$, whose kernel is ${\mathfrak g}$.
On account of the surjectivity of the restriction
$\phi\mapsto\phi\cdot\sigma\colon T^*\!M\otimes
TM\rightarrow\symmetric^2 (TM)$ of this morphism, the lemma delivers
an exact sequence
\begin{equation}
    0\rightarrow{\mathfrak h}\rightarrow{\mathfrak g}\rightarrow 
    TM\rightarrow 0.\mathlab{church}
\end{equation}
Thus ${\mathfrak g}\subset J^1(TM)$ is surjective and has 
constant rank (making it a subalgebroid and thus an infinitesimal 
geometric structure).

The lemma just applied is very useful in determining the image and
structure kernel of infinitesimal geometric structures defined by
isotropy.  Less trivial applications of Lemma \ref{algebraicLemma} are
made in \ref{vectorField} and \ref{almostComplex} below.  Minimal
comment will accompany subsequent applications.

The symmetries of ${\mathfrak g}$ (in the sense of \ref{structure})
are the vector fields along which $\sigma$ has vanishing Lie
derivative, i.e., its Killing fields.  A generator of $\mathfrak g$
(see \ref{neoGenerator}) is a linear connection $\nabla$ on $TM$ such
that $\sigma$ is $\bar\nabla$-parallel, where $\bar\nabla$ denotes the
dual of $\nabla$.  The Levi-Cevita connection is thus the unique
torsion-free generator of ${\mathfrak g}$.

From ${\mathfrak g}$ one can recover the metric $\sigma$ up to a 
positive {\em constant} (not merely its conformal class).  In the 
simply-connected case, slightly more is true:
\begin{proposition}
  Let ${\mathfrak h}\subset T^*\!M\otimes TM$ denote the ${\mathfrak
    o}(n)$-bundle associated with an arbitrary conformal structure.
  Then on simply-connected open subsets of $ M$, every surjective
  infinitesimal geometric structure ${\mathfrak g}\subset J^1(TM)$
  having structure kernel ${\mathfrak h}$ is the isotropy subalgebroid
  of some Riemannian structure $\sigma$ in the given conformal class.
  This structure is uniquely determined up to a constant.
\end{proposition}
\begin{proof}
         Suppose ${\mathfrak g}\subset J^1(TM)$ has structure kernel 
         ${\mathfrak h}$ and let $ L\subset\symmetric^2 (TM)$ be the line 
         bundle determined by the conformal structure. That is, $L$ is the 
         bundle of ${\mathfrak h}$-invariant elements of $\symmetric^2 (TM)$. 
         The non-vanishing elements of  $L$ are either positive or negative 
         definite. By  Lemma \ref{extensionLemma}, $ L$ has a 
         non-vanishing ${\mathfrak g}$-invariant section $\sigma$, unique up 
         to constant. Changing the sign of $\sigma$ if necessary, we obtain 
         the sought after metric.
\end{proof}

The application of Cartan's method to Riemannian structures is given in
\ref{illustration}.

In analogy with the Riemannian case, the isotropy subset ${\mathfrak
  g}\subset J^1 (TM)$ of an arbitrary tensor on $M$ is an
infinitesimal geometric structure, whenever this isotropy has constant
rank.  Moreover, in many cases, this structure encodes all useful
information (i.e., some analogue of the preceding proposition
applies.)

For structures defined by more than one tensor one considers the {\em
  joint} isotropy, defined by intersecting the individual isotropies.
For example, for almost K\"ahler structures, one considers the joint
isotropy of the symplectic and almost complex structures. Here is
another example:

\subsection{Vector fields on a Riemannian manifold}\lab{vectorField}
Let ${V}$ be a {\em nowhere vanishing} vector field on a Riemannian 
manifold $M$ with metric $\sigma$.  The vector fields on $M$ that are 
simultaneously infinitesimal isometries of $\sigma\subset\symmetric^2 
(TM)$ and ${V}\subset TM$ are the symmetries of the joint isotropy of 
$\sigma$ and ${V}$, with respect to representations of $ J^1(TM)$ on 
$\symmetric^2 (TM)$ and $ TM$ respectively.  Denoting the isotropy of 
$\sigma$ alone by ${\mathfrak g}\subset J^1(TM)$ as above, ${\mathfrak 
g}$ acts on $TM$ by restricted adjoint action and the joint isotropy 
is the isotropy ${\mathfrak g}_{V}\subset{\mathfrak g}$ of ${V}$.

The structure kernel of ${\mathfrak g}_{V}$ is the ${\mathfrak 
o}(n-1)$-bundle of skew-symmetric tangent space endomorphisms 
infinitesimally fixing $ V$ (mapping ${V}$ to $ 0$).
\begin{proposition}
         The image of ${\mathfrak g}_{V}$ is the  distribution $ D\subset TM$ 
         tangent to the level sets of $\frac{1}{2}\left\|V\right\|^2$.
\end{proposition}
\noindent In particular, ${\mathfrak g}$ has constant rank (is an
infinitesimal geometric structure) if and only if $V$ has constant
length or $\frac{1}{2}\left\|V\right\|^2$ is a free of critical
points; ${\mathfrak g}_V$ is surjective (transitive) in the former
case only.
\begin{proof}
Applying Lemma \ref{algebraicLemma} to the morphism $X\mapsto X\cdot 
V\colon{\mathfrak g}\rightarrow TM$, of which ${\mathfrak g}_V$ is the 
kernel, we deduce that $D$ is the kernel of the morphism $\Theta$ 
making the following diagram commute:
\begin{equation*}
            \begin{CD}
                        {\mathfrak g} @>\#>> TM\\
                         @V{X\mapsto X\cdot V}VV @VV\Theta V \\
                         TM@>>>  TM/V^\perp
            \end{CD}\enspace.
\end{equation*}
Let $\nabla$ be any generator of ${\mathfrak g}$ (e.g., the
Levi-Cevita connection) and $ s\colon TM\rightarrow{\mathfrak g}$ the
corresponding splitting of \eqref{church} above.  Then $\Theta
(U)=sU\cdot V\modulo V^\perp=\bar\nabla_UV\modulo V^\perp$, where
$\bar\nabla$ is the dual connection.  Or, identifying $ TM/V^\perp$
with the trivial line bundle ${\mathbb R}\times M$, using $V$, we have
$\Theta (U)=\sigma (V,\bar\nabla_UV)=\frac{1}{2}{\mathcal L}_U(\sigma
(V, V))=d(\frac{1}{2}\left\|V\right\|^2) (U)$.  Here we have used
$\bar\nabla\sigma=0$, which holds because $\nabla$ is a generator.
\end{proof}

\subsection{Parallelism}\lab{parallelisms}
The simplest non-trivial example of an infinitesimal geometric 
structure is a transitive infinitesimal geometric structure 
${\mathfrak g}\subset J^1 (TM)$ having trivial structure kernel.  In 
other words, ${\mathfrak g}$ is a subalgebroid of $J^1 (TM)$ mapped 
isomorphically onto $TM$ by the anchor $\#\colon J^1 (TM)\rightarrow 
TM$.  According to Theorem \ref{neoGenerator}, ${\mathfrak g}$ has a 
unique generator $\nabla$ that is a Cartan connection on $TM$. From 
our observations in \eqrefs{basics}{placid}, the dual connection $ 
\bar\nabla $ is flat, i.e., an infinitesimal parallelism on $M$, and 
conversely all infinitesimal parallelisms arise in this way.

When $ M$ is simply-connected the Lie algebroid morphism 
$\bar\nabla\colon TM\rightarrow{\mathfrak{gl}}(TM)$ integrates to a 
Lie groupoid morphism $ M\times M\rightarrow\operatorname{GL}(TM)$, 
i.e., to an {\em absolute} parallelism on $ M$ (a trivialization of 
the tangent bundle).  

% The symmetries of ${\mathfrak g}$ are the $\nabla$-parallel vector
% fields. 
When $\bar\nabla$ comes from an absolute parallelism, viewed
instead as some non-degenerate $V$-valued 1-form $\omega$ on $M$
($\dimension V=\dimension M$), then ${\mathfrak g}\subset J^1 (TM)$ is
the isotropy subalgebroid of $\omega\subset T^*\!M \otimes  V$, and a
symmetry of ${\mathfrak g}$ is a vector field along which $\omega$ has
vanishing Lie derivative.

The infinitesimal geometric structure associated with a {\em conformal} parallelism is described in \ref{conformalParallelism}.

\subsection{Almost complex structures}\lab{almostComplex}
The infinitesimal geometric structure associated with an almost complex structure ${\mathbf J}$  on $M$ is generically {\em intransitive}. Indeed, let ${\mathfrak g} \subset J^1(TM) $ denote the isotropy of $\mathbf J \subset T^*\!M \otimes TM$ under the representation of $J^1(TM) $ on $T^*\!M \otimes TM$ induced by the adjoint representation, and define
\begin{equation*}
  [T^*\!M \otimes TM, {\mathbf J}]:=\cup_{m\in M}\{\phi {\mathbf J}(m)-{\mathbf J}(m)\phi \suchthat \phi \in T_m^*M \otimes T_mM\};
\end{equation*}
let $N \subset \alternating^2(TM)\otimes TM$ denote the Nijenhuis
torsion of ${\mathbf J}$:
\begin{equation*}
  N(U,V)=\frac{1}{4}\Big(\,[{\mathbf J} U, {\mathbf J} V]-[U,V]-{\mathbf J}[U,{\mathbf J} V]-{\mathbf J}[{\mathbf J} U, V]\,\Big).
\end{equation*}
Then:
\begin{proposition}
  The structure kernel of $\mathfrak g$ is $T^*\!M \otimes_{\mathbb C} TM$ and the image of $\mathfrak g$ is the kernel of the morphism
  \begin{gather*}
    \Theta \colon TM \rightarrow (T^*\!M \otimes TM)/[T^*\!M \otimes TM, {\mathbf J}] \\
    \Theta(U)=-4N({\mathbf J} U,\,\cdot,\,) \modulo [T^*\!M \otimes TM, {\mathbf J}].
  \end{gather*}
\end{proposition}
\noindent In particular, ${\mathfrak g} $ is transitive if and only if the section $N(U,\,\cdot\,) \subset T^*\!M \otimes TM $ lies in $[T^*\!M \otimes TM, {\mathbf J}]$ for all vector fields $U$.
\begin{proof}
  First, note that
  \begin{equation}
    \frac{1}{2}\Big[\adjoint_{J^1({\mathbf J} U)}{\mathbf J},{\mathbf J}\Big]V=-{\mathbf J}[{\mathbf J} U,{\mathbf J} V]-[{\mathbf J} U, V];\qquad U,V \subset TM.\mathlab{satori}
  \end{equation}
Next observe that ${\mathfrak g} $ is the kernel of the morphism
\begin{gather*}
  \theta \colon J^1(TM) \rightarrow T^*\!M \otimes  TM \\
  \theta(X)=\adjoint_X {\mathbf J},\\
\text{i.e., }\enspace \theta(X)U=\adjoint_X({\mathbf J} U)-{\mathbf J}(\adjoint_X U).
\end{gather*}
Applying Lemma \ref{algebraicLemma} to this morphism, we obtain a morphism $\Theta $ with the requisite kernel, and satisfying
\begin{align*}
  \Theta(U)&=\adjoint_{J^1 U}{\mathbf J} \enspace\modulo [T^*\!M \otimes TM, {\mathbf J}]\\
  &=\adjoint_{J^1 U}{\mathbf J} -  \frac{1}{2}\Big[\adjoint_{J^1({\mathbf J} U)}{\mathbf J},{\mathbf J}\Big]\enspace\modulo [T^*\!M \otimes TM, {\mathbf J}].
\end{align*}
By \eqref{satori}, we have
\begin{multline*}
 \left(\, \adjoint_{J^1 U}{\mathbf J} - \frac{1}{2}\Big[\adjoint_{J^1({\mathbf J} U)}{\mathbf J},{\mathbf J}\Big] \,\right)V
=[U,{\mathbf J} V]-{\mathbf J} [U,V]\\ + {\mathbf J}[{\mathbf J} U, {\mathbf J} V]+[{\mathbf J} U, V]=-4N({\mathbf J} U, V).
\end{multline*}
\end{proof}

\subsection{Poisson structures}\lab{Poisson}
Although not of finite-type, Poisson structures furnish us with
another interesting example of an infinitesimal geometric
structure. Now generally the isotropy $\mathfrak g\subset J^1(TM)$ of
a Poisson tensor on $M$ fails to have constant rank, and so fails to
be an infinitesimal geometric structure on $TM$. However, one can
define an infinitesimal geometric structure on the {\em cotangent}
bundle $T^*\!M$, which the Poisson tensor makes into a Lie algebroid
(see below). Although not transitive, this structure {\em is}
surjective.

Let $\omega$ be a symplectic structure on $M$ and let $\#\colon 
T^*\!M\rightarrow TM$ denote the inverse of 
$v\mapsto\omega(v,\,\cdot\,)$.  Since $\#$ is an isomorphism, there is 
a unique bracket on $\Gamma(T^*\!M)$ making $T^*\!M$ into a Lie algebroid 
with anchor $\#$.  This bracket is given by
\begin{equation}
  [\alpha,\beta]_{T^*\!M}=\mathcal L_{\#\alpha}\beta-\mathcal L_{\#\beta}\alpha 
                       + d(\Pi(\alpha,\beta)),\qquad\alpha,\beta\in\Gamma(T^*\!M),
  \mathlab{one}
\end{equation}
where $\mathcal L$ denotes Lie derivative and $\Pi$ is the Poisson
tensor. This tensor is defined by
$\Pi(\alpha,\beta):=\omega(\#\alpha,\#\beta)$ and so satisfies
\begin{equation}
  \langle\alpha,\#\beta\rangle=\Pi(\alpha,\beta)
  \qquad\alpha,\beta\in\Gamma(T^*\!M).\mathlab{two}
\end{equation}
More generally, \eqref{one} defines a Lie algebroid structure on
$T^*\!M$ for {\em any} Poisson manifold $(M,\Pi)$, with anchor $\#$
defined by \eqref{two}. The symplectic leaves of $\Pi$ are precisely
the orbits of the Lie algebroid $T^*\!M$.

An infinitesimal isometry of a Poisson manifold $(M,\Pi)$ is a vector
field $V$ on $M$ such that $\mathcal L_V\Pi=0$. Poisson manifolds have
an abundance of infinitesimal isometries. In particular, every closed
1-form $\alpha$ on $M$ determines an infinitesimal isometry $\#\alpha$
tangent to the symplectic leaves known as a {\df local Hamiltonian
  vector field}, or a {\df Hamiltonian vector field} if $\alpha$ is
exact.

It is not too difficult to establish the following result; see
\cite{Blaom_05} for some details:
\begin{proposition}
  Let $\mathfrak g\subset J^1(T^*\!M)$ denote the kernel of the vector 
  bundle morphism $J^1(T^*\!M)\rightarrow\alternating^2(TM)$ whose 
  corresponding map on sections sends $J^1\alpha$ to $d\alpha$.  Then 
  $\mathfrak g$ is a surjective infinitesimal geometric structure on 
  $T^*\!M$, with structure kernel $\symmetric^2(TM)$, whose symmetries 
  are the closed one-forms on $M$.
  
  A linear connection $\nabla$ on $T^*\!M$ is a generator of 
  $\mathfrak g$ if and only if the corresponding linear connection on 
  $TM$ is torsion free. Such a generator is a Cartan connection on 
  $T^*\!M$ if and only if
  \begin{equation*}
    \curvature\nabla\, (V,\#\alpha)\beta
        -\curvature\nabla\, (V,\#\beta)\alpha-(\nabla_V(\nabla\Pi)) 
        (\alpha,\beta)=0,
  \end{equation*}
  for all sections $\alpha,\beta\subset T^*\!M; V\subset TM$.
\end{proposition}
\noindent If $M$ is the dual of a Lie algebra, equipped with its 
Lie-Poisson structure (see, e.g., \cite[\S 10.1]{Marsden_Ratiu_94}),
then the canonical flat linear connection $\nabla$ on $T^*\!M\cong
M\times M^*$ is an example of a Cartan connection as described in the
proposition.  Up to certain momentum map equivariance obstructions,
this is locally the only flat example \cite[Corollary 3.4]{Blaom_05}.

\subsection{Subgeometries of an Abelian Lie group}
We now give a simple example of a geometric structure which has, in
general, trivial isotropy.  Let $E$ be an Abelian Lie group, $V$ its
Lie algebra, and $M\subset E$ a codimension-one submanifold. Let $\omega$ be the $V$-valued one-form on $M$ obtained by restricting the Mauer-Cartan form 
on $E$. Then $d\omega=0$ and $\dimension V
= \dimension M +1$.

Let $\omega(TM)$ denote the tangent bundle of M, viewed as a subbundle
of $V\times M$, so that $N\equiv (V \times M)/\omega(TM)$ is a model
of the rank-one normal bundle of $M$ in $E$. A basic invariant of $M$
associated with its embedding into $E$ is the quadratic form $\ii
\subset \symmetric^2(TM) \otimes N$ defined by $$\ii(U_1,U_2):=(\nabla
\omega)_\sym(U_1,U_2) \modulo \omega(TM),$$ where $(\nabla
\omega)_\sym(U_1,U_2):=(\nabla_{U_1}
\omega)U_2+(\nabla_{U_2}\omega)U_1$. Here $\nabla $ is an arbitrary
linear connection on $M$, the particular choice not effecting the
definition of $\ii$.
\begin{proposition}
  The isotropy ${\mathfrak g} \subset J^1(TM) $ of $\omega$ has
  trivial structure kernel, and has image $\kernel \ii$. In
  particular, if the quadratic form $\ii$ is non-degenerate, then
  ${\mathfrak g}=0$ and $\omega$ has no infinitesimal isometries.
\end{proposition}
\begin{proof}
  Apply Lemma \ref{algebraicLemma} to the morphism $\theta \colon J^1(TM) \rightarrow T^*\!M \otimes V$; $\theta(X):=X \cdot \omega $. Here $(X\cdot \omega)(U)={\mathcal L}_{\#X}(\omega(U))-\omega (\adjoint_XU)$. Then the morphism $\Theta$ in that lemma is given by 
\begin{align*}
 \Theta (U)&={\mathcal L}_U \omega \modulo \omega(TM)\\
           &= \frac{1}{2}(\nabla \omega)_\sym(U,\,\cdot\,)\modulo \omega(TM),\enspace\text{because of the identity}\\
({\mathcal L}_{U_1}\omega)(U_2)&=\omega(\nabla_{U_2}U_1)+\frac{1}{2}d \omega (U_1,U_2)+\frac{1}{2}(\nabla \omega)_\sym(U_1,U_2).
\end{align*}
\end{proof}

\subsection{Affine structures}\lab{linear2}
Any suitably non-degenerate, $k$th-order, linear, differential 
operator on $M$, defines an infinitesimal geometric structure 
${\mathfrak g}\subset J^{k+1}(TM)\subset J^1(J^k(TM))$.  As a simple 
example, which will suffice to illustrate the general principle, we 
consider an affine structure on $M$, i.e., an arbitrary linear 
connection $\nabla$ on $M$, in which case $k=1$.  The relevant 
non-degeneracy condition is that the isotropy of the torsion of 
$\nabla$ should have constant rank; see below.

View an affine structure $\nabla$ as a section of $ J^1(TM)^*\otimes 
T^*\!M\otimes TM$, via
\begin{equation*}
        \nabla (J^1 W, V):=\nabla_VW;\qquad V, W\subset TM.
\end{equation*}
In order to associate a natural isotropy subalgebroid with $\nabla$, 
we  begin with two observations. First, $ J^1 (J^1(TM))$ acts on $ 
J^1(TM)^*\otimes  T^*\!M\otimes TM$ because $ J^1 (J^1(TM))$ acts 
on $  J^1(TM)$ via adjoint action, and on $ TM$ via the composite
\begin{gather*}
                 J^1(J^1 (TM))\xrightarrow{ p}{ J^1 
                 (TM)}\xrightarrow{\adjoint^{TM}}\gl (TM),\\
                 \text{i.e.,}\quad J^1 X\cdot W=\adjoint_{pX}^{TM} W;\qquad
                 X\subset J^1 (TM),\,W\subset TM.
\end{gather*}
Secondly, $ J^2 (TM)$ may be identified with a subalgebroid of $ J^1 
(J^1(TM))$ via the canonical embedding $J^2(TM)\monomorphism 
J^1(J^1(TM))$ whose corresponding map on sections sends $J^2 V$ to 
$J^1(J^1 V)$.  Combining the two observations, we obtain a natural 
action of $ J^2 (TM)$ on $J^1 (TM)^*\otimes T^*\!M\otimes TM$.
\begin{proposition}
         Let ${\mathfrak g}\subset J^2(TM)$ denote the isotropy of 
         $\nabla\subset J^1(TM)^*\otimes T^*\!M\otimes TM$, and ${\mathfrak 
         t}\subset J^1 (TM)$ the isotropy of 
         $\torsion\nabla\subset\alternating^2 (TM)\otimes TM$. Then:
         \begin{conditions}
                \item\lab{aff5} The symmetries of ${\mathfrak g}$ are the prolonged 
                infinitesimal isometries of  $\nabla$.
                
                \item\lab{frex} The image of ${\mathfrak g}\subset J^1 (J^1(TM))$ is 
                ${\mathfrak t}$ and ${\mathfrak g}$ has trivial structure kernel.
         \end{conditions}
\end{proposition}
\noindent In particular, \eqref{frex} implies that ${\mathfrak 
g}\subset J^2 (TM)$ has constant rank (and is therefore an 
infinitesimal geometric structure on $ J^1 (TM)$) if and only if 
${\mathfrak t}\subset J^1 (TM)$  has constant rank.

Now ${\mathfrak t}=J^1(TM)$ if and only if $\nabla$ is torsion-free 
(because $\identity_{TM}$ is a section of $ T^*\!M\otimes TM\subset 
J^1(TM)$).  On account of \eqref{frex}, ${\mathfrak g}$ is surjective if 
and only if $\nabla$ is torsion-free.  Applying Theorem 
\ref{neoGenerator}, we obtain:
\begin{corollary}
        If $\torsion\nabla=0$, then the unique generator $\nabla^{(1)}$ of 
        ${\mathfrak g}$ is a Cartan connection on $J^1 (TM)$ whose 
        parallel sections are the prolonged infinitesimal isometries of 
        $\nabla$.
\end{corollary}
\noindent From an explicit formula for $\nabla^{(1)}$ one may completely 
characterize the obstructions to the existence of infinitesimal 
isometries; see \ref{illustration2}.
\begin{proof}[Proof of proposition]
Recall that a vector field $ U$ on $M$ is an {\df infinitesimal 
isometry} of $\nabla$ if
\begin{equation}
        [U,\nabla_VW]-\nabla_{[U,V]}W-\nabla_V[U,W]=0;
        \qquad V, W\subset TM,\mathlab{aff1}
\end{equation}
a condition that is second-order in $ U$.  Unravelling the definition 
of the representations defined above, we may write this condition as
\begin{equation*}
        J^1 (J^1(U))\cdot\nabla=0.
\end{equation*}
It easily follows that $ J^1 U$ is a symmetry of ${\mathfrak g}$ 
whenever $ U$ is an infinitesimal isometry of $\nabla$.

Suppose, conversely, that $ X\subset J^1 (TM)$ is a symmetry of 
${\mathfrak g}$, i.e., that $ J^1 X$ lies in ${\mathfrak g}$.  This 
means:
\begin{align}
                J^1 X&\subset J^2 (TM),\mathlab{affa}\\
                \text{and}\enspace J^1 X&\cdot\nabla=0.\mathlab{affb}
\end{align}
It is well known that \eqref{affa} is equivalent to $ X\subset J^1
(TM)$ being holonomic (see Lemma \ref{prolongation}).  So $X=J^1 U$,
where $ U$ is an infinitesimal isometry, on account of \eqref{affb},
which reads $ J^1 (J^1 U)\cdot\nabla=0$.  This completes the proof of
\eqref{aff5}.

Let $\xi$ be any section of $J^2 (TM)$.  It is easy to check that the 
section $\xi\cdot\nabla\subset J^1 (TM)^*\otimes T^*\!M\otimes TM$ is 
tensorial, i.e., drops to some section $(\xi\cdot\nabla)^\vee\subset 
T^*\!M\otimes T^*\!M\otimes TM$.  Noting that ${\mathfrak g}$ is then 
the kernel of the morphism
\begin{align*}
                \xi&\mapsto (\xi\cdot\nabla)^\vee\\
                J^2 (TM)&\rightarrow  T^*\!M\otimes T^*\!M\otimes TM,
\end{align*}
whose domain $ J^2 (TM)$ fits into an exact sequence
\begin{equation*}
                 0\rightarrow\symmetric^2 (TM)\otimes TM\monomorphism J^2 
                 (TM)\rightarrow J^1 (TM)\rightarrow 0,
\end{equation*}
one shows, by applying  Lemma \ref{algebraicLemma}, that 
${\mathfrak g}$ fits into a corresponding exact sequence
\begin{equation*}
                 0\rightarrow 0\rightarrow{\mathfrak g}\xrightarrow{b} 
                 {\mathfrak  t}\rightarrow 0.
\end{equation*}
Here $ b$ is the restriction of the canonical projection $ J^2 
(TM)\rightarrow J^1 (TM)$.  This establishes \eqref{frex}.
                \end{proof}

\subsection{Projective structures}\lab{projective}
Recall that two linear connections $\nabla,\nabla'$ are {\df 
projectively equivalent} if their geodesics coincide as 
unparameterized curves.  Equivalently their difference 
$\nabla-\nabla'$, which may be viewed as a section of $$T^*\!M\otimes 
T^*\!M\otimes TM\subset J^1(TM)^*\otimes T^*\!M\otimes TM,$$ should 
take its values in the subbundle $\Sigma_0:=(\alternating^2 
(TM)\otimes TM)\oplus j_{\mathrm S}( T^*\!M)$ of 
\begin{equation*}
  T^*\!M\otimes T^*\!M\otimes TM\cong
  \Big(\,\alternating^2 (TM)\otimes TM\,\Big)\oplus 
  \Big(\,\symmetric^2 (TM)\otimes TM\,\Big).
\end{equation*}
Here $j_\mathrm{S}\colon T^*\!M\rightarrow\symmetric^2 (TM)\otimes 
TM$ is the embedding defined by
\begin{equation*}
                j_\mathrm{ S}(\alpha) (V_1, V_2):=\alpha (V_1)V_2+\alpha 
                (V_2)V_1.
\end{equation*}
A {\df projective structure} is a projective equivalence class of 
linear connections; since $\Sigma_0$ contains $\alternating^2 
(TM)\otimes TM$, each such class has a torsion-free representative 
$\nabla$.

Let $\nabla$ be a torsion-free linear connection on $ TM$ and 
$\langle\nabla\rangle$ the corresponding projective structure.  To 
specify the structure it suffices to specify the affine subbundle
\begin{equation*}
        \nabla+\Sigma_0\subset J^1 (TM)^*\otimes T^*\!M\otimes TM,
\end{equation*}
which we denote by $\langle\nabla\rangle$ also.  As explained in 
\ref{linear2} above, we have $ J^2 (TM)\subset J^1 (J^1 (TM))$ acting 
on $J^1 (TM)^*\otimes T^*\!M\otimes TM$ and can therefore define the 
isotropy ${\mathfrak g}\subset J^2 (TM)$ of $\langle\nabla\rangle$; 
see \ref{isotropy}.  Arguing as in the proof of Lemma 
\eqrefs{linear2}{aff5}, one shows that the symmetries of ${\mathfrak 
g}$ are the prolonged infinitesimal isometries of 
$\langle\nabla\rangle$.

It is not hard to see that ${\mathfrak g}$ has $ j_\mathrm{S}( 
T^*\!M)\cong T^*\!M$  as structure kernel and in fact that
\begin{equation*}
        {\mathfrak g}={\mathfrak g}_\nabla\oplus j_\mathrm{S}( T^*\!M),
\end{equation*}
where ${\mathfrak g}_\nabla\subset J^2 (TM)$ denotes the isotropy of 
$\nabla$ (denoted ${\mathfrak g}$ in \ref{linear2}).  The Cartan 
connection $\nabla^{(1)}$ on $ J^1 (TM)$ in Corollary \ref{linear2} is 
a generator of ${\mathfrak g}_\nabla$ and consequently a generator of 
${\mathfrak g}$ as well.  An explicit formula appears in 
\ref{illustration2}.

\subsection{$G$-structures}\lab{Gstructure}Let $G$ be a subgroup of 
$\operatorname{GL}(n,{{\mathbb R}})$, where $ n$ is the dimension of 
$M$.  A {\df $G$-structure} on $M$ is a $G$-reduction $ P$ of the 
bundle of (absolute) frames on $ M$; see, e.g., \cite{Kobayashi_72}.  
In particular, $P$ is a principal $G$-bundle, so that ${\mathfrak 
g}:=TP/G$ is a transitive Lie algebroid over $M$, and the associated 
vector bundles of $P$ are representations of ${\mathfrak g}$; see, 
e.g., \cite{Mackenzie_05}.  As $P$ is a frame bundle, $TM$ will be 
such a representation (see below).  That is, we have a Lie 
algebroid morphism $${\mathfrak g}\rightarrow\gl 
(TM)\overset{\adjoint}{\cong} J^1 (TM).$$ This turns out to be 
injective, identifying ${\mathfrak g}$ with a subalgebroid of $J^1 
(TM)$.  This infinitesimal geometric structure on $ TM$ is surjective 
because ${\mathfrak g}$ is transitive.

The representation of ${\mathfrak g}$ on $ TM$ may be described as
follows.  Identify sections $X$ of ${\mathfrak g}:=TP/G$ with
$G$-invariant vector fields on $P$, and use the tautological 1-form
on $P$ to identify sections $V$ of $ TM$ with $G$-invariant ${{\mathbb
    R}^n}$-valued functions on $P$.  Then $ X\cdot V:={\mathcal
  L}_XV$, where ${\mathcal L}$ denotes Lie derivative on $P$.

\subsection{Cartan algebroids as infinitesimal geometric structures}%
\lab{cartan}%
We have seen that all surjective infinitesimal geometric structures
with trivial structure kernel define Cartan algebroids (Theorem
\ref{neoGenerator}).  Conversely, if ${\mathfrak t}$ is a Cartan
algebroid with Cartan connection $\nabla$, then ${\mathfrak
  g}:=s_\nabla ({\mathfrak t})\subset J^1{\mathfrak t}$ is a
surjective infinitesimal geometric structure generated by $\nabla $
with trivial structure kernel.  Here $ s_\nabla \colon{\mathfrak t}\rightarrow
J^1{\mathfrak t}$ is the splitting of
\begin{equation*}
                 0\rightarrow T^*\!M\otimes{\mathfrak t}\monomorphism 
                 J^1{\mathfrak t}\rightarrow{\mathfrak t}\rightarrow 0
\end{equation*}
determined by $\nabla$.

\section{Generators, associated operators, and Bianchi identities}\lab{sectionR}
Picking a generator for an infinitesimal geometric structure
${\mathfrak g}\subset J^1{\mathfrak t}$ allows us to identify
${\mathfrak g}$ with the direct sum ${\mathfrak t}_1\oplus{\mathfrak
  h}$ of its image ${\mathfrak t}_1$ and its structure kernel
${\mathfrak h}$. This greatly facilitates computations. Generators are
also the appropriate connections for which to develop all the usual
formalisms of differential geometry: covariant differentiation,
covariant exterior differentiation, Bianchi identities, etc.  (It will
be natural, however, to use the more encompassing descriptor
`associated' in place of `covariant.')  By virtue of \ref{cartan}, we
obtain formalism for Cartan algebroids as a special case.

The present section, rather formal in nature, can be scanned on a
first reading.  In \ref{generator} we address the existence and
uniqueness of generators and prove the theorem in \ref{neoGenerator},
where generators were defined.  In \ref{reconstruction} we see how
information about ${\mathfrak g}$ is encoded in ${\mathfrak t}_1$,
${\mathfrak h}$, and $\nabla$. Basic algebraic invariants of an
infinitesimal geometric structure ${\mathfrak g}\subset J^1{\mathfrak
  t}$ are the vector bundles occurring as representations of
${\mathfrak g}$.  Associated with these representations, and a choice
of generator $\nabla$, are the {\df associated connections} and {\df
  associated differential operators}, described in
\ref{neoAssociated}.  The latter generalize the divergence, gradient,
etc.  of Riemannian geometry when $\nabla$ is the Levi-Cevita
connection.  (In Sect.\,\ref{contact} we describe these objects for
subriemannian contact three-manifolds.)  In principle, any invariant
differential operator may be expressed in terms of associated
differential operators, at least in the case of surjective
infinitesimal geometric structures. In \ref{exterior} we describe the
{\df associated exterior derivative}, and in \ref{Bianchi} analogues
of the classical Bianchi identities.

\subsection{Basic properties of generators}\lab{generator} 
Let ${\mathfrak g}\subset J^1{\mathfrak t}$ be an infinitesimal
geometric structure, with structure kernel ${\mathfrak h}$, and image
${\mathfrak t}_1\subset{\mathfrak t}$.  The projection ${\mathfrak
  g}\xrightarrow{a}{\mathfrak t}$ is of constant rank if and only if
${\mathfrak h}\subset T^*\!M\otimes{\mathfrak t}$ (or equivalently
${\mathfrak t}\subset{\mathfrak t}_1$) has constant rank (i.e., are
subalgebroids).
\begin{proposition}
                If ${\mathfrak g}\xrightarrow{a}{\mathfrak t}$ has constant 
                rank then:
        \begin{conditions}
                \item\lab{py1} ${\mathfrak g}$ admits a generator $\nabla$.
                
                \item\lab{py2} $\nabla$ is unique if and only if ${\mathfrak g}$ is 
                surjective and ${\mathfrak h}=0$.
                
              \item\lab{py3} Every $\nabla$-parallel section of
                ${\mathfrak t}_1$ is a symmetry of ${\mathfrak g}$.
        \end{conditions}
\end{proposition}
\begin{proof}[Proof of proposition and Theorem \ref{neoGenerator}]
          The constant rank hypothesis means that
\begin{equation*}
         0\longrightarrow{\mathfrak h}\longrightarrow{\mathfrak 
         g}\xrightarrow{~~~a~~~}{\mathfrak t}_1\longrightarrow 0
\end{equation*}
         is an exact sequence of vector bundles. Assuming 
         $M$ is paracompact,  it possesses a splitting $ s\colon{\mathfrak 
         t}_1\rightarrow{\mathfrak g}$ which can be extended to a 
         splitting $ s\colon{\mathfrak t}\rightarrow J^1{\mathfrak t}$ of
\begin{equation}
         0\rightarrow T^*\!M\otimes{\mathfrak t}\monomorphism J^1{\mathfrak 
         t}\rightarrow{\mathfrak t}\rightarrow 0.\mathlab{kak}
\end{equation}
         To prove \eqref{py1}, let $\nabla$ be the corresponding linear 
         connection on ${\mathfrak t}$.  
         
Conclusion \eqref{py2} follows readily from the correspondence between 
connections on ${\mathfrak t}$ and splitting of \eqref{kak}.  To prove
\eqref{py3},  let $ s\colon{\mathfrak t}\rightarrow J^1{\mathfrak t}$ 
be the splitting corresponding to a generator $\nabla$, i.e., $ sV=J^1 
V+\nabla V$. Then if $ V\subset{\mathfrak t}_1$ is $\nabla$-parallel 
then $ J^1 V=sV$. Since $ sV$ lies in ${\mathfrak g}$, by the 
definition of generators, we conclude $ V$ is a symmetry.

Assume $\nabla $ is a generator and ${\mathfrak h}=0$. Suppose that
$V\subset{\mathfrak t}$ is a symmetry, i.e., $J^1 V= sV-\nabla V$ is a
section of ${\mathfrak g}$. Then $ sV\subset{\mathfrak g}$ because
$\nabla$ is a generator, implying $\nabla V\subset{\mathfrak g}$.  So
$\nabla V\subset (T^*\!M\otimes{\mathfrak t})\cap{\mathfrak
  g}={\mathfrak h}=0$. Symmetries are thus $\nabla$-parallel. This,
together with \eqref{py3}, establishes Theorem \ref{neoGenerator}.
\end{proof}

In the remainder of this section it is tacitly assumed that all 
infinitesimal geometric structures  have constant rank  in the sense above.

\subsection{Reconstructing geometric structures from generators}%
\lab{reconstruction} %
Knowing the structure kernel $\mathfrak h$, image $\mathfrak t_1$, and 
a generator $\nabla$ of an infinitesimal geometric structure 
$\mathfrak g\subset J^1\mathfrak t$ determines it completely: the 
splitting determined by the generator determines a vector bundle 
isomorphism ${\mathfrak g}\cong{\mathfrak t}_1\oplus{\mathfrak h}$ and 
the induced Lie algebroid structure on ${\mathfrak 
t}_1\oplus{\mathfrak h}$ can be explicitly written down; see 
\eqref{jj} below.

It is not difficult to characterize those linear connections $\nabla$ 
on ${\mathfrak t}$ occurring as generators of infinitesimal geometric 
structures.  Let ${\mathfrak t}$ be an arbitrary Lie algebroid, 
${\mathfrak h}$ a subalgebroid of $T^*\!M\otimes{\mathfrak t}\subset 
J^1{\mathfrak t}$, and $\nabla$ an arbitrary linear connection on 
${\mathfrak t}$.  Let ${\mathfrak t}_1\subset{\mathfrak t}$ be an 
arbitrary subalgebroid.  We define a ${\mathfrak t}_1$-connection 
$\bar\nabla$ on $T^*\!M\otimes{\mathfrak t}$  in the obvious way:
\begin{equation}
                (\bar\nabla_V\phi) (U)=\bar\nabla_V(\phi (U))-\phi 
                (\bar\nabla_V(U)); \quad V\subset{\mathfrak 
                t}_1,\,\phi\subset{\mathfrak h},\,U\subset TM.\mathlab{storke}
\end{equation}
Here $\bar\nabla$ on the right-hand side denotes the associated 
${\mathfrak t}$-connections on ${\mathfrak t}$ and $ TM$ respectively, 
as defined in \ref{dual}.  From the formula $ s_\nabla V= J^1 V+\nabla 
V$ and the characterization of cocurvature \eqrefs{cocurvature}{dgi}, 
one readily obtains:
\begin{proposition}
  A linear connection $\nabla$ on a Lie algebroid $\mathfrak t$ is a 
  generator of some infinitesimal geometric structure 
  ${\mathfrak g}\subset J^1 {\mathfrak t}$ with structure kernel 
  ${\mathfrak h}$ and image ${\mathfrak t}_1\subset{\mathfrak t}$ if 
  and only if:
  \begin{conditions}
   \item\lab{gt} ${\mathfrak h}\subset T^*\!M\otimes{\mathfrak t}$ is 
   $\bar\nabla$-invariant, i.e., $\bar\nabla_V\phi\subset{\mathfrak 
   h}$ for all sections $ V\subset{\mathfrak t}_1$ and 
   $\phi\subset{\mathfrak h}$; and
                
   \item\lab{ala} $\cocurvature\nabla\, (V_1, V_2)\subset{\mathfrak 
   h}$ for all sections $ V_1, V_2\subset{\mathfrak t}_1$.
  \end{conditions}
  If ${\mathfrak g}\subset J^1{\mathfrak t} $ is such an infinitesimal 
  geometric structure, then the  Lie algebroid structure of 
  ${\mathfrak g}\cong{\mathfrak t}_1\oplus{\mathfrak h}$ is given by
  \begin{equation}\left\{
        \begin{split}
                \#(V\oplus\phi)&=\#V\\
                [V_1\oplus\phi_1,V_2\oplus\phi_2]&=\\
                [V_1, V_2]_{\mathfrak t_1}\oplus 
                ([\phi_1,\phi_2]_{\mathfrak 
                h}&+\bar\nabla_{V_1}\phi_2-\bar\nabla_{V_2}\phi_1-\cocurvature\nabla 
                (V_1, V_2)).
        \end{split}\right.\mathlab{jj}
  \end{equation}
\end{proposition}
\noindent We recall that cocurvature was defined in \ref{cocurvature}.

\subsection{Associated connections and differential operators}%
\lab{neoAssociated}%
Let ${\mathfrak g}\subset J^1{\mathfrak t}$ be an infinitesimal
geometric structure with structure kernel ${\mathfrak h}$, image
${\mathfrak t}_1\subset{\mathfrak t}$, and $\nabla$ a generator of
${\mathfrak g}$.  Then for each representation $E$ of ${\mathfrak g}$,
i.e., for each Lie algebroid morphism $\rho\colon{\mathfrak
  g}\rightarrow\gl(E)$, we have an {\df associated} ${\mathfrak
  t}_1$-connection $\bar\nabla$ on $ E$ (a ${\mathfrak t}$-connection
on $E$ if ${\mathfrak g}$ is surjective).  By definition, this is the
composite ${\mathfrak t}_1\xrightarrow{ s_\nabla}{\mathfrak
  g}\xrightarrow{\rho}\gl(E),$ where $ s_\nabla\colon{\mathfrak
  t}\rightarrow J^1{\mathfrak t}$ is the splitting of
\eqrefs{generator}{kak} corresponding to $\nabla$.

\begin{examples}\mbox{}
        \begin{conditions}
        \item Taking ${\mathfrak g}:=J^1{\mathfrak t}$ and
          $\rho=\adjoint^{\mathfrak t}$, we obtain
        \begin{gather*}
           \bar\nabla_ UV=\adjoint^{\mathfrak t}_{s_\nabla 
           U}V=\adjoint^{\mathfrak t}_{J^1 
           U}V+\adjoint^{\mathfrak t}_{\nabla U}V,\\
           \text{i.e.,}\quad 
           \bar\nabla_UV=\nabla_{\#V}U+[U,V]_{\mathfrak t};\quad U, 
           V\subset{\mathfrak t}.
        \end{gather*}
        This is the associated ${\mathfrak t}$-connection on
        ${\mathfrak t}$ defined already in \ref{dual}.
                 
      \item Let ${\mathfrak g}:=J^1{\mathfrak t}$ act on $TM$ via the
        composite
        \begin{equation*}
          J^1{\mathfrak t}\xrightarrow{ J^1\#} 
          J^1(TM)\xrightarrow{\adjoint^{TM}}\gl(TM).
        \end{equation*}
        Then we similarly compute
        \begin{equation*}
          \bar\nabla_UW=\#\nabla_WU+[\#U, W]_{TM};\qquad U\subset{\mathfrak 
            t},\, W\subset TM.
        \end{equation*}
        This is the associated ${\mathfrak t}$-connection on $TM$
        defined in \ref{dual}.
                 
      \item\lab{kiwi} An arbitrary infinitesimal geometric structure
        ${\mathfrak g}\subset J^1{\mathfrak t}$ acts on its structure
        kernel ${\mathfrak h}$ via bracket: $\rho_X Y:=[X,
        Y]_{\mathfrak g}$. (Here we are assuming that ${\mathfrak h}$
        has constant rank.) The associated ${\mathfrak
          t}_1$-connection on ${\mathfrak h}$ is simply the connection
        $\bar\nabla$ appearing in Proposition \ref{reconstruction} and
        satisfying \eqrefs{reconstruction}{storke}.
				 
      \item\lab{kiwi2}If ${\mathfrak g}\subset J^1{\mathfrak t}$ is an
        infinitesimal geometric structure and $\nabla$ a linear
        connection on ${\mathfrak g}$ (rather than ${\mathfrak t}$)
        then $\nabla$ generates $ J^1{\mathfrak g}$, which acts on
        ${\mathfrak t}$ via the composite
	\begin{equation*}
          J^1{\mathfrak g}\xrightarrow{J^1 a}J^1{\mathfrak 
            t}\xrightarrow{\adjoint^{\mathfrak t}}\gl({\mathfrak 
            t}),
	\end{equation*}
	where $ a\colon{\mathfrak g}\rightarrow{\mathfrak t}$ is the
        projection. The associated ${\mathfrak g}$-connection
        $\bar\nabla$ on ${\mathfrak t}$ is given by
        $\bar\nabla_XW=a\nabla_{\#W}X+[aX, W]_{\mathfrak t}$.  For an
        application, see \ref{natural}.
    \end{conditions}
\end{examples}

Let a  {\df ${\mathfrak g}$-tensor} be any section $\sigma\subset E$ 
of a ${\mathfrak g}$-representation $E$. Then, by  the definition of the
associated connections, we have:
\begin{proposition}
         A ${\mathfrak g}$-tensor is ${\mathfrak g}$-invariant if and only 
         if it is simultaneously ${\mathfrak h}$-invariant and 
         $\bar\nabla$-parallel.
\end{proposition}

The {\df associated derivative} of a ${\mathfrak g}$-tensor
$\sigma\in\Gamma( E)$ is defined to be $\bar\nabla\sigma \in \Gamma
({\mathfrak t}_1^*\otimes E)$, where $\bar\nabla$ is the associated
${\mathfrak t}_1$-connection on $E$.  Assume ${\mathfrak g}$ is {\df
  image-reduced}; by this we mean that ${\mathfrak
  t}_1\subset{\mathfrak t}$ is invariant under the adjoint action of
${\mathfrak g}\subset J^1{\mathfrak t}$, which is true, for example,
if ${\mathfrak g}$ is surjective. (Image reduction is described in
\ref{imageReduction}.) Then ${\mathfrak t}_1$ is a ${\mathfrak
  g}$-representation, implying $\bar\nabla\sigma$ is another
${\mathfrak g}$-tensor.  That is, the ${\mathfrak g}$-tensors will be
closed under associated derivative.  In particular, the derivative can
be iterated to obtain higher order differential operators.

Additionally supposing that all ${\mathfrak g}$-representations may be
direct-sum decomposed into $\mathfrak g$-representations coming from
some collection $ E_i$ ($i\in I$) of irreducible ones, we have
\begin{equation*}
         {\mathfrak t}_1^*\otimes E_i\cong E_{n_{i1}}\oplus E_{n_{i2}}\oplus 
         E_{n_{i3}}\oplus\cdots\quad\text{(finitely many non-zero terms),}
\end{equation*}
for some $ n_{ij}\in I$, and obtain a corresponding decomposition,
\begin{equation*}
        \bar\nabla|\Gamma(E_i)=\partial_{i1}\oplus
        \partial_{i2}\oplus\partial_{i3}\oplus\cdots\enspace.
\end{equation*}
We call the differential operators $ \partial_{ij}\colon 
\Gamma (E_i)\rightarrow \Gamma( E_{n_{ij}})$ ($i,j\in I$) the {\df associated 
differential operators}; all differential operators which can be 
constructed algebraically out of associated connections $\bar\nabla$ 
are combinations of these basic ones.

If there is a {\em canonical} way in which to choose the generator 
$\nabla$ then the associated differential operators become {\em 
invariant} differential operators associated with the infinitesimal 
geometric structure ${\mathfrak g}$.  Significant cases in point are:
\begin{conditions}
        \item\lab{hodge} The case where ${\mathfrak t}$ is a Cartan 
        algebroid discussed in \ref{cartan}.  Here ${\mathfrak 
        g}$-representations are just ${\mathfrak t}$-representations 
        because ${\mathfrak g}\cong{\mathfrak t}$.
        
        \item The case where the generator $\nabla$ of ${\mathfrak g}$ is 
        {\em unique}, i.e., Theorem \ref{neoGenerator} applies, reducing 
        the situation to case \eqref{hodge} above.
        
        \item The case where torsion $\torsion\bar\nabla$ has a natural 
        `normalization'; see \ref{normalizing}.
\end{conditions}
\noindent For invariant differential operators associated with 
subriemannian contact threee-manifolds, see Sect.\,\ref{contact}.

\subsection{The associated exterior derivative}
\lab{exterior}%
Let ${\mathfrak g}\subset J^1{\mathfrak t}$ be an infinitesimal 
geometric structure with structure kernel ${\mathfrak h}$.  Then a 
{\df differential form of type ${\mathfrak g}$ and degree $k$} is a 
section $\theta\subset\alternating^k({\mathfrak t}_1)\otimes E$, where 
${\mathfrak t}_1\subset{\mathfrak t}$ is the image of ${\mathfrak g}$ 
and $E$ some ${\mathfrak g}$-representation.  (We use ${\mathfrak 
t}_1$, rather than ${\mathfrak t}$, to ensure \eqref{tu} below.)  The 
{\df exterior derivative} 
$d_{\bar\nabla}\theta\subset\alternating^k({\mathfrak t}_1)\otimes E$ 
of $\theta$ is defined in the obvious way.  For example,
\begin{align*}
         d_{\bar\nabla}\theta\,(U_1)&:=\bar\nabla_{U_1}\theta,\quad\text{for 
         $k=0$,}\\
         \text{and}\enspace 
         d_{\bar\nabla}\theta\,(U_1,U_2)&:=\bar\nabla_{U_1}(\theta (U_2))- 
         \bar\nabla_{U_2}(\theta (U_1))-\theta (\,[U_1,U_2]_{{\mathfrak 
         t}_1}),\quad\text{for $ k=1$.}
\end{align*}
Wedge products of ${\mathfrak g}$-type differential forms are likewise 
defined by familiar formulas.

The principal invariants of the generator $\nabla$ are the ${\mathfrak 
g}$-type differential forms $T\subset\alternating^2 ({\mathfrak 
t}_1)\otimes{\mathfrak t}$ and $\Omega\subset\alternating^2 
({\mathfrak t}_1)\otimes{\mathfrak h}$, obtained by restricting 
$\torsion\bar\nabla\subset\alternating^2 ({\mathfrak 
t})\otimes{\mathfrak t}$ and $-\cocurvature\nabla\subset\alternating^2 
({\mathfrak t})\otimes T^*\!M\otimes{\mathfrak t}$ to ${\mathfrak 
t}_1$ (note carefully the minus sign).
\begin{proposition}\mbox{}
        \begin{conditions}
                \item\lab{won} If $\theta_1$ and $\theta_2$ are ${\mathfrak 
                g}$-type differential forms, then 
                $$d_{\bar\nabla}(\theta_1\wedge\theta_2)=d_{\bar\nabla}\theta_1 
                \wedge\theta_2+(-1)^k\,\theta_1\wedge
                d_{\bar\nabla}\theta_2,$$
                where $ k$ is the degree of $\theta_1$.
                
                \item\lab{tu} For any ${\mathfrak g}$-type differential form 
                $\theta$, we have 
                $$d_{\bar\nabla}^2\theta=\Omega\wedge\theta.$$ Here the wedge 
                implies a contraction 
                $\phi\otimes\sigma\mapsto\phi\cdot\sigma\colon{\mathfrak 
                h}\otimes E\rightarrow E$, defined by the representation of 
                ${\mathfrak h}$ on $E$.
        \end{conditions}
\end{proposition}
\begin{proof}[Proof of \eqref{tu}]
The general case can easily be reduced to the $ k=0$ case that we 
prove now.  Letting $ s\colon{\mathfrak t}\rightarrow J^1{\mathfrak 
t}$ denote the splitting of \eqrefs{generator}{kak} determined by 
$\nabla$, we compute, for arbitrary $ U_1, U_2\subset{\mathfrak t}_1$,
\begin{align*}
   d_{\bar\nabla}^2\theta\,(U_1,U_2)&=
     \bar\nabla_{U_1}\bar\nabla_{U_2}\theta-
     \bar\nabla_{U_2}\bar\nabla_{U_1}\theta
     -\bar\nabla_{[U_1,U_2]_{{\mathfrak 
     t}_1}}\theta\\
   &=sU_1\cdot(sU_2\cdot\theta)-
     sU_2\cdot(sU_1\cdot\theta)-
         s[U_1,U_2]_{{\mathfrak 
     t}_1}\cdot\theta\\
   &=(sU_1\cdot(sU_2\cdot\theta)-
      sU_2\cdot(sU_1\cdot\theta)-
         [sU_1,sU_2]_{\mathfrak g}\cdot\theta)\\
   &-\cocurvature\nabla\,(U_1,U_2)\cdot\theta, \quad\text{by 
   \eqrefs{algebra}{dgi}}\\
   &=0+\Omega (U_1,U_2)\cdot\theta.
\end{align*}
\end{proof}

\subsection{Bianchi identities}\lab{Bianchi}%
Generalizing the classical situation, the Bianchi identities below exhibit 
certain algebraic and differential dependencies between $ T$ and 
$\Omega$, rooted in the equality of mixed partial derivatives.  First, 
since $T=d_{\bar\nabla}i$, where $i\subset{\mathfrak 
t}_1^*\otimes{\mathfrak t} $ denotes the inclusion ${\mathfrak 
t}_1\subset{\mathfrak t}$, We deduce from \eqref{tu} above,
\begin{equation}
         d_{\bar\nabla}T=\Omega\wedge i.\mathlab{num1}
\end{equation}
Next, assume ${\mathfrak g}$ admits a representation $E$ for which the 
restricted representation ${\mathfrak h}\rightarrow\gl(E)$ is faithful 
(injective), and let $\theta\subset E$ be a section, viewed as a 
${\mathfrak g}$-type differential form of degree zero.  Then, 
combining parts \eqref{won} and \eqref{tu} of the preceding 
proposition, we obtain 
$d_{\bar\nabla}^3\theta=d_{\bar\nabla}\Omega\wedge\theta+\Omega\wedge 
d_{\bar\nabla}\theta$.  Applying part \eqref{tu} again, we conclude 
that $ d_{\bar\nabla}\Omega\wedge\theta=0$.  Since $\theta$ is 
arbitrary and ${\mathfrak h}$ acts faithfully on $E$, we obtain
\begin{equation}
         d_{\bar\nabla}\Omega=0.\mathlab{num2}
\end{equation}
A little manipulation allows us to write \eqref{num1} and \eqref{num2} 
in the form
\begin{multline}
        (\bar\nabla_{U_3} T) (U_1, U_2)+
                T(T(U_1,U_2), U_3)\\
        +\Omega(U_1, U_2)\#U_3+\text{1-2-3-cyclic terms}=0\qquad
        \text{(Bianchi I),}\mathlab{BianchiI}
\end{multline}
\begin{multline}
         (\bar\nabla_{U_3}\Omega) (U_1, U_2)+\Omega 
         ( T (U_1, U_2), U_3)\\
         +\text{1-2-3-cyclic terms}=0\qquad\text{(Bianchi II).}
         \mathlab{BianchiII}
\end{multline}
\begin{example}
If $ \nabla $ is the Cartan Connection on some Lie algebroid
${\mathfrak t} $ and ${\mathfrak g} \subset J^1 {\mathfrak t} $ the corresponding infinitesimal geometric structure  (see
\ref{cartan}) then $T = \torsion \bar\nabla $ and $\Omega=0$. Bianchi I becomes 
\begin{equation*}
  (\bar\nabla_{U_3}{\torsion \bar\nabla})(U_ 1, U_2)+ \torsion  \bar\nabla (\torsion \bar\nabla(U_1,U_2),U_3)+\text{1-2-3-cyclic terms}=0.
\end{equation*}
\end{example}

\section{Elementary reduction and image reduction}%
\lab{elementary} In this section we study elementary reduction, as
well as a cruder alternative we call image reduction. These techniques
are useful when an infinitesimal geometric structures fails to be
surjective, and in particular to intransitive infinitesimal geometric
structures on $TM $. A simple application to smooth functions on a
Riemannian three-manifold is included.

\subsection{Image reduction}\lab{imageReduction} Let ${\mathfrak g}
\subset J^1 {\mathfrak t} $ be an infinitesimal geometric structure
with structure kernel ${\mathfrak h} \subset T^*\!M \otimes {\mathfrak
  t} $ and image ${\mathfrak t}_1 \subset {\mathfrak t} $. Assume
${\mathfrak h}$ (or equivalently ${\mathfrak t}_1$) has constant
rank. Then the {\df image reduction} of $\mathfrak g$ is simply the
isotropy ${\mathfrak g}_{{\mathfrak t}_1} \subset J^1 {\mathfrak g} $
of ${\mathfrak t}_1 \subset {\mathfrak t}$, under the adjoint
representation of $\mathfrak g \subset J^1 {\mathfrak t}$ on
${\mathfrak t}$. It is not hard to show that image reduction is cruder
than elementary reduction, as defined in \ref{cradle}, and described
further below. Nevertheless, it is usually easier to apply image
reduction and this may simplify the subsequent application of
elementary reduction.

\subsection{Elementary reduction}\lab{element}
With $\mathfrak g \subset J^1 {\mathfrak t}$, ${\mathfrak h}$, and
${\mathfrak t}_1$ as above, let ${\mathfrak g}_1$ denote the
elementary reduction of $\mathfrak g$ (see \ref{cradle}). The
structure kernel of ${{\mathfrak g}}_1$ is $${\mathfrak h}_1\equiv{\mathfrak
  h}\cap(T^*\!M \otimes {\mathfrak t}_1).$$ One can compute the image
${\mathfrak t}_2 \subset {\mathfrak t}_1$ of ${\mathfrak g}_1$ if one
knows a generator $\nabla $ of $\mathfrak g$:
\begin{proposition}
  There is vector bundle morphism
  \begin{equation*}
    {\mathfrak t}_1 \xrightarrow{b}(T^*\!M \otimes {\mathfrak t})/(T^*\!M \otimes {\mathfrak t}_1 + {\mathfrak h}),
  \end{equation*}
  whose corresponding map of section spaces is 
\begin{equation*}
  U \mapsto \nabla U \mod (T^*\!M \otimes {\mathfrak t}_1 + {\mathfrak h}).
\end{equation*}
The morphism $b$ is independent of the choice of generator $\nabla $
and $\kernel b={\mathfrak t}_2$.
\end{proposition}
\begin{proof}
  Begin by observing that the one-jet $J^1 U(m)\in J^1 {\mathfrak t}_1$ lies in $\mathfrak g$ if and only if $\nabla U(m)$ lies in ${\mathfrak h}$. So we define a morphism 
\begin{equation*}
  J^1 {\mathfrak t}_1 \xrightarrow{B}
  (T^*\!M \otimes {\mathfrak t})/{\mathfrak h},
\end{equation*}
which on sections is the map $J^1 U\mapsto \nabla U \mod {\mathfrak
  h}$, and have ${\mathfrak g}_1= \kernel B$. The proposition now
follows from an application of Lemma \ref{algebraicLemma} to the
morphism $B$; one uses the fact that the sequence
\begin{equation*}
  0 \rightarrow \frac{T^*\!M \otimes {\mathfrak t}_1}{{\mathfrak h}}
  \rightarrow \frac{T^*\!M \otimes {\mathfrak t}}{{\mathfrak h}}
  \rightarrow \frac{T^*\!M \otimes {\mathfrak t}}{T^*\!M \otimes {\mathfrak t}_1 + {\mathfrak h}} \rightarrow 0
\end{equation*}
is exact.
\end{proof}
\begin{remark}
  A significant simplification occurs if ${\mathfrak h}_1={\mathfrak
    h}$, i.e., if ${\mathfrak h} \subset T^*\!M \otimes {\mathfrak t}
  $ lies entirely within $T^*\!M \otimes {\mathfrak t}_1$. Then we may
  view $b$ as the map
  \begin{gather*}
    {\mathfrak t}_1 \xrightarrow{b}T^*\!M \otimes 
    ({\mathfrak t}/{\mathfrak t}_1)\\
    b(U)V=\nabla_VU \mod {\mathfrak t}_1.
  \end{gather*}
  In particular, the proposition will imply that ${\mathfrak
    g}_1={\mathfrak g}$ if and only if ${\mathfrak t}_1$ is
  $\nabla$-invariant. In that case $\nabla$ drops to a linear
  connection on ${\mathfrak t}_1$ which generates $\mathfrak g$, as a
  {\em surjective} infinitesimal geometric structure on ${\mathfrak
    t}_1$.
\end{remark}

\subsection{Functions on a Riemannian three-manifold}\lab{eggs}
As an illustration, consider the (infinitesimal) symmetries of a
smooth function $f$ on an oriented Riemannian three-manifold $M$, with
metric $\sigma$. By {\df symmetries}, we mean the Killing fields of
$\sigma$ preserving $f$. In the terminology of \ref{structure}, these
are the symmetries of the joint isotropy
\begin{equation*}
  (J^1 (TM))_{\sigma,f} \subset J^1 (TM)
\end{equation*}
of $\sigma$ and $f$, under the relevant representations determined by
the adjoint representation of $J^1(TM)$ on $TM$. Any such symmetry must
also preserve $df$, and so we have an immediate reduction,
\begin{equation}
    (J^1 (TM))_{\sigma,f,df} \subset J^1 (TM).\mathlab{herewith}
\end{equation}

Let $E\equiv \frac{1}{2}\|\grad f\|^2$ denote the `energy' of $f$
($\grad f\equiv \sigma^{-1}(df)$), and assume that $df$ and $dE$ are
everywhere linearly independent. It follows that the connected
components of the joint level-sets of $f$ and $E$ constitute a
rank-one foliation on $M$. We denote by ${T}$ the unit vector field
tangent to this foliation, directed so as to make $\{{T}, \grad f
,\grad E\}$ positively oriented.

Define $J \subset T^*\!M \otimes TM$ by $JU:= {\mathbf n} \times U$, where ${\mathbf n}:=\gradient f / \|\gradient f\|$ (so that $J$ restricts to a complex structure on level sets of $f$). Then the reduction in \eqref{herewith} has a rank-one structure kernel spanned by $J$.
Using Lemma \ref{algebraicLemma}, it is not hard to see that its image
is $\langle {T} \rangle=\kernel df \cap \kernel dE$. We therefore pass
to its image-reduction, i.e., the joint isotropy,
%\ref{imageReduction},
\begin{equation*}
  {\mathfrak g}\equiv (J^1 (TM))_{\sigma,f,df,\langle {T}  \rangle}
  \subset J^1 (TM).
\end{equation*}
We observe that $\mathfrak g$ has trivial structure kernel,
${\mathfrak h}=0$, and image ${\mathfrak t}_1:=\langle {T}\rangle$. To
see this one applies Lemma \ref{algebraicLemma} to the morphism,
\begin{gather*}
J^1(TM)_{\sigma,f,df}\rightarrow TM/\langle {T} \rangle\\
X \mapsto \adjoint_X {T} \mod \langle {T} \rangle,  
\end{gather*}
which has $\mathfrak g$ as kernel.

As $\mathfrak g$ itself is evidently stable under image-reduction, we
now turn to elementary reduction. By Proposition \ref{generator},
$\mathfrak g$ has a generator $\nabla$. Because ${T}$ spans the image of ${\mathfrak g}=(J^1(TM))_{\sigma,f,df,\langle{T}\rangle}$, we must have 
\begin{equation*}
  \bar\nabla_{T} \sigma=0,\quad,\bar\nabla_{T}\grad f=0\quad\text{and}\quad \bar\nabla_{T} {T} \subset \langle {T}\rangle,
\end{equation*}
where $\bar\nabla_U V\equiv \nabla_VU+[U,V]$. From these identities we may deduce
\begin{align}
  \nabla_{T}{T}&=0\mathlab{keeda1},\\
  \nabla_{J {T}}{T}&=[J {T},{T}]\mathlab{keeda2},\\
   \nabla_{\grad f}{T}&=[\grad f, {T}]\mathlab{keeda3}.  
\end{align}
% Indeed, $\bar\nabla_{T} \sigma=0$ and $\sigma({\mathbf
%   m},{T})=1$ (differentiated along ${T} $) already
% imply that $\bar\nabla_{T} {T} $ and ${T} $
% are perpendicular. The requirement $\bar\nabla_{T} {T}
% \subset \langle {T} \rangle$ forces $\bar\nabla_{T}
% {T}=0$ and \eqref{keeda1} follows. Using the properties of
% $\bar\nabla $ above, it is not difficult to also show that
% $\bar\nabla_{T}(J {T})=0$.  This implies
% \eqref{keeda2}. Equation \eqref{keeda3} follows immediately from
% $\bar\nabla_{T} \grad f=0$.

Assume that the image ${\mathfrak t}_2 \subset \langle{T}\rangle$ of
the elementary reduction ${\mathfrak g}_1$ of $\mathfrak g$ has
constant rank. Then either ${\mathfrak t}_2=0$, in which case
${\mathfrak g}_1=0$ and there is no possibility for symmetry, or
${\mathfrak t}_2=\langle{T}\rangle$, in which case ${\mathfrak
  g}_1={\mathfrak g}$.  (Note here that the structure kernel of
$\mathfrak g$ is already trivial). According to Remark \ref{element},
we are in the latter case when $\langle{T}\rangle $ is $\nabla
$-invariant. Necessary and sufficient conditions, following from
\eqref{keeda1}--\eqref{keeda3}, are:
\begin{equation}
  [J{T},{T}] \subset \langle{T}\rangle\quad\text{and}
  \quad [\grad f, {T}] \subset \langle{T}\rangle.\mathlab{gui}
\end{equation}
Under these conditions, $\nabla$ descends to a Cartan connection on
the line bundle $\langle{T}\rangle $, whose curvature is the
obstruction to the existence of a Killing filed preserving $f$ (lying
necessarily along $\langle{T}\rangle $). As it turns out, this
curvature generally has a single non-trivial component, $\curvature
\nabla (JT,\grad f)$, which the interested reader is invited to
compute.

The brackets in condition \eqref{gui} can be expressed in terms of the
Levi-Cevita connection (e.g., $[JT,T]=\levi_{J{T}}T-\levi_TJ{T}$),
which leads to the equivalent condition that the rank-two
distributions $(JT)^\perp$ and $T^\perp$ be integrable and geodesic,
respectively. (A distribution is {\df geodesic} if it has trivial
second fundamental form).

\section{Prolongation and torsion}\lab{sectionX}
In this section we characterize the prolongation ${\mathfrak g}
^{(1)}:= J^1 {\mathfrak g} \cap J ^2 {\mathfrak t}$ as the joint
isotropy of a tautological one-form $a$ and its `torsion' $da$. These
tensors are analogues of classical objects bearing the same name but
the characterization is only valid when $\mathfrak g$ is
transitive. We begin, however, with a useful characterization of the
subbundle $J^2 {\mathfrak g} \subset J^1(J^1{\mathfrak g})$ that is
completely general.

This section concludes with the reformulation of several classical
constructions associated with torsion.

\subsection{Prolongation}\lab{prolongation}
Let ${\mathfrak t}$ be an {\em arbitrary vector bundle} over $M$. Then
there is a natural inclusion of vector bundles $J^2 {\mathfrak t}
\monomorphism J^1(J^1 {\mathfrak t}) $ which takes $J^2W(m)$ to
$J^1(J^1 W)(m)$; $W \subset J^1{\mathfrak t}$. As a basic fact one has
the following:
\begin{lemma}
  For any section $X \subset J^1 {\mathfrak t}$, $X$ is holonomic if
  and only if $J^1 X \subset J^2 {\mathfrak t} $.
\end{lemma}
\begin{proof}
  See Appendix \ref{pro}.
\end{proof}
\noindent%
Now the definition of prolongation, ${\mathfrak g} ^{(1)}\equiv J^1
{\mathfrak g} \cap J ^2 {\mathfrak t}$, of some subbundle ${\mathfrak
  g} \subset J^1{\mathfrak t} $ makes sense in general, but suppose
for the moment that ${\mathfrak t}$ is a Lie algebroid and $\mathfrak
g \subset J^1 {\mathfrak t} $ is an infinitesimal geometric structure
on ${\mathfrak t}$. Then $J^2 {\mathfrak t} \subset J^1 (J^1
{\mathfrak t})$ is a subalgebroid, implying ${\mathfrak g} ^{(1)} $ is
an infinitesimal geometric structure on $\mathfrak g$ whenever
${\mathfrak g} ^{(1)} $ has constant rank. Let $W \subset {\mathfrak
  t} $ be a symmetry of $\mathfrak g$. Then a straightforward
consequence of definitions is that $J^1 W$ is a section of $\mathfrak
g$ and, moreover, a symmetry of ${\mathfrak g} ^{(1)} $. In fact, it
is a consequence of the lemma above that {\em all} symmetries of
${\mathfrak g} ^{(1)} $ arise in this way:
\begin{proposition}
  If ${\mathfrak g} \subset J^1{\mathfrak t} $ is an infinitesimal
  geometric structure, then a section $W \subset {\mathfrak t} $ is a
  symmetry of $\mathfrak g$ if and only if $J^1 W \subset {\mathfrak
    g} $ is a symmetry of ${\mathfrak g} ^{(1)} $.
\end{proposition}
\noindent%
Since $J^1\colon \Gamma({\mathfrak t})\rightarrow \Gamma(J^1
{\mathfrak t})$ is injective, this establishes a one-to-one
correspondence between the symmetries of $\mathfrak g$ and those of
$\mathfrak g ^{(1)} $.

\subsection{More on $J^2 {\mathfrak t} $ as a subbundle of 
$J^1(J^1 {\mathfrak t}))$}%
\lab{boot} The definition of ${\mathfrak g}^{(1)}$ is difficult to
work with unless one has the right characterization of $J^2 {\mathfrak
  t}$, {\em as a subbundle of} $J^1(J^1{\mathfrak t})$. We now
characterize $J^2 {\mathfrak t} \subset J^1 (J^1 {\mathfrak t})$ as
the kernel of a natural morphism $J^2_+{\mathfrak t}\rightarrow
\alternating^2(TM)\otimes {\mathfrak t} $, where $J^2_+{\mathfrak t}
\subset J^1(J^1{\mathfrak t}) $ is a larger subbundle that is itself
the kernel of a natural morphism $J^1(J^1{\mathfrak t})\rightarrow
T^*\!M \otimes {\mathfrak t}$. This characterization holds for an
arbitrary vector bundle ${\mathfrak t}$ over $M$.

Define a differential operator
\begin{gather*}
  {\mathcal D}:\Gamma (J^1 {\mathfrak
  t})\rightarrow \Gamma(T^*\!M \otimes {\mathfrak t})    \\
  \text{by}\quad{\mathcal D}X\equiv X - J^1 (aX), 
\end{gather*}
where $a\colon J^1 {\mathfrak t} \rightarrow {\mathfrak t} $ is the
natural projection. We call $\mathcal D$ the {\df deviation} of
$X$. Notice that $X$ is holonomic if and only if ${\mathcal D}
X=0$. Also, ${\mathcal D} \phi =\phi $ for all sections $\phi \subset
T^*\!M \otimes {\mathfrak t} \subset J^1{\mathfrak t}$. We write
${\mathcal D}_VX\equiv ({\mathcal D} X)V$, for $V\subset TM$ and have
a Leibnitz-type identity
\begin{equation}
  {\mathcal D}_V(fX)=f {\mathcal D}_VX + df(V)aX,\mathlab{Leib}
\end{equation}
for arbitrary smooth functions $f$ on $M$.

The above construction, holding for arbitrary ${\mathfrak t}$, may be
applied in the particular case that ${\mathfrak t}$ is replaced by
$J^1{\mathfrak t} $. This delivers an operator
$\Gamma(J^1(J^1{\mathfrak t}))\rightarrow T^*\!M \otimes J^1{\mathfrak
  t}$ which will also be denoted $\mathcal D$. In the formulas above
the projection $a$ gets replaced by the natural projection $p\colon
J^1(J^1{\mathfrak t})\rightarrow J^1{\mathfrak t}$.

\begin{proposition}%
[Characterization of $J^2 {\mathfrak t}\subset J^1(J^1{\mathfrak t})$ as a kernel]
For an arbitrary vector bundle ${\mathfrak t}$, one has $J^2
{\mathfrak t}=\kernel \omega_2$, where
\begin{equation*}
  \omega_2\colon J^2_+{\mathfrak t} \rightarrow \alternating^2(TM)\otimes {\mathfrak t} 
\end{equation*}
is a vector bundle morphism well defined by 
\begin{equation*}
  (\omega_2 \xi)(V_1,V_2)\equiv{\mathcal D}_{V_1}{\mathcal D}_{V_2}\xi 
                        -{\mathcal D}_{V_2}{\mathcal D}_{V_1}\xi 
                        -a {\mathcal D}_{[V_1,V_2]}\xi.  
\end{equation*}
Here $J^2_+{\mathfrak t} \subset J^1(J^1 {\mathfrak t}))$ is the
kernel of the vector bundle morphism
\begin{equation*}
  \omega_1\colon J^1 (J^1{\mathfrak t})\rightarrow 
  T^*\!M \otimes J^1{\mathfrak t},
\end{equation*}
well defined by 
\begin{equation*}
  (\omega_1 \xi)V\equiv {\mathcal D}_V(p \xi)-a {\mathcal D}_V \xi.
\end{equation*}
\end{proposition}
\noindent%
In this proposition some ${\mathcal D} $'s are operators
$\Gamma(J^1{\mathfrak t})\rightarrow \Gamma(T^*\!M \otimes {\mathfrak
  t})$, while others are operators $\Gamma(J^1(J^1{\mathfrak
  t}))\rightarrow \Gamma(T^*\!M \otimes J^1{\mathfrak t})$.  All
ambiguity is mitigated by the context.

Since the proposition above is just a general fact about vector
bundles, its proof is relegated to Appendix \ref{pro}.

\subsection{Torsion}\lab{tautological}
We now return to the case that ${\mathfrak g}$ is an infinitesimal
geometric structure on a Lie algebroid ${\mathfrak t}$. Applying the
general results above we obtain a characterization of ${\mathfrak
  g}^{(1)} $ as an isotropy subalgebroid.

Regard the restriction $a\colon{\mathfrak g}\rightarrow{\mathfrak t}$
of $ J^1{\mathfrak t}\rightarrow{\mathfrak t}$ as a ${\mathfrak
  t}$-valued ${\mathfrak g}$-form of degree one; this is the {\df
  tautological one-form}.  The adjoint representation of
$J^1{\mathfrak t}$ on ${\mathfrak t}$ restricts to a representation of
${\mathfrak g}$ on ${\mathfrak t}$.  So the exterior derivative $ da$
is a well defined ${\mathfrak t}$-valued ${\mathfrak g}$-form, of
degree two.  This is the {\df torsion} of the structure.  Explicitly,
\begin{equation*}
         da(X_1,X_2)=\adjoint^{\mathfrak t}_{X_1}(aX_2)
 -\adjoint^{\mathfrak t}_{X_2}(aX_1)-a[X_1, X_2].
\end{equation*}
Now $ a$ and $ da$ are sections of ${\mathfrak 
g}^*\otimes{\mathfrak t}$ and $\alternating^2 ({\mathfrak 
g})\otimes{\mathfrak t}$, respectively. We get representations of $ 
J^1{\mathfrak g}$ on these spaces by taking $ J^1{\mathfrak g}$ to act 
on ${\mathfrak g}$ via adjoint action, and on ${\mathfrak t}$ via the 
composite
\begin{gather*}
         J^1{\mathfrak g}\xrightarrow{ p}{\mathfrak g}\monomorphism 
         J^1{\mathfrak t}\xrightarrow{\adjoint^{\mathfrak 
         t}}\gl ({\mathfrak t}),\\
                 \text{i.e.,}\quad \xi\cdot W=\adjoint_{p\xi}^{\mathfrak t} W;\qquad
                 X\subset{\mathfrak g},\,W\subset{\mathfrak t}.
\end{gather*}
Here $ p\colon J^1{\mathfrak g}\rightarrow{\mathfrak g}$ denotes the
canonical projection.  We can accordingly define subsets of $J^1\mathfrak g$,
\begin{equation*}
         (J^1{\mathfrak g})_a:=\text{isotropy of $a$,}\quad\text{and}\quad
         (J^1{\mathfrak g})_{da}:=\text{isotropy of $da$.}
\end{equation*}
\begin{proposition}
If $\mathfrak g$ is transitive, then its prolongation ${\mathfrak g} ^{(1)} $ is the joint isotropy of $a$ and $da$,  i.e., 
  ${\mathfrak g}^{(1)}=(J^1 {\mathfrak g})_{a,da}
  \equiv(J^1{\mathfrak g})_a\cap (J^1 {\mathfrak g})_{da}\subset J^1\mathfrak g.$
\end{proposition}
\noindent%
\begin{remark}
  If ${\mathfrak g}$ is intransitive, then $(J^1 {\mathfrak
    g})_{a,da}$ generally has too much symmetry to be the prolongation
  of $\mathfrak g$: every section of ${\mathfrak h}\subset
  T^*\!M\otimes{\mathfrak t}$ annihilating tangent vectors in the
  image of the anchor $\#\colon{\mathfrak g}\rightarrow TM$ turns out
  to be a symmetry of $(J^1 {\mathfrak g})_{a,da}$ that is {\em not} a
  symmetry of ${\mathfrak g}^{(1)}$.
\end{remark}
%\noindent
The proposition is an easy corollary of Proposition \ref{boot} and the
following observation:
\begin{lemma}
  Let ${\mathfrak g} \subset J^1{\mathfrak t} $ be a (possible
  intransitive) infinitesimal geometric structure on ${\mathfrak t}$
  and let $\mathcal D$ denote the deviation operator discussed in
  \ref{boot}. Then for an arbitrary section $\xi \subset J^1
  {\mathfrak g} $, one has
  \begin{conditions}
  \item\lab{flute1} $(\xi \cdot a)X={\mathcal D}_{\#X}(p \xi)-a
    {\mathcal D}_{\#X}\xi$, and
  \item\lab{flute2} $(\xi \cdot da)(X_1,X_2)=d(\xi \cdot a)(X_1,X_2) +
    {\mathcal D}_{\#X_1}{\mathcal D}_{\#X_2}\xi -{\mathcal
      D}_{\#X_2}{\mathcal D}_{\#X_1}\xi -a {\mathcal
      D}_{[\#X_1,\#X_2]}\xi$.
  \end{conditions}
\end{lemma}
\noindent%
Here $X,X_1,X_2 \subset {\mathfrak g}$ are arbitrary sections.
\begin{proof}[Proof of lemma]
  Begin by observing that
  \begin{equation*}
    (\xi\cdot a)(X)=
    \adjoint_{p\xi}^{\mathfrak t}(aX)-a(\adjoint_\xi^{\mathfrak g}X).
  \end{equation*}
  Since $ a\colon{\mathfrak g}\rightarrow{\mathfrak t}$ is a
  Lie algebroid morphism, the identity \eqrefs{adjective}{America}
  gives us $a(\adjoint_\xi^{\mathfrak g}X)=\adjoint_{(J^1
    a)\xi}^{\mathfrak t} (aX)$, and so
  \begin{equation*}
      (\xi\cdot a)(X)=\adjoint^{\mathfrak t}_{p\xi-(J^1 a)\xi}(aX).
  \end{equation*}
  Note here that $ J^1 a\colon J^1{\mathfrak g}\rightarrow
  J^1{\mathfrak t}$ is the morphism sending $ J^1 X(m)$ to $ J^1
  (aX)(m)$. Because $ p\xi-(J^1 a)\xi$ is a section of the kernel of
  $J^1{\mathfrak t}\rightarrow {\mathfrak t}$, we may view it as a
  section of $ T^*\!M\otimes{\mathfrak t}$ and, applying
  \eqrefs{adjective}{gull}, obtain
  \begin{equation}
      (\xi\cdot a)(X)=(p\xi-(J^1a)\xi)(\#X).\mathlab{ballet}
  \end{equation}
  On the other hand, since $\xi=J^1(p \xi)+{\mathcal D} \xi$, we have
  $(J^1 a)\xi=J^1(ap \xi)+(J^1 a)({\mathcal D}\xi)$, implying
  \begin{align*}
    p \xi-(J^1 a)\xi&={\mathcal D}(p \xi)-(J^1 a)({\mathcal D} \xi),\\
    \implies (p \xi - (J^1 a)\xi)V&={\mathcal D}_V(p \xi)-a
    {\mathcal D}_V \xi, \qquad V \subset  TM.
  \end{align*}
  Combining this with \eqref{ballet} gives \eqref{flute1}.

  It is not too difficult to show that $\xi \cdot da=d(\xi \cdot da)$
  whenever $\xi$ is {\em holonomic}. Therefore 
\begin{equation*}
  J^1(p \xi) \cdot da = d(J^1(p \xi) \cdot a)
  =d(\xi \cdot a)-d({\mathcal D} \xi \cdot a),
\end{equation*}
because $J^1(p \xi)=\xi- {\mathcal D} \xi$. We now compute
\begin{align*}
  \xi \cdot da&=(J^1(p \xi)+ {\mathcal D} \xi)\cdot da\\
  &=J^1(p \xi)\cdot da + ({\mathcal D} \xi) \cdot da\\
  &=d(\xi \cdot  a) - d({\mathcal D} \xi \cdot a) + {\mathcal D} \xi \cdot da.
\end{align*}
Note that $({\mathcal D} \xi \cdot a)X=-a{\mathcal D}_{\#X}\xi $. 

The proof of \eqref{flute2} now proceeds in a straightforward way,
applying the kind of manipulations encountered already in the proof of
\eqref{flute1} and is omitted.
\end{proof}
\vspace{\baselineskip}%
The remainder of the section describes, in Lie algebroid language,
some classical constructions related to torsion.  We will require the
additional assumption that ${\mathfrak g}\subset J^1{\mathfrak t}$ is
{\em surjective}.  We continue to denote the structure kernel of
${\mathfrak g}$ by ${\mathfrak h}$.

\subsection{Normalizing torsion and the upper coboundary morphism}
\lab{normalizing}%
Identify ${\mathfrak g}$ with ${\mathfrak t}\oplus{\mathfrak h}$ by
choosing a generator $\nabla$ of ${\mathfrak g}$. Then we obtain a
corresponding identification,
\begin{equation*}
        \alternating^2 ({\mathfrak g})\otimes{\mathfrak 
        t}\cong\Big(\,\alternating^2 ({\mathfrak t})\otimes{\mathfrak 
        t}\,\Big)\oplus\Big(\,{\mathfrak t}^*\otimes{\mathfrak 
        h}^*\otimes{\mathfrak t}\,\Big)\oplus\Big(\,\alternating^2 
        ({\mathfrak h})\otimes{\mathfrak t}\,\Big),
\end{equation*}
and a corresponding splitting of the torsion
\begin{equation*}
         da=\torsion\bar\nabla\oplus\operatorname{ev}\oplus\,0.
\end{equation*}
Here $\bar\nabla$ denotes the associated ${\mathfrak t}$-connection on 
${\mathfrak t}$ and $\operatorname{ev}$ is just evaluation, 
$\operatorname{ev}(V\otimes\phi):=\phi (V)$.  Notice that 
$\torsion\bar\nabla$ is the only component of $ da$ depending on the 
choice of generator.  Given two generators $\nabla^1$ and $\nabla^2$, 
their difference $\nabla^2-\nabla^1$ may be viewed as a section of 
${\mathfrak t}^*\otimes{\mathfrak h}$ and one readily computes
\begin{equation}
        \torsion\bar\nabla^2=\torsion\bar\nabla^1+\Delta 
        (\nabla^2-\nabla^1),\mathlab{deek}
\end{equation}
where $\Delta$ denotes the {\df  upper coboundary morphism}, defined 
as the composite
\begin{equation*}
        {\mathfrak t}^*\otimes{\mathfrak h}\monomorphism{\mathfrak 
        t}^*\otimes T^*\!M\otimes{\mathfrak 
        t}\xrightarrow{\identity\otimes\#^*\otimes\identity}{\mathfrak 
        t}^*\otimes{\mathfrak t}^*\otimes{\mathfrak 
        t}\xrightarrow{A\otimes\identity}\alternating^2 ({\mathfrak 
        t})\otimes{\mathfrak t}.
\end{equation*}
Here $\#^*\colon T^*\!M\rightarrow{\mathfrak t}^*$ is the dual of the
anchor $\#\colon{\mathfrak t}\rightarrow TM$ and $
A(\alpha\otimes\beta):=\alpha\wedge\beta$. As an elementary
consequence of \eqref{deek} above, we obtain:
\begin{proposition}
        If $ C\subset\alternating^2 ({\mathfrak 
        t})\otimes{\mathfrak t}$ is a complement for the image of $\Delta$, 
        then there  exists a generator $\nabla$ such that 
        $\torsion\bar\nabla\subset C$. If $\Delta$ is injective, then 
        this generator is unique.
\end{proposition}
\noindent Note that there is no need to require that $C$ be 
${\mathfrak g}$-invariant.
 
\subsection{Intrinsic torsion and torsion reduction}\lab{intrinsic}
Mimicking a classical construction, we define the {\df torsion 
bundle},
\begin{equation*}
         {H}({\mathfrak g}):=\frac{\alternating^2 ({\mathfrak 
         t})\otimes{\mathfrak t}}{\image\Delta},
\end{equation*}
and call the image $\tau$ of $\torsion\bar\nabla$, under the map
$\Gamma (\alternating^2 ({\mathfrak t})\otimes{\mathfrak
  t})\rightarrow \Gamma ({H} ({\mathfrak g}))$, induced by the
projection $\alternating^2 ({\mathfrak t})\otimes{\mathfrak
  t}\rightarrow H(\mathfrak{g})$, the {\df intrinsic torsion} of
${\mathfrak g}$.  By \eqref{deek}, $\tau$ is independent of the choice
of generator, i.e., is an invariant of ${\mathfrak g}$.  Since
$\Delta$ is ${\mathfrak g}$-equivariant, ${H}({\mathfrak g})$ is a
${\mathfrak g}$-representation whenever it is a bona fide vector
bundle (has constant rank).  In that case we can define the isotropy
${\mathfrak g}_\tau\subset{\mathfrak g}$ of $\tau$.  This is the {\df
  torsion reduction} of ${\mathfrak g}$ and is indeed a reduction, as
we show in \ref{qtax}.

\section{$\Theta$-reduction}\lab{sectionY}%
Let ${\mathfrak g}\subset J^1{\mathfrak t}$ be an infinitesimal
geometric structure with structure kernel ${\mathfrak h}$ and let
${\mathfrak g}^{(1)}\subset J^1{\mathfrak g}$ be its prolongation.
Recall that the {\df $\Theta$-reduction} of ${\mathfrak g}$ is simply
the image of ${\mathfrak g}^{(1)}$.   We will denote it
by ${{\mathfrak g}^{(1)}_1}$.  If ${\mathfrak g}^{(1)}$ and
${{\mathfrak g}^{(1)}_1}$ have constant rank then ${{\mathfrak
    g}^{(1)}_1}$ is evidently a subalgebroid of ${\mathfrak g}$.
Proposition \ref{prolongation} shows that the symmetries of
${\mathfrak g}$ are automatically symmetries of ${{\mathfrak
    g}^{(1)}_1}\subset{\mathfrak g}$.  So ${\mathfrak g}^{(1)}_1$ is
indeed a reduction, as claimed in Proposition \ref{koota}.

For our purposes we may as well assume that ${\mathfrak g}$ is
surjective; see \ref{hopeful}. For simplicity, however, we strengthen
this requirement:
\begin{assumption}
		 In this section ${\mathfrak g}\subset J^1{\mathfrak t}$ is a 
		 surjective infinitesimal geometric structure over a 
		 {\em transitive} Lie algebroid ${\mathfrak t}$.  In particular, 
		 ${\mathfrak g}$ is transitive.  Being a surjective, 
		 $\mathfrak g$ has a structure kernel ${\mathfrak h}$ of 
		 constant rank and   ${\mathfrak g}$ admits 
		 generators (Proposition \eqrefs{generator}{py1}).
\end{assumption}
\noindent Our chief objective is a characterization of the
$\Theta$-reduction ${{\mathfrak g}^{(1)}_1}$ that does not require an
explicit knowledge of ${\mathfrak g}^{(1)}$.

\subsection{The lower coboundary morphism}\lab{structure:kernel}%
As in torsion reduction, a `coboundary morphism' plays a central role 
in $\Theta$-reduction.  However, unless ${\mathfrak t}=TM$, the upper 
coboundary morphism $\Delta$,  defined in \ref{normalizing}, is not the 
appropriate one.  Rather, we need the {\df  lower coboundary morphism}  
$\delta$, defined as the composite
\begin{equation*}
                 T^*\!M\otimes{\mathfrak h}\monomorphism T^*\!M\otimes 
                 T^*\!M\otimes{\mathfrak t}\xrightarrow{A\otimes\identity}
                 \alternating^2 (TM)\otimes{\mathfrak t},
\end{equation*}
where $ A(\alpha\otimes\beta):=\alpha\wedge\beta $. This morphism is 
also a morphism of ${\mathfrak g}$-representations.

As we assume ${\mathfrak t}$ is transitive, we may, by dualizing the 
anchor map $\#\colon{\mathfrak t}\rightarrow TM$,  regard $ T^*\!M$ 
as a subbundle of ${\mathfrak t}^*$, and obtain natural inclusions
\begin{align*}
         T^*\!M\otimes{\mathfrak h}&\monomorphism{\mathfrak 
           t}^*\otimes{\mathfrak h}\\
         \alternating^2 (TM)\otimes{\mathfrak t}&\monomorphism\alternating^2 
           ({\mathfrak t})\otimes{\mathfrak t}.
\end{align*}
With this understanding, we may regard $\delta\colon 
T^*\!M\otimes{\mathfrak h}\rightarrow\alternating^2 
(TM)\otimes{\mathfrak t}$ as the restriction of the upper coboundary 
morphism $\Delta\colon{\mathfrak t}^*\otimes{\mathfrak 
h}\rightarrow\alternating^2 ({\mathfrak t})\otimes{\mathfrak t}$ 
defined in \ref{normalizing}.

The analogue of the torsion bundle $ H({\mathfrak g})$ defined in 
\ref{normalizing} is the (variable-rank) bundle
\begin{equation*}
                 {h}({\mathfrak g}):=\frac{\alternating^2 
                 (TM)\otimes{\mathfrak t}}{\image\delta}.
\end{equation*}
Whenever $ h({\mathfrak g})$ is a genuine vector bundle (has constant 
rank), it is a ${\mathfrak g}$-representation. There is evidently a 
natural morphism $\psi\colon h({\mathfrak g})\rightarrow H({\mathfrak 
g})$ making the following diagram commute:
\begin{equation}
   \begin{CD}
      \alternating^2 (TM)\otimes{\mathfrak t}@>/\image\delta>>h({\mathfrak g})\\
      @V\text{inclusion}VV @VV\psi V\\
      \alternating^2 ({\mathfrak t})\otimes{\mathfrak t}@>>/\image\Delta> 
         H({\mathfrak g})
   \end{CD}\enspace.\mathlab{dg7}
\end{equation}
Note that $\psi$ need not be injective.

\subsection{The characterization of ${\mathfrak 
g}^{(1)}_1$}\lab{qlook}%
The significance of the bundle $ h({\mathfrak g})$ is the existence of
a natural morphism $\Theta\colon{\mathfrak g}\rightarrow h({\mathfrak
  g})$ such that ${\mathfrak g}^{(1)}_1=\kernel\Theta$.  (This is the
origin of our terminology `$\Theta$-reduction.')  Recalling that
${\mathfrak g}^{(1)}:=(J^1{\mathfrak g})_a\cap (J^1{\mathfrak
  g})_{da}$, our construction of $\Theta$ begins with the following
observation:
\begin{conditions}
        \item\lab{z1z} The isotropy $(J^1{\mathfrak g})_a\subset 
        J^1{\mathfrak g}$ is a surjective infinitesimal geometric 
        structure with structure kernel $ T^*\!M\otimes{\mathfrak h}$.
\end{conditions}
\begin{proof}
It is not difficult to see that $ (J^1{\mathfrak g})_a$ is the kernel 
of
\begin{align*}
     \xi&\mapsto p\xi-(J^1 a)\xi\\
     J^1{\mathfrak g}&\rightarrow T^*\!M\otimes{\mathfrak t},
\end{align*}
where $ p\colon J^1{\mathfrak g}\rightarrow{\mathfrak g}$ is the 
projection.  This follows from the transitivity of ${\mathfrak g}$ 
and, e.g., \eqrefs{tautological}{ballet}.  One establishes \eqref{z1z} 
by applying Lemma \ref{algebraicLemma} to this morphism.
\end{proof}
Now ${\mathfrak g}^{(1)}$ is the kernel of the morphism 
$\xi\mapsto\xi\cdot da\colon (J^1{\mathfrak 
g})_a\rightarrow\alternating^2 ({\mathfrak g})\otimes{\mathfrak t}$. 
However, it follows from \eqrefs{tautological}{flute2} and transitivity that:
\begin{conditions}
                \item \lab{m1}For any $\xi\in (J^1{\mathfrak g})_a$, the 
                element $\xi\cdot da\in\alternating^2 ({\mathfrak 
                g})\otimes{\mathfrak t}$ is tensorial --- i.e., drops to an 
                element $(\xi\cdot da)^\vee\in\alternating^2 
                (TM)\otimes{\mathfrak t}$.
\end{conditions}
This means we may regard ${\mathfrak g}^{(1)}$ as the kernel of a 
morphism
\begin{align*}
        (J^1{\mathfrak g})_a&\xrightarrow{\theta}\alternating^2 
        (TM)\otimes{\mathfrak t}\\
        \xi&\mapsto (\xi\cdot da)^\vee.
\end{align*}
According to \eqref{z1z}, the domain of $\theta$ fits into an  exact 
sequence,
\begin{equation*}
         0\rightarrow T^*\!M\otimes{\mathfrak h}\monomorphism  (J^1{\mathfrak 
         g})_a\rightarrow{\mathfrak g}\rightarrow 0.
\end{equation*}
Applying Lemma \ref{algebraicLemma} to the morphism 
$\theta$, we obtain a corresponding exact sequence,
\begin{equation*}
         0\rightarrow\kernel\delta\monomorphism{\mathfrak 
         g}^{(1)}\xrightarrow{a^{(1)}}\kernel\Theta\rightarrow 0,
\end{equation*}
where $\Theta$ is the unique morphism making the following diagram 
commute:
\begin{equation}
                        \begin{CD}
                                (J^1{\mathfrak g})_a @>>> {\mathfrak g}\\
                                @VV\theta V @VV\Theta V \\
                     \alternating^2 (TM)\otimes{\mathfrak 
                     t}@>/\image\delta>> 
                     {h({\mathfrak g})}
                        \end{CD}\enspace.\mathlab{diag}
\end{equation}
Summarizing:
\begin{proposition}
        If ${\mathfrak g}\subset J^1{\mathfrak t}$ is surjective and 
        ${\mathfrak t}$ is transitive, then there is a natural morphism 
        $\Theta\colon{\mathfrak g}\rightarrow h({\mathfrak g})$, 
        constructed above, such that
    \begin{equation*}
         0\rightarrow\kernel\delta\monomorphism{\mathfrak 
         g}^{(1)}\xrightarrow{ a^{(1)}}{\mathfrak 
         g}\xrightarrow{\Theta}h({\mathfrak g})
    \end{equation*}
        is exact.  In particular, the structure kernel of ${\mathfrak 
        g}^{(1)}$ is the kernel of the lower coboundary morphism $\delta$, 
        while the image ${\mathfrak g}^{(1)}_1$ of ${\mathfrak g}^{(1)}$ 
        (the $\Theta$-reduction of ${\mathfrak g}$) is the kernel of 
        $\Theta$. If $\kernel\delta$ and $\kernel\Theta$ have 
        constant rank then ${\mathfrak g}^{(1)}\subset J^1{\mathfrak g}$ 
        is an infinitesimal geometric structure.
\end{proposition}
\begin{remark}
  By the proposition the structure kernel of ${{\mathfrak g}^{(1)}} $
  lies entirely within $T^*\!M \otimes {\mathfrak h} \subset T^*\!M
  \otimes {\mathfrak g} $ and is consequently {\em commutative}. In
  concrete calculations it is often useful to think of $\kernel\delta$
  as the collection of all $\sigma\in\symmetric^2
  (TM)\otimes{\mathfrak t}$ such that $\sigma (U,\,\cdot\,)\subset
  T^*\!M\otimes{\mathfrak t}$ lies in ${\mathfrak h}$ for all $ U\in
  TM$.
\end{remark}
\subsection{The relationship with torsion}\lab{qtax}
By construction $\Theta$ is an {\em invariant} of ${\mathfrak g}$.  
However, the construction of $\Theta$ given here is not immediately 
useful in computations.  In Sect.\,\ref{method} we describe 
$\Theta\colon{\mathfrak g}\rightarrow h({\mathfrak g})$ explicitly in 
terms of a generator $\nabla$ of ${\mathfrak g}$.  As a byproduct, we 
will obtain a proof of the following link between $\Theta$-reduction 
and intrinsic torsion:
\begin{theorem}
  If ${\mathfrak g}\subset J^1{\mathfrak t}$ is surjective,
  ${\mathfrak t}$ is transitive, and $H(\mathfrak g)$ has constant
  rank, then $\psi (\Theta (X))=X\cdot\tau$.
\end{theorem}
\noindent Here $\psi$ is the natural morphism $\psi\colon h({\mathfrak 
g})\rightarrow H({\mathfrak g})$ defined in \ref{structure:kernel}, 
and $\tau\subset H({\mathfrak g})$ is the intrinsic torsion, defined 
in \ref{intrinsic}. As a corollary we obtain the following result 
showing that torsion reduction is generally cruder than 
$\Theta$-reduction:
\begin{corollary}
         Suppose that the torsion bundle ${H}({\mathfrak g})$ has constant 
         rank, so that the torsion reduction ${\mathfrak g}_\tau$ of 
         ${\mathfrak g}$ is well-defined.  Then ${\mathfrak 
         g}^{(1)}_1\subset{\mathfrak g}_\tau$.  In particular, if the 
         lower coboundary morphism $\delta$ has constant rank, and both 
         ${\mathfrak g}_\tau$ and ${\mathfrak g}^{(1)}_1$ have constant 
         rank, then ${\mathfrak g}_\tau$ is a reduction of ${\mathfrak g}$ 
         in the sense of \ref{reduction}.  If ${\mathfrak t}=TM$ then 
         $\Theta$-reduction and torsion reduction coincide.
\end{corollary}
\noindent Here the rank hypotheses and Proposition \ref{qlook} ensure 
that ${\mathfrak g}^{(1)}$ has constant rank, so that Proposition 
\ref{koota} applies.  However, the result presumably holds with a 
constant rank hypothesis on ${\mathfrak g}_\tau$ alone.

\subsection{Structures both surjective and $\Theta$-reduced}
\lab{goosey}%
We say that ${\mathfrak g}\subset J^1{\mathfrak t}$ is {\df
  $\Theta$-reduced} if it coincides with its $\Theta$-reduction.
\begin{theorem}
  Let ${\mathfrak g}\subset J^1{\mathfrak t}$ be a surjective 
  infinitesimal geometric structure on a transitive Lie algebroid 
  ${\mathfrak t}$.  Assume that ${\mathfrak g}$ is $\Theta$-reduced 
  (equivalently, that the map $\Theta$ defined above vanishes).  
  Assume that the associated lower coboundary morphism $\delta$ is 
  injective.  Then ${\mathfrak g}$ has an associated Cartan algebroid, 
  namely ${\mathfrak g}$ itself, equipped with a canonical Cartan 
  connection $\nabla^{(1)}$.  The $\nabla^{(1)}$-parallel sections of 
  ${\mathfrak g}$ coincide with the prolonged symmetries of 
  ${\mathfrak g}$.
\end{theorem}
\begin{proof}
Proposition \ref{koota} implies the prolongation ${\mathfrak 
g}^{(1)}$ of ${\mathfrak g}$ is surjective.  In  addition,  
${\mathfrak g}^{(1)}$ has trivial structure kernel, because we 
suppose $\delta$ is injective (Proposition \ref{qlook}).  Applying 
Theorem \ref{neoGenerator} to the infinitesimal geometric structure 
${\mathfrak g}^{(1)}$, we obtain a Cartan connection $\nabla^{(1)}$ on 
${\mathfrak g}$ whose parallel sections are the symmetries of 
${\mathfrak g}^{(1)}$.  These are nothing but the {\em prolonged} 
symmetries of ${\mathfrak g}$, by Proposition \ref{prolongation}.
\end{proof}

In Proposition \ref{natural} we characterize $\nabla^{(1)}$ as the
unique `natural' connection on ${\mathfrak g}$ whose curvature
$\curvature\nabla^{(1)}\subset\alternating^2 (TM)\otimes{\mathfrak
  g}^*\otimes{\mathfrak g}$ takes values in ${\mathfrak
  h}\subset\mathfrak g$. A general formula expressing $\nabla^{(1)}$
in terms of a generator of ${\mathfrak g}$ will appear in \ref{jit}.  

\subsection{The special case ${\mathfrak t}=TM$}\lab{goosey2}
When ${\mathfrak t}=TM$, $\Theta$-reduction and torsion reduction are
the same thing, as are the upper and lower coboundary morphisms,
$\delta$ and $\Delta$.  We now rewrite the above theorem accordingly,
adding explicit information about the Cartan connection that we
establish later in \ref{caviar}.

Here $\bar\nabla$ will denote the dual of $\nabla$, i.e.,
$\bar\nabla_UV=\nabla_VU+[U,V]$.  We call $\mathfrak g\subset J^1(TM)$
{\df reductive} if $\Delta$ has constant rank and if the image of
$\Delta$ admits a ${\mathfrak g}$-invariant complement $ C$.  We call
the generator $\nabla$ {\df normal} if $\torsion\bar\nabla\subset C$
for some such $ C$.  Proposition \ref{normalizing} guarantees the
existence of normal generators when $\mathfrak g$ is reductive.

We call an arbitrary structure $\mathfrak g\subset\mathfrak t$ {\df
  torsion-reduced} if ${\mathfrak g}={\mathfrak g}_\tau$, i.e., if the
intrinsic torsion $\tau$ is $\mathfrak g$-invariant.

\begin{theorem}
  Let ${\mathfrak g}\subset J^1 (TM)$ be a transitive, torsion-reduced
  infinitesimal geometric structure, and suppose that the associated
  upper coboundary morphism $\Delta$ is injective.  Then ${\mathfrak
    g}$ has an associated Cartan algebroid, namely ${\mathfrak g}$
  itself, equipped with a Cartan connection $\nabla^{(1)}$ described
  below.  The $\nabla^{(1)}$-parallel sections of ${\mathfrak g}$
  coincide with the prolonged symmetries of ${\mathfrak g}$.
                
  Choose a generator $\nabla$ for ${\mathfrak g}$ (equivalently,
  choose a complement for the image of $\Delta$; see Proposition
  \ref{normalizing}).  Then:
         \begin{conditions}
                \item\lab{qd1} The equation
                \begin{equation*}
                        \Delta (\epsilon 
                        (V\oplus\phi))=\bar\nabla_V\torsion\bar\nabla
                        +\phi\cdot\torsion\bar\nabla
                        \quad\text{(\,$V\subset TM$, 
                          $\phi\subset{\mathfrak h}$ arbitrary)}
                \end{equation*}
                has a unique solution morphism $\epsilon\colon
                TM\oplus{\mathfrak h}\rightarrow
                T^*\!M\otimes{\mathfrak h}$.
                
              \item\lab{scar} Identifying ${\mathfrak g}$ with $
                TM\oplus{\mathfrak h}$ using the generator, we have
                \begin{equation*}
                        \nabla_U^{(1)} (V\oplus\phi)=(\nabla_UV+\phi (U))\oplus 
                        (\bar\nabla_U\phi+\epsilon 
                        (V\oplus\phi)U+\curvature\bar\nabla (U, V)).
                \end{equation*}
                
              \item\lab{qd3} If ${\mathfrak g}$ is reductive and
                $\nabla$ is normal, or if $\tau =0$ and $\nabla$ is
                torsion-free, then $\epsilon=0$ for any normal
                generator $\nabla$.
         \end{conditions}
\end{theorem}

When one of the conditions in \eqref{qd3} holds, obstructions to
symmetry are particularly simple to describe, as is the symmetry Lie
algebra ${\mathfrak g}_0$ in the globally flat case.  Indeed, one then
computes, with the help of \eqrefs{reconstruction}{jj} and the Bianchi
identity $d_{\bar\nabla}\curvature \bar\nabla =0$,
\begin{gather*}
  \curvature\nabla^{(1)} (U_1,U_2)(V\oplus\phi)=0\oplus
  \Big(\,-(\bar\nabla_V\curvature\bar\nabla+\phi\cdot\curvature\bar\nabla)
  (U_1,U_2)\,\Big),\\
  \torsion\overline{\nabla^{(1)}} (V_1\oplus\phi_1,
  V_2\oplus\phi_2)=\\
  (\,\torsion \bar\nabla(V_1, V_2)+\phi_1(V_2)-\phi_2(V_1)\,) \oplus
  (\,[\phi_1,\phi_2]_{\mathfrak h}-\curvature\bar\nabla (V_1, V_2)\,).
\end{gather*}
Here $\overline{\nabla^{(1)}}$ denotes the representation of 
${\mathfrak g}$ on itself associated with the Cartan connection $\nabla^{(1)}$ 
on ${\mathfrak g}$. Applying Theorem \ref{dissymmetry}:
\begin{corollary}
  Let ${\mathfrak g}\subset J^1 (TM)$ be an infinitesimal geometric
  structure satisfying the hypotheses of the above theorem, and assume
  either that ${\mathfrak g}$ is reductive and $\nabla$ is normal, or
  that $\tau =0$ and $\nabla$ is torsion-free.  Let \,${\mathcal
    U}\subset M$ be an arbitrary open set and ${\mathfrak g}_0$ be the
  Lie algebra of all symmetries of ${\mathfrak g}|\,{\mathcal U}$.
  Then $\dimension{\mathfrak g}_0\le\rank{\mathfrak g}$. If $U$ is
  simply-connected, then equality holds if and only if
  $\curvature\bar\nabla$ is both ${\mathfrak h}$-invariant and
  $\bar\nabla$-parallel.  In that case ${\mathfrak g}_0$ is naturally
  isomorphic to $ T_mM\oplus{\mathfrak h}(m)$ ($ m\in {\mathcal U}$
  arbitrary) with Lie bracket given by
        \begin{multline*}
          [V_1\oplus\phi_1,V_2\oplus\phi_2]\\=(\,\torsion
          \bar\nabla(V_1, V_2)+\phi_1(V_2)-\phi_2(V_1)\,)\oplus
          (\,[\phi_1,\phi_2]_{\mathfrak h}-\curvature\bar\nabla (V_1,
          V_2)\,).
        \end{multline*}
\end{corollary}

\subsection{The symmetries of Riemannian structures}
\lab{illustration} Let $\mathfrak g\subset J^1 (TM)$ denote the bundle
of 1-symmetries of a Riemannian metric $\sigma$, as described in
detail in \ref{Riemannian}.  The upper coboundary 
morphism for ${\mathfrak g}$ is a map
\begin{equation*}
        T^*\!M \otimes{\mathfrak h}\xrightarrow{\Delta}\alternating^2 
        (TM)\otimes TM,
\end{equation*}
where ${\mathfrak h}\subset T^*\!M\otimes TM$ is the
$\mathfrak{o}(n)$-bundle of all skew-symmetric tangent space
endomorphisms.  This morphism is well known to be an isomorphism.  For
example, using the metric $\sigma$ to make various identifications, we
may view $\Delta$ as the map
\begin{align*}
        T^*\!M\otimes\alternating^2 (TM)&\rightarrow\alternating^2 
        (TM)\otimes T^*\!M\\
        \alpha\otimes (\beta_1\wedge\beta_2)&\mapsto 
        (\alpha\wedge\beta_1)\otimes\beta_2-(\alpha\wedge\beta_2)\otimes\beta_1,
\end{align*}
which has explicit inverse
\begin{equation*}
         (\beta_1\wedge\beta_2)\otimes\alpha\mapsto\frac{ 
         1}{2}\Big(\,\beta_1\otimes (\beta_2\wedge\alpha)-\beta_2\otimes 
         (\beta_1\wedge\alpha)-\alpha\otimes (\beta_1\wedge\beta_2)\,\Big).
\end{equation*}

Since $\Delta$ is an isomorphism, $\mathfrak g$ has, by Proposition 
\ref{normalizing}, a unique torsion-free generator $\nabla$; this is 
the Levi-Cevita connection.  The intrinsic torsion vanishes because 
$H(\mathfrak g)=0$, making $\mathfrak g$ torsion-reduced.  Also, 
$\mathfrak g$ is trivially reductive.  Applying Theorem \ref{goosey2}, 
we obtain the Cartan connection $\nabla^{(1)}$ and related claims
in \ref{cakeWalk}.  (We have $\bar\nabla=\nabla$ since 
$\torsion\nabla=0$.)

According to Corollary \ref{goosey2}, we are in the maximally
symmetric case when $\curvature\nabla$ is both ${\mathfrak
  h}$-invariant and $\nabla$-parallel.  According to a well-known
representation-theoretic analysis of the curvature module, this
happens if and only if
\begin{equation*}
        \curvature\nabla (V_1, V_2)=s\Big(\,\sigma (V_1)\otimes V_2-\sigma 
        (V_2)\otimes V_1\,\Big);\qquad V_1,V_2\in TM,
\end{equation*}
for some constant $ s\in{\mathbb R}$ (the scalar curvature). The Lie
algebra $\mathfrak g_0$ of symmetries described in the corollary is
then isomorphic to the Lie algebra of infinitesimal isometries of
Euclidean space, hyperbolic space, or the sphere, according to whether
$s=0$, $s<0$, or $ s>0$.

\subsection{The symmetries of a conformal parallelism}\lab{conformalParallelism}
Let $\omega \colon TM \rightarrow V$ be global parallelism ($V$ a
vector space with the dimension of $M$) and let $\langle\omega\rangle
\subset T^*\!M \otimes V$ be the line bundle spanned by $\omega$. A
{\df conformal parallelism} is an equivalence class of absolute
parallelisms, where $\omega,\omega' \colon TM \rightarrow V$ are
considered equivalent if $\omega'=f \omega $ for some positive
function $f$. The infinitesimal isometries of the conformal
parallelism having $\omega$ as representative coincide with the
symmetries of the isotropy ${\mathfrak g} \subset J^1(TM) $ of
$\langle\omega\rangle \subset T^*\!M \otimes V$.

A straightforward application of Lemma \ref{algebraicLemma} shows that
${\mathfrak g} $ is transitive, with rank-one structure kernel
$\langle\identity\rangle \subset T^*\!M \otimes TM $. 

The upper boundary morphism, given by
\begin{gather*}
  \Delta \colon T^*\!M \rightarrow \alternating^2(TM)\otimes TM\\
  \Delta(\beta)(U_1,U_2)=\beta(U_1)U_2-\beta(U_2)U_1,
\end{gather*}
is evidently injective ($\dimension M \ge 2$). We leave it to
the reader to verify that $\mathfrak g$ has vanishing intrinsic
torsion $\tau$ precisely when
\begin{equation*}
  d \omega= \alpha \wedge \omega,
\end{equation*}
for some one-form $\alpha$. While $\alpha$ depends on the choice of representative $\omega$, the two-form $d \alpha $ does not. 

Assuming $\tau=0$, $\mathfrak g$ has a unique torsion-free generator
$\nabla$ (by \eqrefs{normalizing}{deek} and the definition of
$\tau$). Moreover, it is not hard to show that
\begin{equation*}
\bar\nabla_U \omega = \alpha(U) \omega, \qquad U \subset TM,
\end{equation*}
and accordingly that
\begin{equation*}
  \curvature \bar\nabla=d \alpha \otimes \identity.
\end{equation*}

Applying Corollary \ref{goosey2}, we are in the maximally symmetric case when $d \alpha=0$. The Lie algebra ${\mathfrak g}_0$ is then isomorphic to the semi-direct product of the commutative Lie algebras $T_mM$ and ${\mathbb R}$, with ${\mathbb R} $ acting on $T_mM$ by scalar multiplication.

\section{Application: subriemannian contact three-manifolds}%
\lab{contact}
This section is an extended application of the general theory
developed in the two preceding sections. It also includes a concrete
construction of invariant differential operators, from the Lie
algebroid representation viewpoint.

Equip a three-dimensional manifold $M$ with a rank-two subriemannian
structure, i.e., a rank-two distribution $\mathcal{H}\subset TM $,
together with a real fibre bundle inner product $\sigma \subset
\symmetric^2(\mathcal{H})$ on $\mathcal{H}$. We assume that
$\mathcal{H}$ is a {\em contact} distribution (see below); so equipped,
$M$ becomes a {\df subriemannian contact three-manifold}. This section
describes the application of Cartan's method to such manifolds, \`a la
Lie algebroids. For the $G$-structure approach, see K. Hughen's thesis
\cite{Hughen_95} or \cite{Montgomery_02}.

Here we shall understand $\mathcal{H}$ to be transversally orientable,
and understand the specification of the subriemannian structure to
include a specification of transverse orientation. This amounts to the
choice of a non-unique, non-vanishing one-form $\theta$
annihilating $\mathcal{H}$. The contact hypothesis means that $\theta$
is contact, i.e., $d\theta$ restricts to a symplectic structure
on $\mathcal{H}$.

The infinitesimal isometries of the subriemannian contact structure
are the symmetries of the infinitesimal geometric structure $J^1
(TM)_{\mathcal{H},\sigma}$, where $J^1(TM)_{\mathcal{H} }\subset
J^1(TM)$ is the isotropy of $\mathcal{H}$ and
$J^1(TM)_{\mathcal{H},\sigma}\subset J^1(TM)_{\mathcal{H}}$ is the
isotropy of $\sigma$; see \ref{isotropy}.

\subsection{Preliminary reduction}\lab{prered}%
The symplectic structure $d\theta|\mathcal{H}$ orients
$\mathcal{H}$. Consequently, there is a well defined area form $dA$
determined by the subriemannian metric $\sigma$ on $\mathcal{H}$. In
fact, rescaling $\theta$ by a positive function if necessary, we may
arrange $dA=d\theta| \mathcal{H}$. A contact form $\theta $ normalized
in this way is evidently an invariant of the subriemannian contact
structure, implying that the isotropy
\begin{equation*}
  {\mathfrak g}:= J^1(TM)_{{\mathcal H},\sigma,d \theta}\subset J^1(TM)_{{\mathcal H},\sigma}
\end{equation*}
of $d \theta $ is a reduction of $J^1(TM)_{{\mathcal
    H},\sigma}$. (This reduction is in fact the first
torsion-reduction of $J^1(TM)_{{\mathcal H},\sigma}$.) 

The subriemannian metric $\sigma$ has a canonical extension to a bona
fide Riemannian metric, defined as follows: Let $\mathbf n$ be the
{\df Reeb vector field} associated with the normalized contact form
$\theta$. That is,
\begin{equation*}
  d\theta(\mathbf{n},\,\cdot\,)=0,\qquad \theta(\mathbf{n})=1. 
\end{equation*}
One extends $\sigma$ so as to make $\mathbf n$ orthogonal to
$\mathcal{H}$ and have unit length; then
$\theta=\sigma(\mathbf{n}):=\sigma(\mathbf{n},\,\cdot\,)$. The easy
proof of the following is left to the reader.
\begin{proposition}
  The reduction $\mathfrak g\subset J^1(TM) $ above coincides with the
  joint isotropy of the extended metric $\sigma$ and $\mathbf{n}$,\enspace
  ${\mathfrak g}=J^1(TM)_{\sigma,{\mathbf n}}$.
\end{proposition}
% \noindent In particular, a linear connection $\nabla $ generates
% $\mathfrak g$ if and only if $\bar\nabla \sigma =0$ and $\bar\nabla
% {\mathbf n}=0$, where $\bar\nabla_UV=\nabla_VU+[U,V]$.

\subsection{The complex structure on ${\mathcal H} $}\lab{cr}
Let $\times$ denote the usual cross product determined by the extended
metric $\sigma$ and define $J\subset T^*\!M\otimes TM$ by
$JU=\mathbf{n}\times U$. Then $J$ has kernel $\langle {\mathbf
  n}\rangle$, image $\mathcal{H}$, and the restriction of $J$ to
$\mathcal{H}$ is the complex structure on $\mathcal{H}$ relating the
area form $dA$ to the subriemannian metric $\sigma$:
\begin{equation*}
  dA(U_1,U_2)=\sigma(JU_1,U_2);\qquad U_1,U_2\subset \mathcal{H}.
\end{equation*}
\begin{proposition}%
  $\mathfrak{g}\subset J^1(TM)$ is a surjective infinitesimal geometric
  structure whose structure kernel $\mathfrak{h}\subset T^*\!M\otimes
  TM$ is the globally trivial $\mathfrak{o}(2)$-bundle spanned by the
  $\mathfrak{g}$-invariant tensor $J$.
\end{proposition}
\begin{proof}
  Combine the characterization of $\mathfrak{g}$ in Proposition
  \ref{prered} with the general observations of \ref{vectorField}
  (with $V:=\mathbf{n}$). In particular, the surjectivity of
  $\mathfrak{g}$ follows from Proposition
  \ref{vectorField}. Alternatively, one may directly apply Lemma
  \ref{algebraicLemma}, in the manner already demonstrated  several times.
\end{proof}

\subsection{Normalizing torsion}\lab{normtor}%
By Proposition \ref{prered}, a linear connection $\nabla$ on
$TM$ generates $\mathfrak{g}$ if and only if ${\bar\nabla} \sigma=0$
and ${\bar\nabla} \mathbf{n}=0$; here $\sigma $ denotes the extended
metric and ${\bar\nabla}_UV:=\nabla_VU +[U,V]$. Following the general
discussion of \ref{normalizing}, we now fix a natural choice of
generator. As a byproduct, we obtain a concrete expression for the
intrinsic torsion $\tau \subset H(\mathfrak{g})$.

Let $({\mathcal H}^*\otimes {\mathcal H})_\sym \subset {\mathcal
  H}^* \otimes {\mathcal H} $ denote the rank-three subbundle consisting of
endomorphisms of fibres of ${\mathcal H} $ that are symmetric with
respect to the metric $\sigma $.
\begin{proposition}\mbox{}%
\begin{conditions}
\item\lab{mice1}$\mathfrak{g}\subset J^1(TM)$ is reductive, in the
  sense of \ref{goosey2}, and has injective upper coboundary morphism
  $\Delta \colon T^*\!M \otimes {\mathfrak h} \rightarrow
  \alternating^2(TM)\otimes TM$.

\item\lab{mice2} For any generator $\nabla $ of $\mathfrak{g}$, and
  all vector fields $V \subset TM$, we have $\nabla_\mathbf{n}
  \mathbf{n} =0$ and $\nabla_V \mathbf{n}\subset \mathcal{H} $,
  allowing us to view $\nabla \mathbf{n} $ as a section of
  $\mathcal{H}^*\otimes \mathcal{H}$.

\item\lab{mice3}There exists a unique and normal generator $\nabla$
  such that
  \begin{equation*}
    \nabla\mathbf{n}\subset ({\mathcal H}^*\otimes {\mathcal H})_\sym\quad\text{and}\quad
    \nabla \sigma | \mathcal{H}=0.
  \end{equation*}
  Here $\nabla \sigma | \mathcal{H} \subset \mathcal{H}^*\otimes
  \symmetric^2(\mathcal{H})$ denotes the restriction of \,$\nabla \sigma
  \subset T^*\!M\otimes \symmetric^2(TM)$.
\end{conditions}
With $\nabla $ so fixed, we have:
\begin{conditions}
  \item\lab{mice3b}The torsion $\torsion {\bar\nabla}=-\torsion \nabla $ is given by the formula,
  \begin{equation*}
      \torsion {\bar\nabla}(U_1+a_1 \mathbf{n},U_2+a_2\mathbf{n})
      =(a_1 \nabla_{U_2}\mathbf{n} -a_2 \nabla_{U_1}\mathbf{n}) 
      + dA(U_1,U_2)\mathbf{n}.
  \end{equation*}
  Here $U_1,U_2 \in \mathcal{H},\,a_1,a_2\in {\mathbb R}$.

  \item\lab{mice4}There exists a natural isomorphism of $\mathfrak{g}$-representations,
    \begin{equation*}
      H(\mathfrak{g})\cong \alternating^2 (TM)\oplus({\mathcal H}^*\otimes {\mathcal H})_\sym,
    \end{equation*}
  with respect to which the intrinsic torsion of $\mathfrak{g}$ takes
  the form,
  \begin{equation*}
   \tau=d\theta \oplus \nabla \mathbf{n}.
  \end{equation*}

\item\lab{mice5}The intrinsic torsion component $\nabla \mathbf{n}$
  can be interpreted as a transverse derivative of the subriemannian
  metric $\sigma$:
  \begin{equation*}
    \sigma(\nabla_{U_1}\mathbf{n},U_2)=(\nabla_\mathbf{n}\sigma) (U_1,U_2);
    \qquad U_1,U_2 \subset \mathcal{H}.    
  \end{equation*}
\end{conditions}
\end{proposition}
The proposition is established by analyzing the coboundary morphism
$\Delta $ in detail, identifying a natural $\mathfrak g$-invariant
complement for its image, and applying Proposition
\ref{normalizing}. This analysis is not hard but a little tedious, and
is relegated to Appendix \ref{possum}. For the interested reader, we include
there a formula for the normalized generator $\nabla $ in terms of the
Levi-Cevita connection associated with the extended metric $\sigma$.

\subsection{Bianchi Identities and low weight differential operators}%
\lab{weight}
We shall write down Bianchi identities for the normalized generator
$\nabla $ using certain natural invariant differential operators, the
systematic construction of which is deferred to \ref{reps}. For now,
we merely list those operators of `weight' two or less.

Let ${\mathcal F}(M)=\Gamma({\mathbb R}\times M)$ denote the space of
all smooth functions. Then we define $\partial_{\mathbf n}\colon {\mathcal
  F}(M)\rightarrow {\mathcal F}(M)$, $\gradient_{\mathcal H}\colon
{\mathcal F}(M)\rightarrow \Gamma({\mathcal H})$, and
$\divergence_{\mathcal H}\colon \Gamma({\mathcal H})\rightarrow
{\mathcal F}(M)$ in a familiar way:
\begin{gather*}
  \partial_{\mathbf n} f:=df({\mathbf n}),\\
  \gradient_{\mathcal H}f:=\sigma^{-1}(df|{\mathcal
    H}),\quad\divergence_{\mathcal H} U := \trace(\bar\nabla U|{\mathcal H}).
\end{gather*}
Additionally, it is convenient to write $\curl_{\mathcal H} U:=\divergence_{\mathcal H}(JU)$.

Next, let ${\mathcal H}_2 \subset ({\mathcal H}^*\otimes {\mathcal
  H})_\sym$ denote the rank-two subbundle of trace-free
elements. Equivalently, ${\mathcal H}_2$ consists of those elements of
${\mathcal H}^*\otimes {\mathcal H}$ that are ${\mathbb
  C}$-antilinear: ${\mathcal H}_2=\bar {\mathcal H}^{*{\mathbb
    C}}\otimes_{\mathbb C} {\mathcal H}$. Like ${\mathcal H}$,
${\mathcal H}_2$ is a ${\mathfrak g}$-representation and a complex
line-bundle. However, while $\adjoint_JU=iU$, for $U \subset {\mathcal
  H}$ (by definition), we have, for $q \subset {\mathcal H}_2$,
$(\adjoint_Jq)U=J(qU)-q(JU)=2J(qU)$, i.e., $\adjoint_Jq=2iq$.

We are now ready to define two invariant operators $\partial_+\colon
\Gamma({\mathcal H})\rightarrow \Gamma({\mathcal H}_2)$ and
$\partial_-\colon \Gamma({\mathcal H}_2)\rightarrow
\Gamma({\mathcal H})$ according to 
\begin{gather*}
  (\partial_+U)V=\frac{1}{2}\Big(\,\bar\nabla_VU+J \bar\nabla_{JV}U\,\Big),\\
  (\bar\nabla_{U_1}q)U_2-(\bar\nabla_{U_2}q)U_1=dA(U_1,U_2)\,\partial_-q.
\end{gather*}

Associated with the normalized generator $\nabla $ of $\mathfrak g$
are its two fundamental invariants $T := \torsion \bar\nabla $ and
$\Omega := \cocurvature \nabla =-\curvature \bar\nabla $, which
satisfy the Bianchi identities \eqrefs{Bianchi}{BianchiI} and
\eqrefs{Bianchi}{BianchiII}. Of course these are also invariants of
the subriemannian contact structure. According to
\eqrefs{normtor}{mice3b}, $T$ depends only on the previously
identified invariant $\nabla {\mathbf n} $. As it turns out, one
component of $\Omega$ is a new invariant function. Recalling that
$\Omega(U_1, U_2)\subset {\mathfrak h} $ for all $U_ 1, U_2 \subset
TM$ (Proposition \eqrefs{reconstruction}{ala}) and that $dA$ spans
$\alternating^2({\mathcal H})$, there is a real-valued function $\kappa
$ well defined by
\begin{equation}
  \Omega(U_1,U_2) U_3 =-\kappa\,dA (U_1,U_2)JU_3;\qquad U_1,U_2,U_3 
  \subset {\mathcal H}.\mathlab{lemonade}
\end{equation}
\begin{proposition}[Bianchi identities]\mbox{}
\begin{conditions}
\item\lab{dooda} $\trace(\nabla {\mathbf n})=0$, i.e., $\nabla
  {\mathbf n} \subset {\mathcal H}_2$.
\item\lab{sweet}  $\partial_{\mathbf n} \kappa=-\frac{1}{2}\curl_{\mathcal H}(\partial_-(\nabla {\mathbf n}))$.
\item\lab{tootsweet} The cocurvature of $\nabla $ is given by
\begin{multline*}
  \Omega(U_ 1+ a_1 {\mathbf n} ,U_2 +a_ 2 {\mathbf n} )(U_3+ a_3 {\mathbf n} )=\\
  \Big(-\kappa \,dA(U_1,U_2) +\frac{1}{2}\sigma \big(\partial _- ( \nabla  {\mathbf n}  ),\, a _1 U_2 -a_2U_1 ) \big)\,\Big) J U_3;\\U_ 1,U_2,U_3 \in {\mathcal H},\,a_ 1, a _2, a_ 3\in {\mathbb R}.
\end{multline*}
\end{conditions}
\end{proposition}
\begin{proof}
Proposition \eqrefs{normtor}{mice3b} states that 
\begin{equation*}
  T(U_1+a_1 \mathbf{n},U_2+a_2\mathbf{n})
      =(a_1 \nabla_{U_2}\mathbf{n} -a_2 \nabla_{U_1}\mathbf{n}) 
      + dA(U_1,U_2)\mathbf{n}.
\end{equation*}
A little multilinear algebra determines that  $\Omega $ has the general form 
\begin{equation*}
  \Omega(U_ 1+ a_1 {\mathbf n} ,U_2 +a_ 2 {\mathbf n} )(U_3+ a_3 {\mathbf n} )=-\Big(\,\kappa \,dA(U_1,U_2) +\omega(a _1 U_2 -a_2U_1 )\,\Big) J U_3,
\end{equation*}
for some section $\omega \subset {\mathcal H}^*$ and some $\kappa$ as
above.  The Bianchi identities \eqrefs{Bianchi}{BianchiI} and
\eqrefs{Bianchi}{BianchiII} are equations in bundle-valued
three-forms. An arbitrary three-form $\lambda$ on $M$ vanishes if and
only if $\lambda (U_1, U_2, {\mathbf n})=0$ for all sections $U_1,U_ 2
\subset {\mathcal H}$. Applying this fact to the Bianchi identities
gives
\begin{gather*}
        \left.
        \begin{split}
          dA(\nabla_{U_1}{\mathbf n},U_2)+dA(U_1,\nabla_{U_2}{\mathbf n})&=0\\
          (\bar\nabla_{U_1}(\nabla {\mathbf n} ))U_2
          -(\bar\nabla_{U_2}(\nabla {\mathbf n} ))U_1&
          =2J(\omega(U_1)U_2-\omega(U_2)U_1)
        \end{split}
        \enspace\right\}\quad\text{(Bianchi I)}\\
      (\partial_{\mathbf n} \kappa) \, dA(U_ 1, U_ 2)=
      (\bar\nabla_{U_1} \omega )U_2- (\bar\nabla_{U_2}\omega )U _1
      \quad\text{(Bianchi II),}
\end{gather*}
for arbitrary $U_1,U_2 \subset {\mathcal H} $. From the first of these
equations one deduces \eqref{dooda}; from the second equation, that
$\omega=-\frac{1}{2}\sigma(\partial_-(\nabla {\mathbf n}))$; and from
the third, that $\partial_{\mathbf n} \kappa=\curl_{\mathcal
  H}(\sigma^{-1}(\omega ))$.
\end{proof}
\subsection{The maximally symmetric case}
Since $d \theta$ is already $\mathfrak g$-invariant,
\eqrefs{normtor}{mice4} implies that $\mathfrak g$ is torsion-reduced
if and only if $\nabla {\mathbf n}$ is $\mathfrak g$-invariant. But
$\mathfrak g$-invariance implies ${\mathfrak h}$-invariance, which, by
Proposition \ref{cr}, is the same thing as ${\mathbb
  C}$-linearity. However, $\nabla {\mathbf n}$ is ${\mathbb C}$-{\em
  anti}\/linear, by \eqref{dooda} above, and we conclude that
$\mathfrak g$ is torsion-reduced if and only if $\nabla {\mathbf n}=0
$.

Suppose that $\nabla {\mathbf n} = 0$. Then ${\mathbf n} $ is
automatically a symmetry of ${\mathfrak g} $, and hence an
infinitesimal isometry of the subriemannian contact structure. This is
a consequence of \eqrefs{generator}{py3}. Note that if the rank-one
foliation generated by ${\mathbf n} $ fibrates over some surface
$\Sigma$, then the invariant function $\kappa$ drops to a function on
$\Sigma$, by \eqref{sweet} above. In any case, Theorem \ref{goosey2}
applies, because of \eqrefs{normtor}{mice1}. Using the formula for
$\Omega =-\curvature \bar\nabla $ above, one applies this theorem and
its corollary to obtain:
\begin{proposition}[Compare with \cite{Hughen_95}]
  Suppose $\nabla {\mathbf n} = 0$. Then ${\mathfrak g} \subset
  J^1(TM)$ has an associated Cartan algebroid, namely $\mathfrak g$
  itself. If ${\mathcal U}\subset M$ is an arbitrary open subset, and
  $\mathfrak g_0$ the Lie algebra of all infinitesimal isometries of
  the subriemannian contact structure on ${\mathcal U}$, then
  $\dimension {\mathfrak g} _ 0\leq\rank {\mathfrak g} = 4$. If
  ${\mathcal U} $ is simply connected, then equality holds if and only
  if the function $\kappa $ defined by \eqref{lemonade} above is
  constant. In that case ${\mathfrak g}_ 0 \cong {\mathfrak b} \times
  {\mathbb R} $ (direct product) where ${\mathfrak b} $ is the Lie
  algebra of infinitesimal isometries (Killing fields) of the
  Euclidean plane, hyperbolic plane, or sphere, according to whether
  $\kappa= 0$, $\kappa < 0$, or $\kappa > 0$.
\end{proposition}
\subsection{Invariant differential operators}\lab{reps}
The normalized generator $\nabla $ determines associated invariant
differential operators, as explained in \ref{neoAssociated}. These
operators will be invariants of the subriemannian contact structure
analogous to the divergence, gradient, curl, etc., of a Riemannian
three-manifold.

Noting that the structure kernel $\mathfrak{h}\subset \mathfrak{g} $
of $\mathfrak{g}$ is an $\mathfrak{o}(2)$-bundle, we construct
irreducible representations of $\mathfrak{g}$ by mimicking a known
construction of the irreducible representations of the Lie algebra
$\mathfrak{o}(2)$. At least locally, this construction accounts for
all irreducible representations of $\mathfrak g$.

Define ${\mathcal H}_0:={\mathbb C} \times M$, ${\mathcal
  H}_1:={\mathcal H}$, and define ${\mathcal H}_2$ as in \ref{weight}
above. More generally, we define 
\begin{equation*}
  {\mathcal H}_k:=
  \symmetric^{k-1}_{\mathbb C}(\bar {\mathcal H})\otimes_{\mathbb C} {\mathcal H};\quad k\ge 1,
\end{equation*}
where $\bar {\mathcal H}$ is ${\mathcal H}$ with the complex structure
$-J$. Each ${\mathcal H}_k$ is simultaneously a $\mathfrak
g$-representation and a complex line-bundle, the two structures being
related according to 
\begin{equation*}
  \adjoint_Jq=kiq;\qquad q\subset {\mathcal H}_k,\,k\ge 0.
\end{equation*}
Every ${\mathcal H}_k$ is irreducible as a (real) $\mathfrak
g$-representation, except ${\mathcal H}_0$, which is two copies of the
irreducible trivial representation ${\mathbb R} \times M$.

Recall that for each section $q\subset E$ of an irreducible
${\mathfrak g}$-representation $E$, our objective is to derive the
decomposition of $\bar\nabla q\subset T^*\!M \otimes E$ corresponding
to the decomposition of $T^*\!M \otimes E$ into irreducibles. When
$E={\mathbb R} \times M$, $q$ is just a real function and, because
$TM={\mathcal H}\otimes \langle {\mathbf n} \rangle$, the
decompositions are straightforward:
\begin{equation*}
  T^*\!M \otimes ({\mathbb R} \times M)\cong {\mathcal H} \oplus ({\mathbb R} \times M),\quad
  \bar\nabla q=\gradient_{\mathcal H}q\oplus \partial_{\mathbf n} q,
\end{equation*}
where $\gradient_{\mathcal H}q$ and $\partial_{\mathbf n} q$ are as
defined previously. To analyze the case $E={\mathcal H}_k$ ($k\ge 1)$
we need the decomposition of ${\mathcal H}^*\otimes {\mathcal H}_k$
into irreducibles. To this end, note that for any $k\ge 1$, we have an
exact sequence,
\begin{equation}
  0\rightarrow {\mathcal H}_{k+1}\monomorphism {\mathcal H}^* \otimes {\mathcal H}_k
  \xrightarrow{\pi_-}{\mathcal H}_{k-1}\rightarrow 0,\mathlab{coughs}
\end{equation}
where $\pi_-$ is defined implicitly via
\begin{equation*}
  Q(U_1,U_2,V_1,\ldots,V_{k-2})-Q(U_2,U_1,V_1,\ldots,V_{k-2})=dA(U_1,U_2)(\pi_-Q)(V_1,\ldots,V_{k-2}),
\end{equation*}
for $k\ge 2$, and via
\begin{equation*}
  \langle Q(U_1),U_2\rangle-\langle Q(U_2),U_1\rangle=dA(U_1,U_2)(\pi_-Q),
\end{equation*}
for $k=1$. In the latter case we are using the Hermitian product on
$\mathcal H$ defined by $\langle U,V\rangle=\sigma(U,V)-i\,dA(U,V)$.

A compatible pair of splitting morphisms ${\mathcal
  H}_{k+1}\xleftarrow{\pi_+} {\mathcal H}^* \otimes {\mathcal
  H}_k\xleftarrow{s} {\mathcal H}_{k-1}$ are defined as follows:
\begin{equation*}
  (\pi_+Q)(U,V_1,\ldots,V_{k-1})=\frac{1}{2}\Big(\,Q(U,V_1,\ldots,V_{k-1})+
  iQ(JU,V_1,\ldots,V_{k-1})\,\Big),
\end{equation*}
for all $k\ge 1$, and 
\begin{equation*}
  (sq)(U_1,U_2,V_1,\ldots,V_{k-2})=\frac{i}{2}\Big\langle U_1,U_2\Big\rangle\, q(V_1,\ldots,V_{k-2}),
\end{equation*}
for $k\ge 2$, while 
\begin{equation*}
  (sq)U = \frac{1}{2}qU,
\end{equation*}
for $k=1$ ($q$ a ${\mathbb C} $-valued function).

Let $q$ be a section of ${\mathcal H}_k$. Then we have a restriction
$\bar\nabla q | {\mathcal H}\subset {\mathcal H}^* \otimes {\mathcal
  H}_k$. As \eqref{coughs} splits, we have ${\mathcal H}^* \otimes
{\mathcal H}_k\cong {\mathcal H}_{k-1} \oplus {\mathcal H}_{k+1} $ and
are led to define $\partial_+q:=\pi_+(\bar\nabla q|{\mathcal H} )$ and
$\partial_-q:=\pi_-(\bar\nabla q|{\mathcal H} )$. That is,
$\partial_+q\subset {\mathcal H}_{k+1}$ and $\partial_-q \subset
{\mathcal H}_{k-1}$ are defined by
\begin{equation*} (\partial_+q)(U,V_1,\ldots,V_{k-1})=\frac{1}{2}\Big(\,(\bar\nabla_Uq)(V_1,\ldots,V_{k-1})+
  i(\bar\nabla_{JU}q)(V_1,\ldots,V_{k-1})\,\Big),  
\end{equation*}
for any $k\ge 1$,
\begin{multline*}
    (\bar\nabla_{U_1}q)(U_2,V_1,\ldots,V_{k-2})-(\bar\nabla_{U_2}q)(U_1,V_1,\ldots,V_{k-2})\\=dA(U_1,U_2)(\partial_-q)(V_1,\ldots,V_{k-2}),
\end{multline*}
for $k\ge 2$, and 
\begin{equation*}
  \langle \bar\nabla_{U_1}q,U_2\rangle- \langle \bar\nabla_{U_2}q,U_1\rangle
  =dA(U_1,U_2)\,\partial_-q,
\end{equation*}
for $k=1$. This last formula simply means, for $q \subset {\mathcal H}$, that 
\begin{equation*}
  \partial_-q=\curl_{\mathcal H}(q)-i\divergence_{\mathcal H}(q).
\end{equation*}

Finally, for any section $q \subset {\mathcal H}_k$ and any $k\ge 1$,
we define $\partial_{\mathbf n} q:=\bar\nabla_{\mathbf n} q$, another
section of ${\mathcal H}_k$. 

Combining our observation ${\mathcal H}^* \otimes {\mathcal H}_k\cong
{\mathcal H}_{k-1} \oplus {\mathcal H}_{k+1} $ with the decomposition
$TM={\mathcal H}\otimes \langle {\mathbf n} \rangle$, we now obtain:
\begin{proposition}
  For any $k\ge 1$, we have a natural decomposition 
\begin{equation*}
  T^*\!M \otimes {\mathcal H}_k\cong {\mathcal H}_{k-1} \oplus {\mathcal H}_k \oplus {\mathcal H}_{k+1},
\end{equation*}
and a corresponding decomposition of operators 
\begin{equation*}
  \bar\nabla =\partial_- \oplus \partial_{\mathbf n} \oplus \partial_+,
\end{equation*}
where $\partial_-,\partial_{\mathbf n}$ and $\partial_+$ are defined above.
\end{proposition}

\section{Advanced prolongation theory}\lab{method} 
This section describes in detail the prolongation ${\mathfrak
  g}^{(1)}$ of a {\em surjective} infinitesimal geometric structure
$\mathfrak g\subset J^1{\mathfrak t}$. This description, obtained by
`prolonging' a generator of $\mathfrak g$ to a generator of $\mathfrak
g ^{(1)} $, is concrete enough to permit computations in examples. The
central result, Theorem \ref{jit}, is probably the most widely
applicable in such calculations. This theorem is used to prove several
theoretical results stated earlier; some special cases are also
considered.

We assume throughout that $\mathfrak t$ is a {\em transitive} Lie
algebroid.  For basic implications, see Sect.\,\ref{sectionY} under
`Assumption.'  We continue to denote the structure kernel of
${\mathfrak g}$ by ${\mathfrak h}$, and the associated lower
coboundary morphism, defined in \ref{structure:kernel}, by $\delta$.

\subsection{Natural connections}\lab{natural}%
Call a linear connection $D$ on ${\mathfrak g}\subset J^1{\mathfrak 
t}$ {\df natural} if the associated ${\mathfrak g}$-connection on 
${\mathfrak t}$ --- see Example \eqrefs{neoAssociated}{kiwi2} --- is the 
adjoint representation of ${\mathfrak g}\subset J^1{\mathfrak t}$ on 
${\mathfrak t}$; in symbols, if
\begin{equation*}
		aD_{\#V}X+[aX, V]_{\mathfrak t} 
		  =\adjoint_X^{\mathfrak t} V;\qquad V\subset{\mathfrak 
		  t},\,X\subset{\mathfrak g}.
\end{equation*}
Here $a\colon {\mathfrak g}\rightarrow {\mathfrak t}$ is the
restricted projection $J^1 {\mathfrak t} \rightarrow {\mathfrak
  t}$. The following proposition is not immediately useful in
computations but is a natural intermediate result. It stands between
the rather abstract Proposition \ref{qlook} and the computationally
useful Theorem \ref{jit} given later.
\begin{proposition}
	 Let $D$ be any natural connection on ${\mathfrak g}$ and  consider 
	 the morphism,
	 \begin{gather*}
		{\mathfrak g}\xrightarrow{\dot\Theta}\alternating^2 
	      (TM)\otimes{\mathfrak t}\\
        \dot\Theta (X)(U_1, U_2):=a(\curvature{D}\,(U_1, U_2)X),
	 \end{gather*}
	 where $ a$ is the projection ${\mathfrak g}\rightarrow{\mathfrak 
	 t}$. Then:
	 \begin{conditions}
         \item\lab{prop1} The morphism $\Theta\colon{\mathfrak
             g}\rightarrow h({\mathfrak g})$, defined in \ref{qlook},
           coincides with the composite
           \begin{equation*}
             {\mathfrak g}\xrightarrow{\dot\Theta}\alternating^2 
             (TM)\otimes{\mathfrak t}\xrightarrow{/\image\delta}h({\mathfrak 
               g}).
	   \end{equation*}
	 \end{conditions}
	 Moreover, if $\kernel\delta$ and $\kernel\Theta$ have constant 
	 rank (so that ${\mathfrak g}^{(1)}\subset J^1{\mathfrak g}$ is an 
	 infinitesimal geometric structure, by Proposition \ref{qlook}) then:
	 \begin{conditions}
         \item\lab{prop2} If $\Theta=0$, then all generators of
           ${{\mathfrak g}^{(1)}}$ are natural.
%            There exist natural connections on
%            ${\mathfrak g}$ generating ${\mathfrak g}^{(1)}$, with all
%            generators being natural if $\Theta=0$.
        \item\lab{gup2}A natural connection $D$ on ${\mathfrak g}$
           generates ${\mathfrak g}^{(1)}$ if and only if
           \begin{equation*}
             \curvature D(U_1,U_2) X\in{\mathfrak h}\qquad 
             \text{for all $X\in\kernel\Theta$ and $U_1,U_2\in TM$.}
           \end{equation*}
	 \end{conditions}
\end{proposition}
\begin{corollary}
  If $\mathfrak g \subset J^1{\mathfrak t}$ is a $\Theta $-reduced
  infinitesimal geometric structure, then a linear connection $D $ on
  ${\mathfrak t}$ generates $\mathfrak g$ if and only if $D$ is
  natural and $\curvature D \subset\alternating^2
  (TM)\otimes{\mathfrak g}^*\otimes{\mathfrak h}$.
\end{corollary}
\begin{corollary}
  The Cartan connection $\nabla^{(1)}$ in Theorem \ref{goosey} is the
  unique natural connection on ${\mathfrak g}$ such that
  $\curvature\nabla^{(1)} \subset\alternating^2 (TM)\otimes{\mathfrak
    g}^*\otimes{\mathfrak h}$. \end{corollary} The second corollary
holds because $D:=\nabla^{(1)}$ is the unique generator of ${\mathfrak
  g}^{(1)}$ under the hypotheses of Theorem \ref{goosey}. The
proposition's proof rests on parts \eqref{z1bz} and \eqref{m2} of the
following technical result (\eqref{zzz} is needed in the next
section). \begin{lemma}\mbox{}
  \begin{conditions}
  \item\lab{z1bz} A linear connection $D$ on ${\mathfrak g}$ is
    natural if and only if it is a generator of the isotropy $
    (J^1{\mathfrak g})_a\subset J^1{\mathfrak g}$ of $ a$.
		
  \item\lab{m2} If $D$ is a natural connection on ${\mathfrak g}$,
    then
    \begin{equation*}
      (sX\cdot da)^\vee(U_1,U_2)=a(\curvature D (U_1,U_2) X);\qquad 
      X\subset{\mathfrak g}.
    \end{equation*}
    Here $ s\colon{\mathfrak g}\rightarrow J^1{\mathfrak g}$ denotes
    the splitting of the exact sequence,
    \begin{equation*}
      0\rightarrow T^*\!M\otimes{\mathfrak g}\rightarrow J^1{\mathfrak 
        g}\rightarrow{\mathfrak g}\rightarrow 0.
    \end{equation*}
    corresponding to $D$, and $ (sX\cdot da)^\vee$ is the `reduction'
    of $ sX\cdot da$ as per \eqrefs{qlook}{m1}.
		
  \item\lab{zzz} If $\nabla$ is a connection on ${\mathfrak t}$
    generating $\mathfrak g$ and $\nabla^{\mathfrak h}$ is any linear
    connection on $\mathfrak h$, then, identifying $\mathfrak g$ with
    $\mathfrak t\oplus\mathfrak h$ using the generator $\nabla$, every
    natural connection on $\mathfrak g$ is of the form
    \begin{equation*}
      D_U(V\oplus\phi)=(\nabla_UV+\phi(U))\oplus(\nabla^{\mathfrak h}_U\phi
      +\epsilon(V\oplus\phi)U),
    \end{equation*}
    for some vector bundle morphism $\epsilon\colon\mathfrak
    t\oplus\mathfrak h\rightarrow T^*\!M\otimes\mathfrak h$.
  \end{conditions}
\end{lemma}
\begin{proof}
Let $\mathcal D$ be the deviation operator described in \ref{boot}. Then 
\begin{equation}
  [aX,V]_{\mathfrak t} = \adjoint^{\mathfrak t}_{J^1 (aX)}V
 =\adjoint^{\mathfrak t}_XV-{\mathcal D}_{\#V}X;
  \qquad X \subset {\mathfrak g}, V \subset {\mathfrak t}.\mathlab{buckle}
\end{equation}
Taking care not to confuse $D$'s with ${\mathcal D}$'s, we also have
\begin{equation*}
  aD_{\#V}X=a(sX-J^1 X)(\#V)=a({\mathcal D}(sX))\#V=a {\mathcal D}_{\#V}(sX).
\end{equation*}
Here $s$ is the splitting in \eqref{m2}. Combining this with \eqref{buckle} gives
\begin{equation*}
  a D_{\#V}X+[aX,V]_{\mathfrak t}-\adjoint_X^{\mathfrak t}V
 =a {\mathcal D}_{\#V}(sX)-{\mathcal D}_{\#V}(\,p(sX)\,).
\end{equation*}
The claim in \eqref{z1bz} now follows from
\eqrefs{tautological}{flute1} and transitivity. Similarly, \eqref{m2}
may be derived as a consequence of \eqrefs{tautological}{flute2} and
transitivity. One checks that the connection in \eqref{zzz} is natural
and, with the help of \eqref{z1bz} and \eqrefs{qlook}{z1z}, that all
possibilities are covered, establishing \eqref{zzz}.
\end{proof}
\begin{proof}[Proof of proposition]
By \eqref{z1bz}, $D$ generates $ (J^1{\mathfrak g})_a$.  Recalling 
that $ (J^1{\mathfrak g})_a$ is surjective (see \eqrefs{qlook}{z1z}), 
let $s\colon{\mathfrak g}\rightarrow (J^1{\mathfrak g})_a$ denote the 
splitting of the exact sequence,
\begin{equation}
         0\rightarrow T^*\!M\otimes{\mathfrak h}\monomorphism  (J^1{\mathfrak 
         g})_a\rightarrow{\mathfrak g}\rightarrow 0,\mathlab{when7}
\end{equation}
determined by the generator $D$.  By the commutativity of the diagram 
\eqrefs{qlook}{diag} defining $\Theta$, we have
\begin{equation*}
        \Theta (X)=\theta (sX)\modulo\image\delta 
		=(sX\cdot da)^\vee \modulo\image\delta.
\end{equation*}
Invoking \eqref{m2}, we prove \eqref{prop1}.  

If $\Theta=0$, then ${{\mathfrak g}^{(1)}} $ is surjective, and so
$s({\mathfrak g})\subset {{\mathfrak g}^{(1)}} $ for any generator
$D$ of ${{\mathfrak g}^{(1)}} $; here $s \colon {\mathfrak g} \rightarrow {{\mathfrak g}^{(1)}}$ is the corresponding splitting of the exact sequence in \eqref{m2}. This being the case we have, in particular, $s({\mathfrak
  g})\subset (J^1{\mathfrak g})_a$, making $D$ a generator of $(J^1{\mathfrak g})_a$ also. By \eqref{z1bz}, $D$ is natural. This proves \eqref{prop2}. 

By \eqref{z1bz}, a natural connection $D$ generates ${\mathfrak
  g}^{(1)}$ if and only if $(sX\cdot da)^\vee=0$ for all $X\subset\kernel \Theta $, where $s$ is the splitting morphism corresponding to $D$. Conclusion \eqref{gup2} is therefore a consequence of \eqref{m2}.  \end{proof}

\subsection{Prolonging a generator}\lab{jit}
We now show how to construct a generator for the prolongation
${\mathfrak g}^{(1)}$ of ${\mathfrak g}$ given a generator of
${\mathfrak g}$. Fix a generator $\nabla$ of ${\mathfrak g}$; recall
that this is a certain linear connection on ${\mathfrak t}$.  Also
fix an arbitrary linear connection $\nabla^{\mathfrak h}$ on
${\mathfrak h}$.  (In practice, there is usually a preferred choice,
and if $\mathfrak t=TM$ this is always true; see \ref{caviar}.  No
canonical choice exists in general, however).

With $\nabla$ and $\nabla^{\mathfrak h}$ fixed, there is vector bundle 
morphism
\begin{equation*}
		{\mathfrak t}\oplus{\mathfrak 
		h}\xrightarrow{\tilde\Theta}\alternating^2 
		(TM)\otimes{\mathfrak t}
\end{equation*}
well defined by its action on sections via
\begin{equation}
   \tilde\Theta\, (V\oplus\phi):=\curvature\nabla 
   (\,\cdot\,,\,\cdot\,)V+d_\nabla\phi-\delta (\nabla^{\mathfrak 
   h}\phi);\qquad V \subset {\mathfrak t},\,\phi \subset {\mathfrak h}.\mathlab{box1}
\end{equation}
Here
\begin{align*}
    d_\nabla\phi\, (U_1, U_2)&:=\nabla_{U_1}(\phi (U_2))
                                                  -\nabla_{U_2}(\phi (U_1))
                                                  -\phi ([U_1,U_2]),\\
  \text{and}\qquad       (\nabla^{\mathfrak h}\phi)U&:=\nabla_U^{\mathfrak h}\phi.
\end{align*}

\begin{theorem}[Prolonging a generator of ${\mathfrak g}\subset 
J^1{\mathfrak t}$]
Let $\mathfrak g$ be a surjective infinitesimal geometric structure on
a transitive Lie algebroid $\mathfrak t$, with $\mathfrak h$,
$\nabla$, $\nabla^\mathfrak h$, and $\tilde\Theta$ as above.  Use
$\nabla$ to identify ${\mathfrak g}$ with ${\mathfrak
  t}\oplus{\mathfrak h}$.  Then:
\begin{conditions}
\item\lab{ty1} The composite morphism,
  \begin{equation*}
    {\mathfrak g}\cong{\mathfrak t}\oplus{\mathfrak 
      h}\xrightarrow{\tilde\Theta}\alternating^2 
    (TM)\otimes{\mathfrak t} \xrightarrow{/\image\delta} 
    h({\mathfrak g}),
  \end{equation*}
  coincides with the morphism $\Theta\colon{\mathfrak g}\rightarrow
  h({\mathfrak g})$ defined in \ref{qlook}.
                
\item\lab{ty1c} $\kernel\Theta\subset{\mathfrak t}\oplus{\mathfrak h}$
  is precisely the set of all $ V\oplus\phi$ for which the {\df
    generator equation},
  \begin{equation*}
    \delta (\epsilon)=\tilde\Theta (V\oplus\phi),%\mathlab{sga}
  \end{equation*}
  admits a solution $\epsilon\in T^*\!M\otimes{\mathfrak h}$.
                          
\item\lab{ty3} Assuming $\kernel\delta$ and $\kernel\Theta$ have
  constant rank (so that ${\mathfrak g}^{(1)}\subset J^1{\mathfrak g}$
  is an infinitesimal geometric structure, by Proposition \ref{qlook})
  a linear connection $\nabla^{(1)}$ on ${\mathfrak g}\cong{\mathfrak
    t}\oplus{\mathfrak h}$ generating ${\mathfrak g}^{(1)}$ is given
  by
  \begin{equation*}
    \boxed{\nabla_U^{(1)} (V\oplus\phi)=(\nabla_UV+\phi 
      (U))\oplus (\nabla^{\mathfrak 
        h}_U\phi+\epsilon(V\oplus\phi)U),}
  \end{equation*}
  where $\epsilon\colon{\mathfrak t}\oplus{\mathfrak h}\rightarrow
  T^*\!M\otimes{\mathfrak h}$ is any of the vector bundle morphisms
  for which $\epsilon:=\epsilon (V\oplus\phi)$ solves the generator
  equation defined in \eqref{ty1c} above, for each $V\oplus\phi$ in
  $\kernel\Theta$.  If ${\mathfrak g}^{(1)}$ is surjective (i.e.,
  $\Theta =0$) then every generator is of the above form.
\end{conditions} \end{theorem} \begin{proof} Let $ D$ denote the
general form of a natural connection on ${\mathfrak g}$ given in
\eqrefs{natural}{zzz}, with $\epsilon\colon{\mathfrak
  t}\oplus{\mathfrak h}\rightarrow T^*\!M\otimes{\mathfrak h}$
completely arbitrary. If $\dot\Theta$ is the morphism defined in
Proposition \ref{natural}, then one computes \begin{equation*}
  \dot\Theta (V\oplus\phi) (U_1, U_2)=
  a\Big(\,\curvature D\, (U_1, U_2)(V\oplus\phi)\,\Big)=\tilde\Theta 
  (V\oplus\phi)-\delta (\epsilon (V\oplus\phi)),
\end{equation*}
where $\tilde\Theta$ is the morphism defined by \eqref{box1}.  
Conclusion \eqref{ty1} of the theorem now follows from Proposition 
\eqrefs{natural}{prop1}.  Conclusion \eqref{ty1c} is just a 
consequence of \eqref{ty1}.  One obtains \eqref{ty3} by taking 
$\nabla^{(1)}:=D$; choosing $\epsilon$ as described guarantees the 
curvature condition in Proposition \eqrefs{natural}{gup2}.  If 
$\Theta=0$ then every generator of ${\mathfrak g}^{(1)}$ is of the 
stated form because every generator is natural (Proposition 
\eqrefs{natural}{prop2}), and all natural connections are in included in the general form \eqrefs{natural}{zzz}.
\end{proof}

\subsection{Torsion revisited}\lab{revisited}
The preceding theorem   will allow us to prove Theorem \ref{qtax} 
relating $\Theta$-reduction and torsion reduction. The appropriate 
argument also leads to an alternative form of Theorem \ref{jit} when 
${\mathfrak t}=TM$. See \ref{caviar} below.

Although the definition of $\tilde\Theta$ in \eqref{box1} above is 
explicit, it depends on a choice of linear connection 
$\nabla^{\mathfrak h}$ on ${\mathfrak h}$.  Here is an implicit 
formula depending only on a choice of generator $\nabla$ for 
${\mathfrak g}$:
\begin{proposition}
  Let $\bar\nabla$ denote the ${\mathfrak t}$-connection on
  ${\mathfrak t}$ associated with a generator $\nabla$ of ${\mathfrak
    g}\subset J^1{\mathfrak t}$. Then for arbitrary sections $
  X\subset{\mathfrak g}$ and $ U_1, U_2\subset{\mathfrak t}$, one has
  \begin{equation*}
    \tilde\Theta(X)(\#U_1,\#U_2)=\Big(\,X\cdot\torsion\bar\nabla+\Delta 
    (\cocurvature\nabla (aX,\,\cdot\,))\,\Big)(U_1,U_2),
  \end{equation*}
  where $\Delta$ is the upper coboundary morphism and $
  a\colon{\mathfrak g}\rightarrow{\mathfrak t}$ the projection.
\end{proposition}
\noindent Note that $\cocurvature\nabla (aX,\,\cdot\,)$ is a section 
of $ T^*\!M\otimes{\mathfrak h}$, by Proposition 
\eqrefs{reconstruction}{ala}.
\begin{proof}
Since ${\mathfrak t}$ is assumed to be transitive, there exists a 
linear connection $\nabla^{\mathfrak h}$ on ${\mathfrak h}$ such that 
$\nabla_{\#U}^{\mathfrak h}=\bar\nabla_U$ for  all $ U\in{\mathfrak 
t}$. Here $\bar\nabla$ denotes the associated ${\mathfrak 
t}$-connection on ${\mathfrak h}$ discussed in 
\eqrefs{neoAssociated}{kiwi}. After a little manipulation, we obtain
\begin{equation}
  (d_\nabla\phi-\delta (\nabla^{\mathfrak h}\phi)) (\#U_1,\#U_2)=
  (\phi\cdot\torsion\bar\nabla) (U_1,U_2);\quad 
  U_1,U_2\subset{\mathfrak t}.\mathlab{bird1}
\end{equation}
Note that $ T^*\!M\otimes{\mathfrak t}$ (of which $\phi$ is a section) 
acts on $\alternating^2 ({\mathfrak t})\otimes{\mathfrak t}$ as a 
subalgebroid of $J^1{\mathfrak t}$ (which acts on ${\mathfrak t}$ via 
adjoint action).

Replace ${\mathfrak g}$ in Proposition \ref{dual} with ${\mathfrak t}$ 
and replace $\nabla$ there by the composite
\begin{equation*}
	{\mathfrak t}\xrightarrow{\#}TM\xrightarrow{\nabla}\gl({\mathfrak t}).
\end{equation*}
Then part  \eqref{dual2} of that proposition delivers the formula,
\begin{align*}
	\curvature\nabla\,(\#U_1,\#U_2)V&=(\bar\nabla_V\torsion\bar\nabla)(U_1,U_2) 
	-\curvature\bar\nabla\,(V, U_1)U_2\\
                       &+\curvature\bar\nabla\,(V, U_2)U_1;\qquad 
                       U_1,U_2,V\subset{\mathfrak t}.
\end{align*}
Applying Proposition \eqrefs{cocurvature}{determine}, we may rewrite 
this as,
\begin{equation}
	\curvature\nabla 
	(\#U_1,\#U_2)V=(\bar\nabla_V\torsion\bar\nabla+\Delta (
	\cocurvature\nabla (V,\,\cdot\,))\,)(U_1,U_2).\mathlab{bird2}
\end{equation}
Substituting \eqref{bird1} and \eqref{bird2} into the definition 
\eqrefs{jit}{box1} of $\tilde\Theta$ gives,
\begin{equation*}
	\tilde\Theta (V\oplus\phi) (\#U_1,\#U_2)=\\
    \Big(\,\bar\nabla_V\torsion\bar\nabla+\phi\cdot\torsion\bar\nabla+\Delta
	  (\cocurvature\nabla (V,\,\cdot\,))\,\Big)(U_1,U_2).
\end{equation*}
Under the identification ${\mathfrak g}\cong{\mathfrak 
t}\oplus{\mathfrak h}$ determined by the generator $\nabla$, we 
obtain the stated formula.
\end{proof}
\begin{proof}[Proof of Theorem \ref{qtax}]
By Theorem \eqrefs{jit}{ty1} and the commutativity of 
\eqrefs{structure:kernel}{dg7}, we have,
\begin{equation*}
	 	\psi (\Theta (X))=i(\tilde\Theta (X))\modulo \image\Delta,
\end{equation*}
where $ i\colon\alternating^2 (TM)\otimes{\mathfrak 
t}\rightarrow\alternating^2 ({\mathfrak t})\otimes{\mathfrak t}$ 
denotes the natural inclusion ($\#^*\colon T^*\!M \rightarrow {\mathfrak t}^* $ being injective).  The proposition above then gives
\begin{equation*}
	 	\psi (\Theta 
	 	(X))=X\cdot\torsion\bar\nabla\modulo\image\Delta=X\cdot\tau.
\end{equation*}
\end{proof}

\subsection{The special case ${\mathfrak g}\subset J^1 
(TM)$}\lab{caviar} We now specialize Theorem \ref{jit} to the case 
${\mathfrak t}=TM$. As an application, we complete the proof of 
Theorem \ref{goosey2}, the last unproven assertion of preceding 
sections.

Let ${\mathfrak g}\subset J^1 (TM)$ be an infinitesimal geometric 
structure with structure kernel ${\mathfrak h}$ and generator $\nabla$. 
Let $\bar\nabla$ denote the dual connection, i.e., 
$\bar\nabla_UV=\nabla_VU+[U, V]$ and define
\begin{gather*}
	 TM\oplus{\mathfrak 
	   h}\xrightarrow{\tilde{\tilde\Theta}}\alternating^2 (TM)\otimes TM\\
     \tilde{\tilde\Theta}(V\oplus\phi)
       :=\bar\nabla_V\torsion\bar\nabla+\phi\cdot\torsion\bar\nabla
	   =(V\oplus\phi)\cdot\torsion\bar\nabla.
\end{gather*}
\begin{theorem}[Prolonging a generator of ${\mathfrak g}\subset J^1(TM)$]
  With ${\mathfrak g}$, ${\mathfrak h}$, $\nabla$, and
  $\tilde{\tilde\Theta}$ as above, we have:
  \begin{conditions}
  \item The composite morphism,
    \begin{equation*}
      {\mathfrak g}\cong TM\oplus{\mathfrak h}\xrightarrow{\tilde{\tilde\Theta}}
      \alternating^2 (TM)\otimes TM\xrightarrow{\image\Delta}
      H({\mathfrak g})
    \end{equation*}
    coincides with the morphism $\Theta\colon{\mathfrak g}\rightarrow
    h({\mathfrak g})=H({\mathfrak g})$ defined in \ref{qlook}.
	
  \item\lab{sga2} $\kernel\Theta\subset{\mathfrak g}\cong
    TM\oplus{\mathfrak h}$ is precisely the set of all $ V\oplus\phi$
    for which the {\df generator equation},
    \begin{equation*}
      \Delta (\epsilon)=\tilde{\tilde\Theta} (V\oplus\phi),
    \end{equation*}
    admits a solution $\epsilon\in T^*\!M\otimes{\mathfrak h}$.
	
  \item\lab{three} Assume $\kernel\Delta$ and $\kernel\Theta$ have
    constant rank, so that ${\mathfrak g}^{(1)}\subset J^1{\mathfrak
      g}$ is an infinitesimal geometric structure (by Proposition
    \ref{qlook}) and that $H({\mathfrak g})$ has constant rank (is a
    ${\mathfrak g}$-representation). Then $\Theta (X)=X\cdot\tau$,
    where $\tau\subset H({\mathfrak g})$ is the intrinsic
    torsion. Also, a linear connection $\nabla^{(1)}$ on ${\mathfrak
      g}\cong TM\oplus{\mathfrak h}$ generating ${\mathfrak g}^{(1)}$
    is given by
    \begin{equation*}
      \boxed{\nabla^{(1)}_U(V\oplus\phi)=(\nabla_UV+\phi (U))\oplus
        (\bar\nabla_U\phi+\epsilon (V\oplus\phi) U+\curvature\bar\nabla (U, 
        V)),}
    \end{equation*}
    where $\epsilon\colon TM\oplus{\mathfrak h}\rightarrow
    T^*\!M\otimes{\mathfrak h}$ is any of the vector bundle morphisms
    for which $\epsilon:=\epsilon (V\oplus\phi)$ solves the generator
    equation defined in \eqref{sga2} above, for each $V\oplus\phi$ in
    $\kernel\Theta$.  If ${\mathfrak g}^{(1)}$ is surjective (i.e.,
    $\Theta =0$) then every generator is of the above form.
  \end{conditions}
\end{theorem}
\begin{note}
  The $\epsilon$'s solving the generator equation above, and the
  generator equation of Theorem \ref{jit}, are different.
\end{note}
\begin{proof}
  In \ref{jit} above take ${\mathfrak t}=TM$ and let
  $\nabla^{\mathfrak h}$ be the $ TM$-connection on ${\mathfrak h}$
  associated with the generator $\nabla$ (given by
  \eqrefs{reconstruction}{storke} with ${\mathfrak t}={\mathfrak t}_1=
  TM $). Then Proposition \ref{revisited} gives
\begin{equation*}
  \tilde\Theta (V\oplus\phi)=\tilde{\tilde\Theta}(V\oplus\phi)+\Delta 
  (\cocurvature\nabla (V,\,\cdot\,).
\end{equation*}
Noting that $\cocurvature\nabla=-\curvature\bar\nabla$ and 
$\delta=\Delta$ (because ${\mathfrak t}=TM$) one obtains the  stated 
results as a special case of Theorem \ref{jit} and Theorem \ref{qtax}.
\end{proof}
\begin{proof}[Proof of Theorem \ref{goosey2}]
The hypothesis that ${\mathfrak g}$ be torsion-reduced means 
$\Theta=0$.  So the generator equation defined in \eqref{sga2} above has a 
solution for all $ V\oplus\phi\in TM\oplus{\mathfrak h}$. The solution 
is unique because $\Delta$ is injective, by hypothesis. Part \eqref{qd1} 
of  \ref{goosey2} follows. The Cartan connection on ${\mathfrak g}$ 
in Theorem \ref{goosey2} is the unique generator of  ${\mathfrak 
g}^{(1)}$; conclusion \eqref{three} above implies that it has the form 
given in part  \eqref{scar} of  \ref{goosey2}.

Suppose ${\mathfrak g}$ is reductive and let $\nabla$ be a normal 
generator. Then $\torsion\bar\nabla\subset C$ for some ${\mathfrak 
g}$-invariant complement $ C\subset\alternating^2 (TM)\otimes TM$ of 
the image of $\Delta$. In particular,
\begin{equation*}
	\bar\nabla_V\torsion\bar\nabla+\phi\cdot\torsion\bar\nabla=
	(V\oplus\phi)\cdot\torsion\bar\nabla
\end{equation*}
is a section of $ C$. However, this last also lies in $\image \Delta$ 
because its image under the projection
\begin{equation*}
	\alternating^2 (TM)\otimes TM\xrightarrow{\image\Delta}H({\mathfrak 
	g})
\end{equation*}
is $ (V\oplus\phi)\cdot\tau$, which vanishes because ${\mathfrak g}$ 
is torsion-reduced.  We conclude that 
$\bar\nabla_V\torsion\bar\nabla+\phi\cdot\torsion\bar\nabla=0$.  Part 
\eqref{qd3} of \ref{goosey2} follows.
\end{proof}

\subsection{Symmetries of torsion-free affine
  structures}\lab{illustration2}
We now offer a simple application of Theorem \ref{caviar}.  Let
$\nabla$ be a torsion-free linear connection on $ TM$.  Then $\nabla$
is a generator of $ J^1 (TM)$.  Taking ${\mathfrak g}=J^1 (TM)$, the
generator equation in \eqrefs{caviar}{sga2} reads
$\Delta(\epsilon)=0$, with trivial solution $\epsilon=0$.  Theorem
\eqrefs{caviar}{three} delivers the following generator $\nabla^{(1)}$
for $ J^2 (TM)=(J^1 (TM))^{(1)}={\mathfrak g}^{(1)}$:
\begin{equation*}
        \nabla^{(1)} (V\oplus\phi)=(\nabla_UV+\phi (U))\oplus 
        (\nabla_U\phi+\curvature\nabla(U, V)).
\end{equation*}
Here we are identifying ${\mathfrak g}=J^1(TM)$ with $ TM\oplus ( 
T^*\!M\otimes TM)$ using $\nabla$.

Now let ${\mathfrak g}\subset J^2 (TM)$ instead denote the isotropy 
subalgebroid of $\nabla$, as described in \ref{linear2}.  Then it is 
not too difficult to check that $\nabla^{(1)}$ is even a generator for 
${\mathfrak g}$.  In this regard, a helpful formula, readily derived, 
is
\begin{equation*}
        \nabla^{(1)} (J^1 V)(U_1, U_2)=0\oplus ((J^2 V)\cdot\nabla)^\vee(U_1, 
        U_2),
\end{equation*}
for any section $ V\subset TM$.

The generator $\nabla^{(1)}$ is necessarily the Cartan connection 
on $ J^1 (TM)$ referred to in Proposition \ref{linear2}. Its 
curvature is given by
\begin{equation*}
  \curvature\nabla^{(1)} (U_1,U_2)(V\oplus\phi)
  =0\oplus\Big(\,-(\nabla_V\curvature\nabla+\phi\cdot\curvature\nabla) 
  (U_1,U_2)\,\Big).
\end{equation*}
In particular, $\curvature\nabla^{(1)}$ vanishes if and only if 
$\curvature\nabla$ is $\nabla$-parallel and $ T^*\!M\otimes 
TM$-invariant.  But as $\identity_{TM}$ is a section of $T^*\!M\otimes 
TM$, this happens if and only if $\curvature\nabla=0$.  In that case 
we obtain
\begin{equation*}
        \torsion\overline{\nabla^{(1)}} (V_1\oplus\phi_1, 
        V_2\oplus\phi_2)=(\,\phi_1(V_2)-\phi_2(V_1)\,)\oplus 
        [\phi_1,\phi_2]_{T^*\!M\otimes TM},      
\end{equation*}
where $\overline{\nabla^{(1)}}$ denotes the representation of 
$J^1(TM)$ on itself associated with the Cartan connection 
$\nabla^{(1)}$ on $ J^1(TM)$.  Applying Theorem 
\ref{dissymmetry}, we recover the following classical result:
\begin{proposition}
  Let $\nabla$ be a torsion-free linear connection on $TM$ and
  ${\mathfrak g}_0$ the Lie algebra of infinitesimal isometries of
  $\nabla$ on some open set \,$\mathcal U\subset M$.  Then
  $\dimension{\mathfrak g}_0\le\rank J^1 (TM)=n(n+1)$, $ n=\dimension
  M$.  If \,$\mathcal U$ is simply-connected then equality holds if
  and only if $\curvature\nabla=0$, in which case ${\mathfrak g}_0$ is
  naturally isomorphic to the semidirect product $
  T_mM\oplus(T^*_mM\otimes T_mM)$, $ m\in U$.
\end{proposition}

\section{Application: conformal structures}\lab{point} 
In this section we turn to the application of Cartan's method to
conformal structures.  Our results are summarized in Theorems
\ref{tha}, \ref{thb} and \ref{thc} below.
\subsection{The Lie algebroid setting} Let $\sigma $ be a Riemannian
metric on a smooth connected manifold $M$, with $n:=\dimension M\ge
3$.  Let $\langle\sigma\rangle$ denote its conformal class, viewed as
the one-dimensional subbundle of $\symmetric^2 (TM)$ generated by the
section $\sigma$.
  
Let ${\mathfrak g}_\sigma\subset J^1 (TM)$ denote the isotropy of 
$\sigma\subset\symmetric^2 (TM)$ (denoted $\mathfrak{g}$ in 
\ref{Riemannian}).  Let ${\mathfrak g}\subset J^1 (TM)$ denote the 
isotropy of $\langle \sigma\rangle\subset\symmetric^2 (TM)$.  
According to \ref{isotropy}, this means that the 1-jet of a vector 
field $ V$ at $ m\in M$ lies in ${\mathfrak g}$ if and only if $ 
({\mathcal L}_V\sigma) (m)\in\langle\sigma\rangle (m)$, were 
${\mathcal L}$ denotes Lie derivative.  Evidently, the symmetries of 
${\mathfrak g}$ are the infinitesimal isometries of the conformal 
structure $\langle\sigma\rangle$, henceforth called the {\df conformal 
Killing fields}.  

The rank of ${\mathfrak g}$ is constant (see below), making it an 
infinitesimal geometric structure on $ TM$.  Observe that ${\mathfrak 
g}$ is surjective because ${\mathfrak g}_\sigma$ is surjective and 
${\mathfrak g}_\sigma\subset{\mathfrak g}$.  If ${\mathfrak h}_\sigma$ and 
${\mathfrak h}$ denote the structure kernels of ${\mathfrak g}_\sigma$ 
and ${\mathfrak g}$ respectively, then
\begin{equation*}
  {\mathfrak h}={\mathfrak
    h}_\sigma\oplus\langle\identity_{TM}\rangle,\enspace
  \text{implying}\enspace\mathfrak g=\mathfrak
  g_\sigma\oplus\langle\identity_{TM}\rangle.
\end{equation*}
Here $\langle\identity_{TM}\rangle$ denotes the one-dimensional
subbundle of $T^*\!M\otimes TM$ generated by the identity section
$\identity_{TM}$.  If vector bundles $E_1$ and $E_2$ are
representations of $\mathfrak g$, then a vector bundle morphism
$\phi\colon E_1\rightarrow E_2$ is a morphism of $\mathfrak
g$-representations if and only if it is a morphism of $\mathfrak
g_\sigma$-representations commuting with the action of
$\identity_{TM}\subset\mathfrak g$.  In particular, this applies to
$\alpha\otimes V\mapsto (\alpha\otimes V)^\transpose:=\sigma
(V)\otimes\sigma^{-1}(\alpha)$ (the transpose involution of
$T^*\!M\otimes TM$) and to the epimorphism, $\sskew\colon
T^*\!M\otimes TM\rightarrow\mathfrak h_\sigma$, defined by $\sskew
(\phi):=(\phi-\phi^\transpose)$.  These morphisms and ${\mathfrak
  h}_\sigma$ depend only on the conformal class of $\sigma$.

Because the Levi-Cevita connection $\nabla$ associated with $\sigma$ 
generates  ${\mathfrak g}_\sigma$ 
(see \ref{illustration}), it also generates ${\mathfrak 
g}\supset{\mathfrak g}_\sigma$.

\subsection{Classical ingredients}\lab{outline}
From well-known representation-theoretic arguments we know that the
curvature of the Levi-Cevita connection $\nabla$ takes values in a
proper subbundle $E_{\mathrm{Weyl}}\oplus
E_{\mathrm{Ricci}}\subset\alternating^2 (TM)\otimes{\mathfrak
  h}_\sigma$, with $ E_\mathrm{Ricci}$ being (for $\dimension M\ge 3$)
the isomorphic image of $\symmetric^2 (TM)$ under the monomorphism of
$\mathfrak g$-representations
\begin{align*}
        T^*\!M\otimes T^*\!M&\xrightarrow{\coRicci}\alternating^2 
        (TM)\otimes{\mathfrak h}_\sigma\\
        \coRicci(\Phi) (V_1, V_2)&:=\sskew(\Phi V_1\otimes V_2-\Phi 
        V_2\otimes V_1).
\end{align*}
$E_{\mathrm{Weyl}} \subset \alternating^2 (TM)\otimes {\mathfrak
  h}_\sigma $ is the intersection of the kernels of the so-called
Bianchi and Ricci morphisms; see, e.g.,
\cite[p.\,230]{Sharpe_97}. Whence,
\begin{equation}
        \curvature\nabla=W+\coRicci (R)\mathlab{fortify}
\end{equation}
for uniquely determined sections $W\subset E_\mathrm{Weyl}$ and $
R\subset\symmetric^2 (TM)$.  These are the {\df Weyl} and {\df
  modified Ricci} curvatures of $\sigma$.  Both $ E_\mathrm{Weyl}$ and
$E_\mathrm{Ricci}$ are ${\mathfrak g}$-subrepresentations of
$\alternating^2(TM)\otimes {\mathfrak{h}_{\sigma}}$ and in particular
we may speak of the isotropy $\mathfrak{g}_W\subset\mathfrak{g}$ of $
W$. If $n=3$, $E_\mathrm{Weyl}=0$.

Also of significance will be the {\df Cotton-York tensor}. This is the
associated exterior derivative $d_\nabla R$ of 
$R\subset\symmetric^2 (TM)\subset T^*\!M\otimes T^*\!M$,  viewed as a 
$T^*\!M$-value one-form on $M$:
\begin{equation*}
                 d_\nabla R(U_1, 
                 U_2):=\nabla_{U_1}(R(U_2))-\nabla_{U_2}(R(U_1))-R([U_1,U_2]).
\end{equation*}
Alternatively, by torsion-freeness, $d_\nabla R$ is the image of 
$\nabla R$ under the composite
\begin{align*}
  T^*\!M\otimes\symmetric^2(TM)\monomorphism T^*\!M\otimes T^*\!M\otimes 
  T^*\!M&\rightarrow\alternating^2(TM)\otimes T^*\!M\\
\alpha\otimes\beta\otimes\gamma&\mapsto\alpha\wedge\beta\otimes\gamma.
\end{align*}

Bianchi's second identity \eqrefs{Bianchi}{BianchiII} for the
generator $\nabla$ enforces a relationship between the Cotton-York
tensor $d_\nabla R$, and the derivative $\nabla W$.  In particular, it
is well known that $W=0$ implies the vanishing of $d_\nabla R$ in all
dimensions except three, where $E_\mathrm{Weyl}=0$ and the values of
$d_\nabla R$ are restricted to a certain six-dimensional subbundle of
$\alternating^2 (TM)\otimes T^*\!M$.

\vspace{\baselineskip} The next two sections summarize the main
conclusions of our application of Cartan's method to conformal
structures.

\subsection{Conformal invariance of $W$}\lab{tha}
\begin{theorem}
        The Weyl curvature $W$ is an invariant of $\mathfrak g$ and therefore 
        a conformal invariant.
\end{theorem}
\noindent Classically, the conformal invariance of $W$ required an 
explicit check.  In our proof of this theorem conformal invariance is
simply a manifestation of the invariance of the morphism $\Theta$
defined in general in \ref{qlook}, as applied to the prolongation of
$\mathfrak{g}$ (rather than to $\mathfrak{g}$ itself). This is
invariant by construction.

\subsection{The $W=0$ case}\lab{thb}
Our second theorem likewise expresses, in invariant Lie algebroid
language, results that are essentially classical:
\begin{theorem}
  Suppose $W=0$. Then $\mathfrak{g}$ has an associated Cartan
  algebroid, namely its prolongation $\mathfrak{g}^{(1)}\subset
  J^2(TM)$, which is surjective and transitive.  Denoting the Cartan
  connection on $\mathfrak{g}^{(1)}$ by $\nabla^{(2)}$, we have:
  \begin{conditions}
  \item The $ \nabla^{(2)}$-parallel sections of
    $\mathfrak{g}^{(1)}\subset J^2(TM)$ are precisely the
    twice-prolonged conformal Killing fields.
    
  \item Each metric $\sigma$ in the conformal class determines natural
    isomorphisms $$\mathfrak{g}^{(1)}\cong\mathfrak{g}\oplus T^*\!M,
    \qquad\mathfrak{g} \cong TM \oplus \mathfrak{h},$$
    and an
    associated explicit formula for $\nabla^{(2)}$ (see
    \eqrefs{dutiful}{honda} and \eqrefs{third}{home2} below).
    
  \item If $ n\geq 4$, then $\nabla^{(2)}$ is automatically flat. If $
    n=3$, then $\nabla^{(2)}$ is flat if and only if  $d_\nabla R=0$. In
    particular, the Lie algebra $\mathfrak{g}_0$ of all conformal
    Killing fields over any simply-connected open set
    $\,\mathcal{U}\subset M$ satisfies
    \begin{equation*}
      \dimension \mathfrak{g}_0\leq\rank \mathfrak{g}^{(1)}
      =\frac{1}{2}(n+1)(n+2),
    \end{equation*}
    with equality holding if and only if $ n\geq 4$ or $d_\nabla R=0$.
  \end{conditions}
\end{theorem}

\subsection{Outline of the application of Cartan's method}
Before describing partial results for the general case $W\ne 0$, we
sketch the arguments leading to the results above. 

Although ${\mathfrak g} \subset J^1(TM)$ is surjective, we have
${\mathfrak h}\ne 0$ and Theorem \ref{neoGenerator} does not apply. In
\ref{first} we show that $\mathfrak g$ is already $\Theta
$-reduced. The associated lower coboundary morphism is not injective
and Theorem \ref{goosey} is therefore not applicable. We turn then, in
\ref{healthy} and \ref{third}, to the prolongation ${\mathfrak
  g}^{(1)}\subset J^1 {\mathfrak g} $ of $\mathfrak g$. Now
${\mathfrak g}^{(1)}$ is surjective (because $\mathfrak g$ is
$\Theta $-reduced) but has non-trivial structure kernel ${\mathfrak
  h}^{(1)}\cong T^*\!M $.

We show in \ref{expectation} that ${\mathfrak g}^{(1)} $ is already
$\Theta$-reduced when $W=0$. The coboundary morphism associated with
${\mathfrak g}^{(1)}$ is injective and Theorem \ref{goosey} applies to
${\mathfrak g}^{(1)} $, making it an associated Cartan algebroid.

\subsection{The $W\neq 0$ case and intransitivity}\lab{thc} 
If $W\ne 0$, then ${\mathfrak g}^{(1)} $ is no longer $\Theta
$-reduced. According to Proposition \ref{expectation} below, the
$\Theta $-reduction of ${\mathfrak g}^{(1)} $ is the preimage of
${\mathfrak g}_W$ under the natural projection ${\mathfrak g}^{(1)}
\rightarrow {\mathfrak g} $. In particular, this reduction has
${\mathfrak g}_W \ne {\mathfrak g}$ as image and is consequently not
surjective. Rather than continue to apply Cartan's method \`a la
\ref{hopeful} (turning next to the elementary reduction of the
preimage of ${\mathfrak g}_W$), we remark that ${\mathfrak g}_W$
itself is necessarily a reduction of ${\mathfrak g} $ (assuming
rank-constancy), and applying the algorithm to ${\mathfrak g}_W$ is an
easier prospect. 

Suppose that $W$ vanishes nowhere. Then the structure kernel of
$\mathfrak g_W$ is contained within that of ${\mathfrak g}_\sigma$
(because $\identity\cdot W=-2W $).  Since the upper (=lower)
coboundary morphism for ${\mathfrak g}_\sigma $ is injective, the same
is true for ${\mathfrak g}_W$; the prolongation of ${\mathfrak g}_W$
(or the prolongation of any reduction of ${\mathfrak g}_W$) will have
trivial structure kernel. Singularities not withstanding, $\mathfrak
g$ will therefore have an associated Cartan algebroid that is some
reduction of ${\mathfrak g}_W$, and in particular will be a
subalgebroid of $J^1(TM)$ (in contrast to the $W=0$ case). A detailed
argument is omitted here.

However we proceed, the following result implies that any associated
Cartan algebroid will be intransitive in general: Call $W$ {\df
  strongly degenerate} if there exists a section $\phi \subset
{\mathfrak h}$ such that $\nabla_VW = \phi \cdot W$, i.e., such that
\begin{multline*}
(\nabla_VW)(U_1,U_2)U_3=\phi W(U_1,U_2)U_3\\
-W(\phi U_1,U_2)U_3 -W(U_1,\phi U_2)U_3-W(U_1,U_2)(\phi U_3),  
\end{multline*}
for all vector fields $V,U_1,U_2,U_3$. We shall see in the proof of
the following that this definition is independent the metric within the conformal class used to fix a Levi-Cevita connection $\nabla $.
\begin{theorem}
  The isotropy ${\mathfrak g}_W \subset J^1(TM)$ is surjective
  (=transitive) if and only if $W$ is strongly degenerate.
\end{theorem}

The remainder of this section is devoted to proofs of the three
preceding theorems.

\subsection{The torsion reduction of $\mathfrak{g}$}\lab{first}
Since $\mathfrak{g}\subset J^1(TM)$, $\Theta$-reduction is the same as
torsion reduction. To compute it, we turn to the upper (=lower)
coboundary morphism for $\mathfrak{g}$,
\begin{equation*}
  T^*\!M\otimes \mathfrak{h}\xrightarrow{\Delta}\alternating^2(TM)\otimes TM.
\end{equation*}
Its restriction to $T^*\!M\otimes{\mathfrak h}_\sigma$ is nothing but
the upper coboundary morphism for ${\mathfrak g}_\sigma$.  Since the
latter is an isomorphism (see \ref{illustration}) the former is
surjective. In particular, $H(\mathfrak{g})=0$, implying $\mathfrak{g}$
is already torsion-reduced. Theorem \ref{goosey2} does not apply,
however, because $\Delta $ has non-trivial kernel. Indeed, counting
dimensions, we have
\begin{equation}
        \rank (\kernel\Delta)=\rank (TM).\mathlab{zx2}
\end{equation}

\subsection{The first prolongation ${\mathfrak g}^{(1)}$}\lab{healthy}
Since $\mathfrak{g}$ is torsion-reduced (and hence $\Theta$-reduced)
the prolongation $\mathfrak{g}^{(1)}$ is surjective (Proposition
\ref{koota}). By Proposition \ref{qlook}, its structure kernel
$\mathfrak{h}^{(1)}$ is $\kernel\delta=\kernel\Delta$.  Define a map
\begin{align*}
        T^*\!M&\xrightarrow{ i}\symmetric^2 (TM)\otimes TM\\
        i(\alpha)&:= j_{\mathrm S}(\alpha)-\sigma\otimes\sigma^{-1} (\alpha),
\end{align*}
where $j_\mathrm{S}\colon T^*\!M\rightarrow\symmetric^2 (TM)\otimes 
TM$ is the canonical morphism defined by
\begin{equation*}
  j_\mathrm{ S}(\alpha) (V_1, V_2)=\alpha (V_1)V_2+\alpha (V_2)V_1.
\end{equation*}
Then $ i$ is a monomorphism of ${\mathfrak g}$-representations 
($\dimension M\ge 2$). Since $$ i(\alpha) V=\sskew (\alpha\otimes 
V)+\alpha (V)\identity_{TM};\qquad V\subset TM,$$ we have $ 
i(T^*\!M)\subset T^*\!M\otimes{\mathfrak h}$. Therefore
\begin{equation*}
         i(T^*\!M)\subset (T^*\!M\otimes{\mathfrak h})\cap (\symmetric^2 
         (TM)\otimes TM)=\kernel\Delta=\mathfrak{h}^{(1)}.
\end{equation*}
Invoking \eqrefs{first}{zx2}, we conclude that
$\mathfrak{h}^{(1)}=i(T^*\!M)$. 

\subsection{A generator for ${\mathfrak g}^{(1)}$}\lab{third}%
To obtain a generator for $\mathfrak{g}^{(1)}$ we apply Theorem
\ref{caviar}. Recall that the Levi-Cevita connection $ \nabla $
generates $\mathfrak{g}$. We have $\tilde{\tilde \Theta }=0$ because
$H({\mathfrak g})=0$. We may therefore take $\epsilon=0$ in
\ref{caviar} and, using $\nabla$ to identify $\mathfrak{g}$ with $
TM\oplus \mathfrak{h} $, obtain \begin{equation}
        \nabla_U^{(1)} (V\oplus\phi):=(\nabla_UV+\phi (U))\oplus 
        (\nabla_U\phi+\curvature\nabla (U, V)).\mathlab{home2}
\end{equation}
We compute, with the help of Bianchi's second identity,
\begin{equation*}
        \curvature\nabla^{(1)} (U_1, U_2)(V\oplus\phi)=0\oplus\Big(\, 
        -(\nabla_V\curvature\nabla+\phi\cdot\curvature\nabla) (U_1, U_2)\,\Big).
\end{equation*}
This formula may also be written
\begin{equation}
  \curvature\nabla^{(1)}(U_1,U_2)X=-(X\cdot\curvature\nabla)(U_1,U_2)
  \subset \mathfrak{h}_\sigma;\quad 
  X\subset{\mathfrak g}.\mathlab{mornings}
\end{equation}

\subsection{The $\Theta$-reduction of $\mathfrak{g}^{(1)}$}\lab{expectation}
Since $\mathfrak{g}^{(1)}\subset J^1 \mathfrak{g} $ is surjective and
$\mathfrak{g}$ is transitive, the $\Theta$-reduction of
$\mathfrak{g}^{(1)}$ is the kernel of a morphism
$\mathfrak{g}^{(1)}\rightarrow h(\mathfrak{g}^{(1)}),$ which we denote
by $\Theta^{(1)}$, to distinguish it from the corresponding morphism
$\Theta\colon \mathfrak{g} \rightarrow h(\mathfrak{g})$ for
$\mathfrak{g}$; see \ref{qlook}.  The definition of
$h(\mathfrak{g}^{(1)})$ depends on the lower coboundary morphism for
$\mathfrak{g}^{(1)}$, which we denote by
\begin{equation*}
        T^*\!M\otimes{\mathfrak 
        h}^{(1)}\xrightarrow{\delta^{(1)}}\alternating^2 
        (TM)\otimes{\mathfrak g}.
\end{equation*}
Identifying ${\mathfrak h}^{(1)}$ with $T^*\!M$ as described above, 
one shows that $\delta^{(1)}$ is the map
\begin{equation*}
        \alpha\otimes\beta\mapsto\coRicci (\alpha\otimes\beta)+ 
        (\alpha\wedge\beta)\otimes\identity_{TM}.
\end{equation*}
Note that the first term on the right belongs to $\alternating^2 
(TM)\otimes{\mathfrak h}_\sigma$ 
and the second to $\alternating^2 
(TM)\otimes\langle\identity_{TM}\rangle$.  In particular, the image of 
$\delta^{(1)}$ lies entirely within $\alternating^2 
(TM)\otimes{\mathfrak h}$.  

Since $\coRicci$ is injective ($n\ge 3$) we have
$\kernel\delta^{(1)}=0$.  Therefore the second prolongation
${\mathfrak g}^{(2)}:=({\mathfrak g}^{(1)})^{(1)}$ of ${\mathfrak g}$
has trivial structure kernel (Proposition \ref{qlook}). In particular,
$h(\mathfrak{g}^{(1)}):=(\alternating^2(TM)\otimes\mathfrak{g})/
\image\delta^{(1)}$ has constant rank.

Next, we observe that the composite morphism of 
${\mathfrak g}$-representations,
\begin{equation*}
         E_\mathrm{Weyl}\monomorphism\alternating^2(TM)
         \otimes{\mathfrak h}_\sigma
         \monomorphism\alternating^2 (TM)\otimes{\mathfrak 
         g}\xrightarrow{/\image\delta^{(1)}}h(\mathfrak{g}^{(1)})
\end{equation*}
is injective.  This follows from the description of $\delta^{(1)}$ 
above, and
\begin{gather*}
         E_\mathrm{Weyl}\cap E_\mathrm{Ricci}=0,\\
         \text{where}\quad E_\mathrm{Ricci}=\coRicci (\symmetric^2 (TM)).
\end{gather*}
Identifying $ E_\mathrm{Weyl}$ with the corresponding ${\mathfrak
  g}$-subrepresentation of $h(\mathfrak{g}^{(1)})$, we have:
\begin{proposition}
The following diagram commutes:
\begin{equation*}
  \begin{CD}
    {\mathfrak g}^{(1)} @>\Theta^{(1)}>>h(\mathfrak{g}^{(1)})\\
    @V\text{{\normalfont projection}}VV @AA\text{{\normalfont inclusion}}A\\
    {\mathfrak g}@>>X\mapsto -X\cdot W> E_\mathrm{Weyl}
  \end{CD}\qquad.
\end{equation*}
\end{proposition}
\noindent Using the fact that ${\mathfrak g}^{(1)}$ and $\Theta^{(1)}$
are invariants of ${\mathfrak g}$, together with the fact that
$\identity_{TM}$ is a section of ${\mathfrak g}$, one deduces Theorem
\ref{tha}. The proposition also shows that the $\Theta$-reduction of
$\mathfrak{g}^{(1)}$ (the kernel of $\Theta^{(1)}$) is the preimage of
$\mathfrak{g}_W$ under the projection
$\mathfrak{g}^{(1)}\rightarrow\mathfrak{g}$. 

To prove Theorem \ref{thc}, apply Lemma \ref{algebraicLemma} to the
morphism $X\mapsto W\cdot X\colon\mathfrak g\rightarrow
E_\mathrm{Weyl}$, to show that ${\mathfrak g}_W$ has image $D\subset
TM$ (a possibly singular distribution on $M$), where $D$ is the kernel
of the morphism
\begin{gather*}
        TM\rightarrow E_\mathrm{Weyl}/(\mathfrak h\cdot W)\\
        V\mapsto\nabla_VW\modulo\mathfrak h\cdot W.
\end{gather*}
Recall here that $\mathfrak h\subset T^*\!M\otimes TM$ is the vector
bundle whose fiber ${\mathfrak h}(m)$ over $ m$ is the Lie algebra of
all infinitesimally conformal endomorphisms of $T_mM$.  This Lie
algebra acts on $ E_\mathrm{Weyl}(m)$ and
\begin{equation*}
                {\mathfrak h}\cdot W=\underset{m\in M}{\bigcup}\{\phi\cdot W(m)
                \suchthat\phi\in{\mathfrak h}(m)\}.
\end{equation*}
Evidently, $W$ is strongly degenerate if and only if $D=TM$.
\begin{proof}[Proof of proposition]
  We will apply part \eqref{ty1} of Theorem \ref{jit}, with the roles
  of
  $\mathfrak{g},\mathfrak{t},\mathfrak{h},\delta,\Theta,\tilde\Theta,
  \mathfrak{g}^{(1)}$ in the theorem being played by
  $\mathfrak{g}^{(1)},\mathfrak{g},\mathfrak{h}^{(1)},\delta^{(1)},
  \Theta^{(1)},\tilde\Theta^{(1)},\mathfrak{g}^{(2)}$.
  
  Our first task is to choose a connection
  $\nabla^{\mathfrak{h}^{(1)}}$ on $\mathfrak{h}$.  The Levi-Cevita
  connection $\nabla$ on $ TM$ determines a linear connection on
  $T^*\!M\otimes (T^*\!M\otimes TM)$ and one has a chain of inclusions
  \begin{equation*}
        {\mathfrak h}^{(1)}\subset T^*\!M\otimes{\mathfrak h}\subset 
        T^*\!M\otimes (T^*\!M\otimes TM),
  \end{equation*}
  which we claim are $\nabla$-invariant.  The $\nabla$-invariance of
  ${\mathfrak h}\subset T^*\!M\otimes TM$ follows from Proposition
  \eqrefs{reconstruction}{gt}. So the second inclusion is indeed
  $\nabla$-invariant. Because $\nabla$ generates $\mathfrak{g}$ and
  because
  \begin{equation*}
        \Delta\colon T^*\!M\otimes{\mathfrak h}\rightarrow\alternating^2 
        (TM)\otimes TM
  \end{equation*}
  is ${\mathfrak g}$-equivariant, it follows that $\Delta$ is
  $\bar\nabla$-equivariant (by Proposition \ref{neoAssociated}). But
  $\nabla$ is torsion free, meaning $\bar\nabla$-invariance is the
  same as $\nabla$-invariance. So the kernel ${\mathfrak
    h}^{(1)}\subset T^*\!M\otimes{\mathfrak h}$ of $\Delta$ must be
  $\nabla$-invariant, as claimed. 
  
  We choose $\nabla^{{\mathfrak h}^{(1)}}$ to be the connection that 
  ${\mathfrak h}^{(1)}$ inherits from $T^*\!M\otimes{\mathfrak h}$ as 
  a $\nabla$-invariant subbundle.  Appealing to \eqrefs{third}{home2} 
  and the fact that ${\mathfrak h}^{(1)}\subset\symmetric^2 
  (TM)\otimes TM$, one can show that
  \begin{equation}
         \Big(\,d_{\nabla^{(1)}}\phi-\delta^{(1)} (\nabla^{{\mathfrak 
        h}^{(1)}}\phi)\,\Big)(U_1, U_2)=\phi (\torsion\nabla (U_1, 
        U_2))=0,\mathlab{laugh}
  \end{equation}
  for all sections $\phi\subset{\mathfrak h}^{(1)}\subset 
  T^*\!M\otimes{\mathfrak h}\subset T^*\!M\otimes{\mathfrak g}$ and $ 
  U_1, U_2\subset TM$.
  
  In the present context \eqrefs{jit}{box1} reads
  \begin{equation*}
        \tilde\Theta^{(1)} (X\oplus\phi):=\curvature\nabla^{(1)}(\,\cdot\,,
        \,\cdot\,)X-d_{\nabla^{(1)}}\phi+\delta^{(1)} (\nabla^{{\mathfrak 
        h}^{(1)}}\phi).
  \end{equation*}
  From \eqrefs{third}{mornings} and \eqref{laugh} above one obtains
  \begin{equation}
        \tilde\Theta^{(1)} (X\oplus\phi)=-X\cdot\curvature\nabla
        =-X\cdot W+\coRicci (X\cdot R),\mathlab{sqa}
  \end{equation}
  for arbitrary sections $ X\subset{\mathfrak g} $ and 
  $\phi\subset{\mathfrak h}^{(1)}$. Since 
  $\Theta^{(1)}\colon{\mathfrak 
  g}^{(1)}\rightarrow h(\mathfrak{g}^{(1)})$ is the composite
  \begin{equation*}
        {\mathfrak g}^{(1)}\cong{\mathfrak g}\oplus{\mathfrak 
        h}^{(1)}\xrightarrow{\tilde\Theta^{(1)}}\alternating^2 
        (TM)\otimes{\mathfrak g}\xrightarrow{/\image\delta^{(1)}}
        h(\mathfrak{g}^{(1)}),
  \end{equation*}
  and because the image of $\delta^{(1)}$ contains the image of 
  $\coRicci$, this completes the proof.
\end{proof}

\subsection{The $W=0$ case}\lab{dutiful}
By the above, $\mathfrak{g}^{(1)}$ is $\Theta$-reduced if and only if
$\mathfrak{g}_W=\mathfrak{g}$.  But as $\identity_{TM}\subset T^*\!M
\otimes TM\subset J^1 (TM)$ is a section of ${\mathfrak g}$, this is
clearly equivalent to the vanishing of $ W$. 

Assuming $W=0$, Theorem \ref{goosey} applies (with $\mathfrak{g}$ and
$\mathfrak{t}$ in the theorem replaced with $\mathfrak{g}^{(1)}$ and
$\mathfrak{g}$), establishing our claim that $\mathfrak{g}^{(1)}$ is
the Cartan algebroid associated with $\mathfrak{g}$ when $W=0$. To
compute the Cartan connection $\nabla^{(2)}$ on $\mathfrak{g}^{(1)}$
we will now apply parts \eqref{ty1c} and \eqref{ty3} of
Theorem \ref{jit}.

In the present context, the generator equation defined in
\eqrefs{jit}{ty1c} reads

\begin{equation*}
        \delta^{(1)} (\epsilon)+\coRicci (X\cdot R)=0.
\end{equation*}
We have used \eqref{sqa} above. Referring to the
description of $\delta^{(1)}$ in \ref{expectation}, we see that a
solution is given by $\epsilon=-X\cdot R$. Using $\nabla^{(1)}$ to
identify $\mathfrak{g}^{(1)}$ with $\mathfrak{g}\oplus
\mathfrak{h}^{(1)}$, and keeping in mind the identification
$\mathfrak{h}^{(1)}\cong T^*\!M $ implicit above, we deduce from
\eqrefs{jit}{ty3},
\begin{equation}
  \nabla_U^{(2)}(X\oplus\alpha)=\Big(\,\nabla_U^{(1)} X+\sskew (\alpha\otimes 
  U)+\alpha (U)\identity_{TM}\,\Big)\oplus (\,\nabla_U\alpha+(X\cdot R)U\,),
  \mathlab{honda}
\end{equation}
for arbitrary sections $ X\subset{\mathfrak g}$ and $\alpha\subset 
T^*\!M$. 

We claim
\begin{equation}
        \curvature\nabla^{(2)}(U_1,U_2)(X\oplus\alpha)=
                (X\cdot d_\nabla R)(U_1,U_2),\mathlab{fledgling}
\end{equation}
where $d_\nabla R$ is the Cotton-York tensor, defined in
\ref{outline}.  In particular, whenever $ W=0$, the tensor $d_\nabla
R$ is a conformal invariant which vanishes if and only if
$\nabla^{(2)}$ is flat, i.e., if and only if the Cartan algebroid
${\mathfrak g}^{(1)}$ is flat (Theorem \ref{dissymmetry}). This
completes the proof of Theorem \ref{thb}.
\begin{proof}[Proof of  \eqref{fledgling}]
  Since $ W=0$ we have $\curvature\nabla=\coRicci (R)$ and, with a 
  little effort, one computes
  \begin{align*}
        \curvature\nabla^{(2)}(U_1,U_2)(X\oplus\alpha)
        &=-((\nabla_{U_1}^{(1)}X)\cdot R)U_2-\nabla_{U_1}((X\cdot R)U_2)\\
        &+((\nabla_{U_2}^{(1)}X)\cdot R)U_1+\nabla_{U_2}((X\cdot R)U_1)\\
        &-(X\cdot R)([U_1, U_2]).
  \end{align*}
  Equation \eqref{fledgling} now follows from the readily verified 
  identities
  \begin{align*}
         (\nabla_U^{(1)}X)\cdot V&=\nabla_U(X\cdot V)-X\cdot 
         (\nabla_UV)+\nabla_{X\cdot U}V\enspace\text{for $ X 
         \subset{\mathfrak g},\, V\subset TM$};\\
         (\nabla_U^{(1)}X)\cdot \alpha&=\nabla_U(X\cdot\alpha)-X\cdot 
         (\nabla_U\alpha)+\nabla_{X\cdot U}\alpha\enspace\text{for $ X 
         \subset{\mathfrak g},\, \alpha\subset T^*\!M$}.
  \end{align*}
  One also makes use of the fact that $\torsion\nabla=0$.
\end{proof}

\appendix%
\section{Cartan groupoids and Lie pseudogroups}\lab{global}
We now explain how {\em flat} Cartan algebroids may be viewed as
infinitesimal versions of Lie pseudogroups; and conversely, how Lie
pseudogroups integrate flat Cartan algebroids.  As a byproduct of this
discussion, we are led to define {\df Cartan groupoids}.  These are
the global versions of Cartan algebroids and may be viewed as
deformations of Lie pseudogroups.  Flat Cartan groupoids appear in
\cite{Tang_04} where they are called `groupoid etalifications.'

\subsection{Lie pseudogroups via pseudoactions}\lab{gaily}
Let us explain, in invariant groupoid language, what it means for a 
pseudogroup to be a {\em Lie} pseudogroup.  For the classical 
description see, e.g., \cite{Stormark_00}.  

A group of transformations in a smooth manifold $M$ is a Lie group of 
transformations if it arises from the (smooth) action of some abstract 
Lie group.  Analogously, we declare an arbitrary pseudogroup of 
transformations in $M$ to be a {\df Lie pseudogroup} if it is arises 
from the pseudoaction of some Lie {\em groupoid}.  It remains to 
explain what we mean by pseudoactions and the pseudogroups of 
transformations they define.  We shall understand all constructions to 
be made in the smooth category.

 Let $G$ be a Lie groupoid over $M$.  Call an immersed submanifold 
 $\Sigma\subset G$ a {\df pseudotransformation} if the restrictions to 
 $\Sigma$ of the groupoid's source and target maps are local 
 diffeomorphisms.  In other words, each point of $\Sigma$ should have 
 an open neighborhood in $\Sigma$ that is a (smooth) local bisection 
 of $G$.  For example, the pseudotransformations of the pair groupoid 
 $ M\times M$ are the local transformations in $M$ taking possibly 
 multiple values.

A {\df pseudoaction} of $G$ on $M$ is  any foliation ${\mathcal 
F}$ on $G$ such that:
\begin{conditions}
        \item  The leaves of ${\mathcal  F}$ are pseudotransformations.
        \item ${\mathcal  F}$ is  {\df multiplicatively closed}.\lab{potato}
\end{conditions}
To define what is meant in \eqref{potato} let $\hat{\mathcal F}$ denote 
the collection of those subsets of $G$ that are simultaneously an open 
subset of some leaf of ${\mathcal F}$, and a local bisection.  Let $ 
\hat G$ denote the collection of {\em all} local bisections of $G$, this 
being a groupoid over the power set of $M$.  Then condition 
\eqref{potato} is the requirement that $\hat{\mathcal F}\subset \hat G$ be 
a subgroupoid.

Given a pseudoaction ${\mathcal F}$ of $G$ on $M$, each element of 
$\hat{\mathcal F}$ defines a local diffeomorphism in $M$ and, by 
\eqref{potato}, the collection of all such local diffeomorphisms 
constitutes a pseudogroup of transformations in $M$.  For example, if 
$G$ is an action groupoid $ G=G_0\times M$, then the canonical 
horizontal foliation ${\mathcal F}$ furnishes us with the usual 
pseudogroup of transformations associated with the prescribed action 
of the Lie group $ G_0$.  

\subsection{The flat Cartan algebroid associated with a Lie 
pseudogroup}\lab{jury}%
Let ${\mathcal G}$ be a Lie pseudogroup of transformations in $M$.  
Then ${\mathcal G}$ is generated by the pseudoaction ${\mathcal F}$ of 
some Lie groupoid $G$ over $M$.  Define $\hat{\mathcal F}$ as in 
\ref{gaily} above.  Then each point $ g\in G$ lies in some bisection $ 
b\in\hat{\mathcal F}$ and all such bisections have the same one-jet at 
$g$.  Thus ${\mathcal F}$ defines a map $ D_{\mathcal F}\colon 
G\rightarrow J^1 G$ into the Lie groupoid of all one-jets of 
bisections of $G$.  This map, which is a right inverse for the natural 
projection $ J^1 G\rightarrow G$, is a groupoid morphism because 
${\mathcal F}$ is multiplicatively closed.

An arbitrary groupoid morphism $D\colon G\rightarrow J^1 G$ furnishing 
a right inverse for $ J^1 G\rightarrow G$ is what we call a {\df 
Cartan connection} on $G$.  These connections may be viewed as certain 
`multiplicatively closed' distributions on $G$.  The connection $ D$ 
is Frobenius integrable precisely when it comes from a pseudoaction 
${\mathcal F}$ as above, in which case $ D$ is simply the tangent 
distribution.  A Lie groupoid equipped with a (possibly 
non-integrable) Cartan connection is a {\df Cartan groupoid}.  Thus 
Cartan groupoids are deformed Lie pseudogroups.

Differentiating a Cartan connection $ D\colon G\rightarrow J^1 G$, we 
obtain a splitting ${\mathfrak g}\rightarrow J^1{\mathfrak g}$ 
for the exact sequence of Lie algebroids
\begin{equation}
         0\rightarrow T^*\!M\otimes{\mathfrak g}\monomorphism J^1{\mathfrak 
         g}\rightarrow{\mathfrak g}\rightarrow 0.\mathlab{here}
\end{equation}
This splitting will be a morphism of Lie algebroids, i.e., amounts to
a Cartan connection $\nabla$ on ${\mathfrak g}$, as defined in
\ref{algebra}. If $ D=D_{\mathcal F}$ for some pseudoaction
${\mathcal F}$ as above, then $D$ is Frobenius integrable and we claim
that this implies $\curvature \nabla = 0$, so that ${\mathfrak g} $ is
flat (Theorem \ref{dissymmetry}); ${\mathfrak g}$ is then a flat
Cartan algebroid associated with the pseudogroup ${\mathcal G}$
generated by ${\mathcal F}$.

\subsection{The Lie pseudogroup integrating a flat Cartan algebroid}
Let ${\mathfrak g}$ be a Cartan algebroid over $M$ with Cartan
connection $\nabla$ and assume ${\mathfrak g}$ is the Lie algebroid of
some Lie groupoid $G$.  For simplicity, suppose $G$ has connected
source-fibers.  The connection $\nabla$ determines a Lie algebroid
morphism ${\mathfrak g}\rightarrow J^1{\mathfrak g}$ splitting the
exact sequence \eqrefs{jury}{here}.  By the groupoid version of Lie's
Second Theorem, this morphism integrates to a groupoid morphism $
D\colon G\rightarrow J^1 G$, i.e., to a Cartan connection on the Lie
groupoid $G$.  Supposing $\mathfrak g$ is flat, we have $\curvature
\nabla = 0$ (Theorem \ref{dissymmetry}), and we claim this guarantees
that $ D$ is Frobenius integrable.  The integrating foliation
${\mathcal F}$ is a pseudoaction generating a Lie pseudogroup
${\mathcal G}$ of transformations in $M$.

For each locally defined $\nabla$-parallel section $X\subset{\mathfrak 
g}$, the vector field $\#X\subset TM$ integrates to a one-parameter 
family of local transformations belonging to $\mathcal G$.  Conversely 
each transformation in the pseudogroup $\mathcal G$ --- or at least 
each transformation `close' to the identity --- arises as the time-one 
map associated with such a vector field.  In this sense $\mathcal G$ 
integrates the flat Cartan algebroid ${\mathfrak g}$.

\section{Miscellany}\lab{miscellany}
\subsection{On morphisms whose domains sit in a short exact sequence}
\lab{algebraicLemma}%
       In the category of vector spaces, or of vector bundles over $M$, let 
         $\theta\colon B\rightarrow B_1$ be an arbitrary morphism, $ B_0$ 
         its kernel, and  suppose $ B$ occurs in some exact sequence, as shown below:
        \begin{equation*}
          \begin{CD}
                 ~ @.   ~ @.   B_0  @.   ~ @.   ~\\
                   @.     @.   @VVV @.     @.    \\
                 0 @>>> A @>>> B    @>>> C@>>>0  \\
                   @.     @.   @VV\Theta V @.     @.   \\
                 ~ @.   ~ @.   B_1  @.   ~ @.   ~
          \end{CD}\enspace.
        \end{equation*}
The proof of the following is a straightforward diagram chase.
\begin{lemma}
        Let $ A_0$ and $ A_1$ denote, respectively, the kernel and
        image of the composite morphism
        $A\monomorphism B\xrightarrow{\theta}B_1$; and define
        $C_1:=B_1/A_1$, so that the sequence
        \begin{equation*}
          0\rightarrow A_1\monomorphism B_1\rightarrow 
          C_1\rightarrow 0          
        \end{equation*}
        is also exact. Then:
        \begin{conditions}
        \item There exists a unique morphism $C\xrightarrow{\Theta}
          C_1$ such that
                \begin{equation*}
                        \begin{CD}
                                B @>>> C\\
                                @VV\theta V @VV\Theta V \\
                                B_1 @>>> C_1
                        \end{CD}
                \end{equation*}
               is commutative.
             \item If $ C_0$ denotes the kernel of $\Theta$, and
               $B_0\monomorphism C_0$ the restriction of
               $B\monomorphism C$, then the top row in the following
               commutative diagram is exact (in addition to the other
               two rows):
               \begin{equation*}
                 \begin{CD}
                   0 @>>> A_0 @>>> B_0 @>>>C_0 @>>> 0\\
                   @. @VVV @VVV @VVV @.\\
                   0 @>>> A @>>> B @>>>C @>>> 0 \\
                   @. @VVV @VV\theta V @VV\Theta V @.\\
                   0 @>>> A_1 @>>> B_1 @>>>C_1 @>>> 0 \enspace.
                 \end{CD}%
               \end{equation*}
        \end{conditions}
\end{lemma}
\subsection{Gluing Lie algebroid `point' invariants into global
invariants}\lab{extensionLemma}%
Let ${\mathfrak g}$ be a transitive Lie algebroid over $M$ and 
${\mathfrak h}\subset{\mathfrak g}$ the kernel of its anchor.  Let $E$ 
be a ${\mathfrak g}$-representation.  Each fiber ${\mathfrak h}(m)$ of 
${\mathfrak h}$ is a Lie algebra acting on the vector space $ E(m)$.  
The following lemma furnishes conditions under which the the existence 
of ${\mathfrak h}(m)$-invariant elements $\sigma (m)\in E(m)$, for 
each $ m\in M$, implies the existence of {\em global} ${\mathfrak 
g}$-invariant sections $\sigma\subset E$.  For applications, see 
\ref{Riemannian}.
\begin{lemma}[Extension Lemma]
  Suppose that $M$ is simply-connected and that the set $E^{\mathfrak
    h}\subset E$ of ${\mathfrak h}$-invariant elements has constant
  rank $ r>0$.  Then $E$ possesses a non-vanishing ${\mathfrak
    g}$-invariant section $\sigma$.  If $ r=1$, then $\sigma$ is
  unique up to constant.
\end{lemma}
\begin{proof}
Noting that $ Y\subset{\mathfrak h}$ implies $[ X, Y]_{\mathfrak 
g}\subset{\mathfrak h}$, the identity
        \begin{equation*}
                 Y\cdot (X\cdot\sigma)=X\cdot (Y\cdot\sigma)-[X, 
        Y]_{\mathfrak g}\cdot\sigma;\qquad X\subset{\mathfrak 
        g},\,Y\subset{\mathfrak h},\,\sigma\subset E, 
\end{equation*}
shows that the rank-$r$ subbundle $ E^{\mathfrak h}\subset E$ is
${\mathfrak g}$-invariant.  Because ${\mathfrak h}$ acts trivially on
$ E^{\mathfrak h}$, the representation ${\mathfrak
  g}\rightarrow\gl(E^{\mathfrak h})$ factors through the anchor,
delivering a representation $ TM\rightarrow\gl(E^{\mathfrak h})$,
i.e., a {\em flat} linear connection $D$ on $E^{\mathfrak h}$.  One
takes $\sigma$ to be any non-vanishing $D$-parallel section of
$E^{\mathfrak h}$, whose existence is guaranteed by flatness and the
simple-connectivity of $M$.  The uniqueness claim is clear.
\end{proof}

\subsection{Proof of Proposition \ref{normtor}}\lab{possum}
Let $\nabla$ be {\em any} generator of $\mathfrak{g}$.  Then as
$\mathcal{H},\sigma,\mathbf{n},\theta$ and $d\theta$ are all
$\mathfrak{g}$-invariant, they are all ${\bar\nabla} $-invariant
(Proposition \ref{neoAssociated}). From the $\bar\nabla$-invariance of
$\sigma$ and $\mathbf{n} $, one immediately computes
\begin{align}
         (\nabla_U\sigma) (V_1, V_2)&=\sigma \Big(\,(\torsion\bar\nabla 
         (U))_\sym V_1, V_2\,\Big)\mathlab{week1}\\
\text{and}\quad\nabla{\mathbf 
         n}&=\torsion\bar\nabla ({\mathbf n}),\mathlab{week2}
\end{align}
where $\torsion\bar\nabla (U):=\torsion\bar\nabla(U,\,\cdot\,)\subset
T^*\!M\otimes TM$ and $ U, V_1, V_2$ are arbitrary vector fields on
$M$.  Here and in the sequel a subscript $\sym$ (or $\alt$) on any
type of 2-tensor indicates its symmetrization
(resp.~skew-symmetrization), defined using the metric $\sigma$ as
appropriate. 
% For example, if $\phi\in T^*\!M\otimes T^*\!M$ (or
% $\phi\in \mathcal{H}^*\otimes \mathcal{H}^*$), then
% \begin{equation*}
%              \phi_\sym(V_1,V_2):=\phi (V_1, V_2)+\phi (V_2, V_1).
% \end{equation*}
% On the other hand, if $\psi\in T^*\!M\otimes TM$, for example, then $\psi_\sym$ denotes the image of $\psi$ under the composite
% \begin{equation*}
%   T^*\!M\otimes TM\xrightarrow{\identity\otimes\sigma}
%   T^*\!M\otimes T^*\!M\xrightarrow{\phi\mapsto\phi_\sym}
%   T^*\!M\otimes T^*\!M\xrightarrow{\identity\otimes\sigma^{-1}}
%   T^*\!M\otimes TM.
% \end{equation*}
For any 2-tensor $\phi$, we have $\phi=(\phi_\sym+\phi_\alt)/2$.

From \eqref{week2} it follows that $\nabla_\mathbf{n}
\mathbf{n}=0$. From \eqref{week2} and the $\bar\nabla $-invariance of $\theta$, we compute,
\begin{equation*}
  \theta(\nabla_V\mathbf n)=\theta(\torsion\bar\nabla(\mathbf n, V))=d\theta(\mathbf n, V)=0;\quad V \subset TM.
\end{equation*}
So $\nabla_V \mathbf{n} $ is $\mathcal{H}$-valued, for any $V \subset TM$. This establishes \eqrefs{normtor}{mice2}. 

Now $\alternating ^2 (\mathcal{H})$ is rank-one and spanned by $dA$,
implying that the restriction of $\torsion \bar\nabla$ to
$\mathcal{H}$ (a section of $\alternating^2(\mathcal{H})$) is of the
form $dA\otimes V$, for some vector field $V \subset TM$. However, the
${\bar\nabla} $-invariance of $\theta$ gives,
\begin{equation*}
  \theta(\torsion {\bar\nabla}(U_1,U_2))=d\theta(U_1,U_2)=dA(U_1,U_2),
\end{equation*}
for arbitrary $U_1,U_2 \subset \mathcal{H} $, implying $\theta(V)=1$. We conclude that,
\begin{equation}
  (\torsion {\bar\nabla})| \mathcal{H}=dA\otimes(\mathbf{b}+\mathbf{n}),\mathlab{xx1}
\end{equation}
for some unique section $\mathbf{b}\subset \mathcal{H} $. Equation \eqref{week1} now gives
\begin{equation*}
  (\nabla_U \sigma)(V_1,V_2)= dA(U,V_1)\sigma(\mathbf{b},V_2)
  +dA(U,V_2)\sigma(\mathbf{b},V_1);\enspace U,V_1,V_2\subset\mathcal{H}.
\end{equation*}
It is not hard to see that this implies
\begin{equation}
  \mathbf{b}=0\iff \nabla \sigma | \mathcal{H}=0.\mathlab{xx}
\end{equation}
From \eqref{week2} and \eqref{xx1} we obtain,
\begin{equation*}
  \torsion {\bar\nabla}(U_1+a_1 \mathbf{n},U_2+a_2 \mathbf{n})=
  a_1 \nabla_{U_2}\mathbf{n}  - a_2 \nabla_{U_1}\mathbf{n} 
  + dA(U_1,U_2)(\mathbf{b}+\mathbf{n}).
\end{equation*}
Therefore, if $\nabla $ is a generator satisfying $\nabla \sigma|
\mathcal{H}=0$, as in \eqrefs{normtor}{mice3}, then
\eqrefs{normtor}{mice3b} becomes a consequence of \eqref{xx} above.

If $\nabla \mathbf{n}\subset ({\mathcal H}^* \otimes {\mathcal
  H})_\sym$ then \eqref{week2} implies that $(\torsion
{\bar\nabla}(\mathbf{n}))_{\sym}=\nabla \mathbf{n}$, so that
\eqrefs{normtor}{mice5} follows from \eqref{week1} (take
$U:=\mathbf{n}$).

We return to supposing that $\nabla $ is an arbitrary generator of
$\mathfrak{g}$. To establish the remaining claims of the proposition we
require a detailed analysis of the upper coboundary morphism,
\begin{equation}
  T^*\!M\otimes \mathfrak{h} \xrightarrow{\Delta}\alternating ^2(TM)\otimes TM.
  \mathlab{exit}
\end{equation}
By Proposition \ref{cr}, we have $\mathfrak{h} \cong
(\mathbb R \times M)$, so that $T^*\!M\otimes \mathfrak{h}\cong T^*\!M$. With
the help of the $\mathfrak{g}$-invariant splitting
$TM=\mathcal{H}\oplus\langle \mathbf{n} \rangle$ and a little
multilinear algebra, one identifies a natural isomorphism of $\mathfrak{g}$-representations,
\begin{equation}
  \alternating ^2(TM)\otimes TM \overset{\phi}{\cong}
        \alternating^2 (TM)\oplus(\mathcal{H}^*\otimes \mathcal{H})_\sym
        \oplus \mathcal{H}  \oplus (\mathbb R \times M),
  \mathlab{hello}  
\end{equation}
where $(\mathcal{H}^*\otimes \mathcal{H})_\sym \subset \mathcal{H}
^*\otimes \mathcal{H} $ denotes the $\mathfrak{g}$-subrepresentation
of symmetric elements. We write $\phi=\phi_1 \oplus \phi_2 \oplus
\phi_3 \oplus \phi_4$ and describe the component morphisms $\phi_j$ at
the end. Knowing the $\phi_j$, one readily establishes the following:
\begin{lemma}
  Under the identifications above, we have:
  \begin{conditions}
    \item\lab{piz1} The torsion $\torsion {\bar\nabla} \subset
  \alternating ^2 (TM)\otimes TM$ is given by, $$\torsion
  {\bar\nabla}=d\theta \oplus (\nabla \mathbf{n})_\sym \oplus
  \mathbf{b} \oplus f,$$ where $f \subset (\mathbb R \times M)$ is the
  function on $M$ defined by $(\nabla \mathbf{n} )_\alt U=fJU$.

    \item\lab{piz2} The upper coboundary morphism $\Delta$ takes the form
    \begin{equation*}
    \begin{matrix}
       T^*\!M & \xrightarrow{~~\Delta~~} & \alternating^2 (TM) & \oplus
              & (\mathcal{H}^*\otimes \mathcal{H})_\sym & \oplus & 
       \mathcal{H}      & \oplus &  (\mathbb R \times M)\\
       \alpha & \mapsto   & 0  & \oplus
         & 0                & \oplus & 
         -2\sigma^{-1}(\alpha)|\mathcal{H} & \oplus & \alpha ({\mathbf n})
    \end{matrix}\enspace.
    \end{equation*}
  \end{conditions}
\end{lemma}
\noindent In particular, $\Delta$ is injective, and its image has complement,
\begin{equation*}
  C:=\alternating ^2(TM)  \oplus (\mathcal{H} ^*  \oplus \mathcal{H})_\sym  \oplus 0  \oplus 0,
\end{equation*}
which is $\mathfrak{g}$-invariant because the splitting \eqref{hello}
is $\mathfrak{g}$-invariant. This establishes
\eqrefs{normtor}{mice1}. Also, we obtain $\mathfrak{g}$-invariant
isomorphisms,
\begin{equation*}
  H(\mathfrak{g})\cong C\cong \alternating ^2(TM)  \oplus (\mathcal{H}^*  \oplus \mathcal{H})_\sym.
\end{equation*}
This proves the first part of \eqrefs{normtor}{mice4}.

Now \eqref{piz1} shows that $\torsion {\bar\nabla} \subset C$ if and
only if $\mathbf{b}=0 $ and $f=0$. These are true if and only if
$\nabla \sigma | \mathcal{H}=0$ (by \eqref{xx}) and $\nabla \mathbf{n}
\subset (\mathcal{H} ^* \otimes \mathcal{H})_\sym. $ This finishes the
proof of \eqrefs{normtor}{mice3}.

With $\nabla$ fixed \`a la \eqrefs{normtor}{mice3}, we have $(\nabla
\mathbf{n})_\sym=\nabla \mathbf{n} $ and \eqref{piz1} gives,
\begin{equation*}
  \torsion {\bar\nabla} =d\theta  \oplus \nabla \mathbf{n}  \oplus 0  \oplus 0.
\end{equation*}
Under our identifications, the canonical projection $\alternating ^2(TM) \otimes TM\rightarrow H(\mathfrak{g})$ is just the map
\begin{align*}
  \alternating ^2(TM)  \oplus (\mathcal{H}^* \otimes \mathcal{H})_\sym \oplus \mathcal{H}  \oplus (\mathbb R \times M) &\rightarrow \alternating ^2(TM)  \oplus (\mathcal{H}^* \otimes \mathcal{H})_\sym\\
  \lambda \oplus \mu \oplus U \oplus s &\mapsto \lambda \oplus
  \mu.
\end{align*}
The corresponding map of section spaces sends $\torsion {\bar\nabla} $
to the intrinsic torsion $\tau$, so that
\begin{equation*}
 \tau=d\theta  \oplus \nabla  \mathbf{n}.
\end{equation*}
This finishes the proof of \eqrefs{normtor}{mice4} and the
proposition.

\vspace{\baselineskip}
{\itshape The definitions of $\phi_1,\phi_2,\phi_3,\phi_4$.} %
The morphism $\phi_1$ is the composite morphism
\begin{equation*}
        \alternating^2 (TM)\otimes TM\rightarrow\alternating^2 
        (TM)\otimes\langle{\mathbf n}\rangle\cong\alternating^2 (TM),
\end{equation*}
where the first arrow is  the identity on $\alternating^2 (TM)$ 
tensored with the orthogonal projection $ TM\rightarrow\langle{\mathbf 
n}\rangle$.  The morphism $\phi_2$ is the composite
\begin{equation*}
        \alternating^2 (TM)\otimes TM\rightarrow T^*\!M\otimes 
        TM\rightarrow  \mathcal{H}^*\otimes\mathcal{H}\rightarrow
        (\mathcal{H}^*\otimes\mathcal{H})_\sym,
\end{equation*}
where the first arrow is contraction $\rho\mapsto\rho ({\mathbf
  n},\,\cdot\,)$, the second arrow is obtained by tensoring the
restriction $ T^*\!M\rightarrow \mathcal{H}^*$ with orthogonal projection $
TM\rightarrow \mathcal{H}$, and the third arrow is symmetrization.  The
morphism $\phi_3$ is the composite
\begin{equation*}
        \alternating^2 (TM)\otimes TM\rightarrow\alternating^2 
        (\mathcal{H})\otimes \mathcal{H}=\langle dA\rangle\otimes\mathcal{H}\cong \mathcal{H},
\end{equation*}
where the first arrow is the restriction $\alternating^2 
(TM)\rightarrow\alternating^2 (\mathcal{H})$ tensored with orthogonal projection 
$TM\rightarrow \mathcal{H}$.  The morphism $\phi_4$ is the composite 
\begin{equation*}
        \alternating^2 (TM)\otimes TM\rightarrow\alternating^2 
        (\mathcal{H})\otimes\langle{\mathbf n}\rangle=\langle 
        dA\rangle \otimes \langle{\mathbf n}\rangle\cong \mathbb R \times M,
\end{equation*}
where the first arrow is restriction tensored with orthogonal 
projection.

\vspace{\baselineskip} {\itshape Relationship with the Levi-Cevita
  connection.} Although we have no need to do so, it is not difficult
to express the generator $\nabla $ fixed in \eqrefs{normtor}{mice3} in
terms of the Levi-Cevita connection $\levi$ associated with $\sigma$:
\begin{equation*}
        \nabla_UV=\levi_UV-\epsilon (V)U,
\end{equation*}
where $\epsilon\subset T^*\!M\otimes(T^*\!M\otimes TM)_\alt$ 
is defined by
\begin{align*}
         \epsilon ({\mathbf n})&=\frac{1}{2}(\levi{\mathbf n})_\alt,\\
         \epsilon (U)&=(\theta\otimes\levi_U{\mathbf 
                     n})_\alt\quad\text{for  $ U\subset \mathcal{H}$,}\\
\text{or}\quad \epsilon (U)V&=(J\levi_U{\mathbf 
n})\times V\quad\text{for  $U\subset \mathcal{H} $ and $V\subset TM$}.                              
\end{align*}
Here $\times$ denotes cross product and $(T^*\!M\otimes TM)_\alt
\subset T^*\!M \otimes TM$ denotes the
$\mathfrak{g}$-subrepresentation of skew-symmetric elements.

\subsection{On $J^2 {\mathfrak t} $ as a subbundle of $J^1(J^1
  {\mathfrak t}))$}\lab{pro}%
Here we prove Proposition \ref{boot} and Lemma
\ref{prolongation} which describe properties of the second jet bundle
$J^2 {\mathfrak t} $ of an arbitrary vector bundle ${\mathfrak t}$.

Evidently the formula for $\omega_1$ in Proposition \ref{boot}
defines, at the very least, a linear map of section spaces,
$\Gamma(J^1(J^1{\mathfrak t}))\rightarrow \Gamma(T^*\!M \otimes J^1
{\mathfrak t})$. Employing \eqrefs{boot}{Leib}, one shows that
$\omega_1(f \xi)=f \omega_1 \xi$ for any smooth function $f$. This
implies that the map of section spaces drops to a well defined
morphism of vector bundles. A similar argument applies to $\omega_2$.

It is not difficult to show that $J^2 {\mathfrak t} \subset
J^2_+{\mathfrak t}$ and, moreover, that $J^2 {\mathfrak t} \subset
\kernel \omega_2$. To finish the proof of the proposition it suffices
to show that $J^2 t$ and $\kernel \omega_2$ have the same rank.

Noting that $J^1(J^1{\mathfrak t})$ sits in an exact sequence,
\begin{equation*}
  0 \rightarrow T^*\!M \otimes J^1{\mathfrak t} \rightarrow J^1(J^1{\mathfrak t}) \rightarrow J^1{\mathfrak t} \rightarrow 0,
\end{equation*}
we apply Lemma \ref{algebraicLemma} to the morphism
$\omega_1\colon J^1(J^1{\mathfrak t})\rightarrow T^*\!M \otimes
{\mathfrak t}$ and derive an exact sequence,
\begin{equation*}
  0 \rightarrow T^*\!M \otimes T^*\!M \otimes {\mathfrak t} \rightarrow J^2_+ {\mathfrak t} \rightarrow J^1{\mathfrak t} \rightarrow 0.
\end{equation*}
Next, applying Lemma \ref{algebraicLemma} to the morphism
$\omega_2\colon J^2_+{\mathfrak t} \rightarrow
\alternating^2(TM)\otimes {\mathfrak t}$, we derive an exact sequence
\begin{equation*}
  0 \rightarrow  \symmetric^2(TM)\otimes {\mathfrak t} \rightarrow \kernel \omega_2 \rightarrow J^1{\mathfrak t} \rightarrow 0.
\end{equation*}
Since $J^2 {\mathfrak t} $ itself occurs in a natural exact sequence 
\begin{equation*}
  0 \rightarrow \symmetric^2(TM) \otimes {\mathfrak t} \rightarrow J^2 {\mathfrak t} \rightarrow J^1{\mathfrak t} \rightarrow 0,
\end{equation*}
the bundles $\kernel \omega_2$ and $J^2 {\mathfrak t}$ have the same rank.

Recalling that $X \subset J^1{\mathfrak t} $ is holonomic if and only
if ${\mathcal D} X=0$, it is not hard to see that $X$ is holonomic if
and only if $J^1 X \subset \kernel \omega_2$. Lemma \ref{prolongation}
is then a corollary of Proposition \ref{boot}.
%%%%%%%%%%%%%%%%%%%%%% body ends here %%%%%%%%%%%%%%%%%%%%%%%%%
%


\begin{thebibliography}{10}

\bibitem{Blaom_05}
A.~D. Blaom.
\newblock {Geometric structures as deformed infinitesimal symmetries}.
\newblock {\em Trans. Amer. Math. Soc.}, 358:3651--3671, 2006.

\bibitem{Bryant_etal_91}
R.~L. Bryant, S.~S. Chern, R.~B. Gardner, H.~L. Goldschmidt, and P.~A.
  Griffiths.
\newblock {\em {Exterior Differential Systems}}, volume~18 of {\em Mathematical
  Sciences Research Institute Publications}.
\newblock Springer-Verlag, New York, 1991.

\bibitem{CannasdaSilva_Weinstein_99}
A.~Cannas~da Silva and A.~Weinstein.
\newblock {\em {Geometric Models for Noncommutative Algebras}}, volume~10 of
  {\em Berkeley Mathematics Lecture Notes}.
\newblock American Mathematical Society, Providence, RI, 1999.

\bibitem{Cap_Gover_02}
A.~{\v{C}}ap and A.~R. Gover.
\newblock {Tractor calculi for parabolic geometries}.
\newblock {\em Trans. Amer. Math. Soc.}, 354(4):1511--1548, 2002.

\bibitem{Crainic_Fernandes_05}
M.~Crainic and R.~L. Fernandes.
\newblock {Secondary characteristic classes of Lie algebroids}.
\newblock In {\em Quantum field theory and noncommutative geometry}, volume 662
  of {\em Lecture Notes in Phys.}, pages 157--176. Springer, Berlin, 2005.

\bibitem{Crampin_09}
Michael Crampin.
\newblock Cartan connections and {L}ie algebroids.
\newblock {\em SIGMA Symmetry Integrability Geom. Methods Appl.}, 5:Paper 061,
  13, 2009.

\bibitem{Fernandes_02}
R.~L. Fernandes.
\newblock {Lie algebroids, holonomy and characteristic classes}.
\newblock {\em Adv. Math.}, 170(1):119--179, 2002.

\bibitem{Gardner_89}
R.~B. Gardner.
\newblock {\em {The Method of Equivalence and its Applications}}.
\newblock SIAM, Philadelphia, 1989.

\bibitem{Hughen_95}
W.~K. Hughen.
\newblock {\em {The sub-Riemannian geometry of three-manifolds}}.
\newblock PhD thesis, Duke University, 1995.

\bibitem{Ivey_Landsberg_03}
T.~A. Ivey and J.~M. Landsberg.
\newblock {\em {Cartan for Beginners: Differential Geometry via Moving Frames
  and Exterior Differential Systems }}, volume~61 of {\em Graduate Studies in
  Mathematics}.
\newblock American Mathematical Society, 2003.

\bibitem{Kobayashi_72}
S.~Kobayashi.
\newblock {\em {Transformation Groups in Differential Geometry}}.
\newblock Springer-Verlag, New York, 1972.

\bibitem{Mackenzie_05}
K.~C.~H. Mackenzie.
\newblock {\em {General Theory of Lie Groupoids and Lie Algebroids}}, volume
  213 of {\em London Mathematical Society Lecture Note Series}.
\newblock Cambridge University Press, Cambridge, 2005.

\bibitem{Marsden_Ratiu_94}
J.~E. Marsden and T.~S. Ratiu.
\newblock {\em {Introduction to Mechanics and Symmetry}}, volume~17 of {\em
  Texts in Applied Mathematics}.
\newblock Springer, 1994.

\bibitem{Montgomery_02}
R.~Montgomery.
\newblock {\em {A tour of subriemannian geometries, their geodesics and
  applications}}, volume~91 of {\em Mathematical Surveys and Monographs}.
\newblock American Mathematical Society, 2002.

\bibitem{Olver_95}
P.~J. Olver.
\newblock {\em {Equivalence, Invariants, and Symmetry}}.
\newblock Cambridge University Press, Cambridge, 1995.

\bibitem{Sharpe_97}
R.~W. Sharpe.
\newblock {\em {Differential Geometry: Cartan's generalization of Klein's
  Erlangen program}}, volume 166 of {\em Graduate Texts in Mathematics}.
\newblock Springer, 1997.

\bibitem{Sternberg_64}
S.~Sternberg.
\newblock {\em {Lectures on Differential Geometry}}.
\newblock AMS Chelsea Publishing, Providence, 1964.

\bibitem{Stormark_00}
O.~Stormark.
\newblock {\em {Lie's Structural Approach to {PDE} Systems}}, volume~80 of {\em
  Encyclopedia of Mathematics and its Applications}.
\newblock Cambridge University Press, Cambridge, 2000.

\bibitem{Tang_04}
X.~Tang.
\newblock {Deformation Quantization of Pseudo Symplectic (Poisson) Groupoids}.
\newblock {\em {\normalfont Preprint}}, 2004.

\end{thebibliography}
\end{document}